\documentclass[12pt,a4paper]{amsart}

\usepackage[american]{babel}
\usepackage{amsmath}
\usepackage{amssymb}
\usepackage{amsbsy}
\usepackage{amsfonts}
\usepackage{graphicx}
\usepackage{subfigure}
\usepackage[normalem]{ulem}
\usepackage{siunitx}
\usepackage{chemfig}
\usepackage{mhchem}
\usepackage{import}

\usepackage{stmaryrd}

\graphicspath{{./figures/}}
\DeclareGraphicsExtensions{.pdf,.jpg,.png,.eps}

\usepackage{url}

\usepackage[format=plain, labelsep=endash, justification=justified, singlelinecheck=true, font={footnotesize,up}]{caption}

\newcommand{\Frac}{\displaystyle\frac}

\newtheorem{remark}{Remark}

\definecolor{MidnightBlue}{rgb}	{0,0.1,.5}
\definecolor{darkgreen}{rgb}	{0,0.6,0}

\begin{document}

\title[Hierarchical Simulation of Bio-Hybrid Interfaces]
{Hierarchical Electrochemical Modeling and Simulation 
of Bio-Hybrid Interfaces}

\author{Emanuela Abbate$^{1}$ \and Matteo Porro$^{2,4}$ \and 
Thierry Nieus$^{3}$ \and Riccardo Sacco$^{4}$}

\address{$^{1}$ Inria Bordeaux Sud-Ouest, \\
200 avenue de la vieille tour, 
33405 Talence Cedex, France\\
		{\tt e-mail:} emanuela.abbate@inria.fr \\
$^{2}$ Center for Nano Science and Technology @PoliMi,
           Istituto Italiano di Tecnologia, \\
					 via Pascoli 70, 20133 Milano, Italy \\	
		{\tt e-mail:} matteo.porro@iit.it            \\
$^{3}$ Nets3 laboratory,
		       Department of Neuroscience and Brain Technologies, \\ 
					 Istituto Italiano di Tecnologia, 
					 via Morego 30, 16163 Genova, Italy \\
		{\tt e-mail:} thierry.nieus@iit.it            \\
$^{4} $Dipartimento di Matematica,
               Politecnico di Milano, \\
	           Piazza Leonardo da Vinci 32, 20133 Milano, Italy \\
{\tt e-mail:} riccardo.sacco@polimi.it}

\date{\today}

\begin{abstract}
In this article we propose and investigate a hierarchy of 
mathematical models based on partial differential equations 
(PDE) and ordinary differential equations (ODE) 
for the simulation of the biophysical phenomena occurring 
in the electrolyte fluid that connects a biological component
(a single cell or a system of cells) and a solid-state device 
(a single silicon transistor or an array of transistors). 
The three members of the hierarchy, ordered by decreasing complexity,
are: ({\it i}) a 3D Poisson-Nernst-Planck (PNP) PDE system for ion
concentrations and electric potential; 
({\it ii}) a 2D reduced PNP system for the same dependent variables 
as in ({\it i}); ({\it iii}) 
a 2D area-contact PDE system for electric potential
coupled with a system of ODEs for ion concentrations.
The backward Euler method is adopted for temporal 
semi-discretization and a fixed-point iteration 
based on Gummel's map is used to decouple system equations. 
Spatial discretization is performed using piecewise linear 
triangular finite elements stabilized via edge-based exponential fitting. 
Extensively conducted simulation results are in excellent agreement 
with existing analytical solutions of the PNP problem in radial coordinates 
and experimental and simulated data using simplified lumped parameter models.
\end{abstract}

\maketitle


\section{Introduction and Motivation}\label{sec:intro}

\begin{figure}[th!]
\centering
\subfigure[Rat neuron on electronic substrate]
{\includegraphics[height=0.32\textwidth]{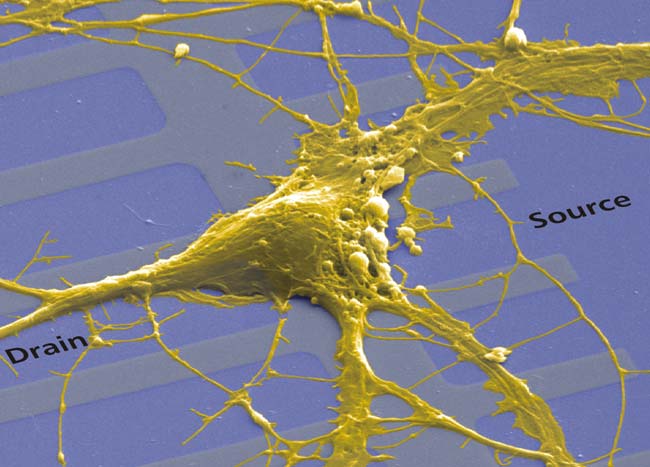}
\label{fig:EOS}}
\hfill
\subfigure[Neuro-chip]
{\includegraphics[height=0.32\textwidth]{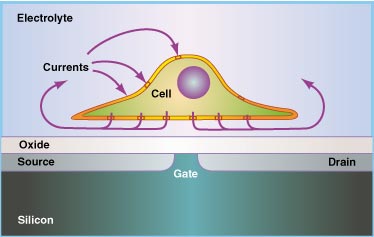}
\label{fig:hybrid_dev}}
\caption{Left: rat neuron grown on an EOSFET, 
image reprinted from~\cite{voelkerfromherz2005small}. 
Right: schematics of a neuro-chip, image 
reprinted from~\cite{neher2001molecular}.}
\label{fig:neurochips}
\end{figure}

In this article we address the study of a class of problems arising
in the context of Bioelectronics, a recently emerged 
discipline at the crossroad among Nanotechnology, Solid-State Electronics, 
Biology and Neuroscience.
The focus of our investigation is on the mathematical and computational
modeling of bioelectronic interfaces 
(see~\cite{willner2006bioelectronics,Fro2012} for a review 
and~\cite{hutzler2006high,Ghezzi2011,Nieus2012,moulin2008new} for 
a selection of significant applications).
Bioelectronic interfaces are bio-hybrid structures constituted by 
living cells attached to an electronic substrate and surrounded 
by an electrolyte bath.
An example can be seen in Fig.~\ref{fig:EOS} which shows
an electronmicrograph of hippocampal neuron 
cultured on an electrolyte-oxide-silicon field-effect transistor 
(EOSFET)~\cite{voelkerfromherz2005small}.
Fig.~\ref{fig:hybrid_dev} reports a schematic cross-section
view of a neuro-chip, which allows to identify 
the main parts of the bio-hybrid system:
the cell, the extracellular bath,
the thin interstitial cleft separating the cell and the
electronic substrate, the protective oxide layer deposited on the top 
of the substrate, 
and the source-to-drain transistor structure. 

In the basic function mode of the EOSFET, 
as a consequence of the cellular activity
elicited by the application of an external stimulus, ionic current 
flows through the adhering cell membrane and along 
the cleft. The resulting extracellular voltage 
turns out to play the role of the gate voltage 
which controls electrical charge flow in the substrate
and, ultimately, the current flowing out from the 
drain terminal of the device.

The interface contact in the scheme illustrated in 
Fig.~\ref{fig:neurochips} is realized by the thin conductive electrolyte
separating the two subsystems, whose amplitude is smaller
than the cell radius by about three orders of magnitude. 
Therefore, the cell-chip junction forms a planar electrical 
core-coat conductor and the main physical phenomena 
(ion electrodiffusion and gate voltage modulation) take 
place in this three-dimensional region whose vertical thickness 
is much smaller than the two-dimensional area where the cell 
adheres to the substrate. 

The above description indicates that a sound mathematical
picture of a bioelectronic interface requires the adoption 
of a genuine multiscale perspective. For this reason, 
in this article 
we continue our analysis started off in~\cite{brera2009multiphysics}
and numerically investigated in~\cite{brera2010conservative}
and propose a hierarchy of models based on PDEs and ODEs 
for the simulation of the biophysical phenomena occurring 
in the 3D interface contact described above.
The hierarchy includes the following three members, ordered by decreasing level of complexity:
({\it i}) a 3D Poisson-Nernst-Planck (PNP) PDE system for ion
electrodiffusion and electric potential dynamics~\cite{rubinstein1990electro};
({\it ii}) a 2D reduced PNP system for the same dependent variables 
and phenomena as in ({\it i}); ({\it iii}) 
a 2D area-contact PDE system for electric potential dynamics
coupled with a system of ODEs for ion dynamics. 
This last member of the hierarchy is a variant of the 
area-contact model proposed and studied in~\cite{Brittinger2005}.

Model ({\it i}) is the most accurate in the hierarchy but, of course,
requires a considerable amount of computational effort for its numerical
simulation. Model ({\it ii}) is obtained 
by averaging the 3D PNP equations in the direction 
$z$ perpendicular to the electrolyte cleft. This model reduction 
procedure leads to a modified PNP system 
to be solved in a 2D plane $x$-$y$ parallel to the substrate.
Model ({\it iii}) is a further reduction of ({\it ii}) obtained
by neglecting spatial dependence of ion concentrations in the electrolyte
cleft. This leads to a time-dependent 2D Poisson equation for 
electric potential coupled with electrolyte cleft ion dynamics 
described by a system of ODEs as done in~\cite{Brittinger2005}.
In all members of the hierarchy, iono-electric coupling between substrate
and electrolyte is accounted for by ``lumped" 
transmission conditions expressing continuity
of dielectric and ionic fluxes across the interfaces. 
Electrodiffusive ionic coupling between cell(s) and 
electrolyte is described through a variety of 
transmembrane currents including the Goldman-Hodgkin-Katz 
and Hodgkin-Huxley models~\cite{mori2006three,moulin2008new}.

The backward Euler method is adopted for temporal 
semi-discretization and a fixed-point iteration 
based on Gummel's map~\cite{jerome1996analysis} 
is used to decouple system equations. 
Spatial discretization is performed using a generalization
to axisymmetric cylindrical coordinates of the
piecewise linear triangular finite element scheme 
stabilized via edge-based exponential fitting 
proposed and analyzed in~\cite{xu1999monotone}. 
 
Extensively conducted simulations 
using the full 3D PNP model reveal that 
ion concentration and electric field variations
mainly occur in the electrolyte cleft.
Sensible results are also obtained addressing non ideal effects occurring
in realistic devices, such as undesired cellular activity 
detection on more than one electrode and cellular cross-stimulation.
Simulation experiments using the 2D formulations
proposed in the present article
demonstrate their excellent agreement with  
experimental and numerical results in the existing literature 
\cite{Pabst2007,Brittinger2005} but with a significant saving of 
computational effort with respect to 
the solution of the full 3D PNP system.

A short outline of the article is as follows.
In Sects.~\ref{sec:models} and~\ref{sec:hierarchy} we illustrate the hierarchy 
of mathematical models used to represent ion electrodiffusion 
and electric potential distribution in the bioelectronic structure.
In Sect.~\ref{sec:methods} we describe the 
numerical techniques adopted for the solution of the discrete problem.
In Sect.~\ref{sec:results} we address the 
validation of the proposed computational model in the
simulation of several test cases of 
biophysical significance. In Sect.~\ref{sec:conclusions} we summarize the
main contents of our analysis and indicate some perspectives for
future research directions.

\section{Three-Dimensional Model}\label{sec:models}

In this section we illustrate a three-dimensional model 
of ion electrodiffusion throughout the interstitial cleft 
separating the cell and the electronic substrate under the
application of an external stimulus.

\subsection{Geometrical model}
\label{sec:3Dgeometrical_model}
\begin{figure}[h!]
\centering
\subfigure[Cell, electrolyte and substrate]
{ \includegraphics[width=0.45\textwidth]{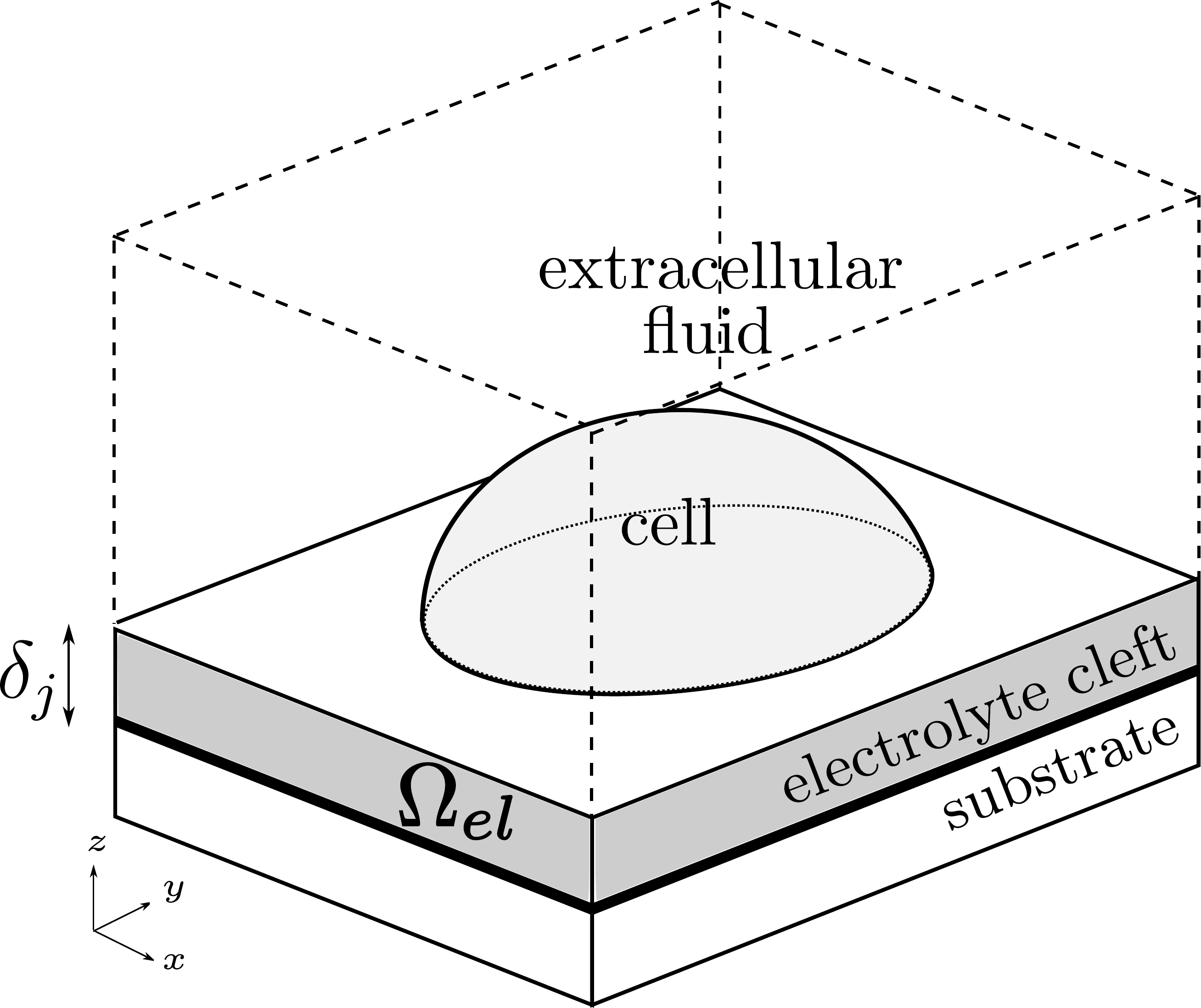}
\label{fig:Physical-model}}
\hfill
\subfigure[3D computational domain and interfaces]
{\includegraphics[width=0.5\textwidth]{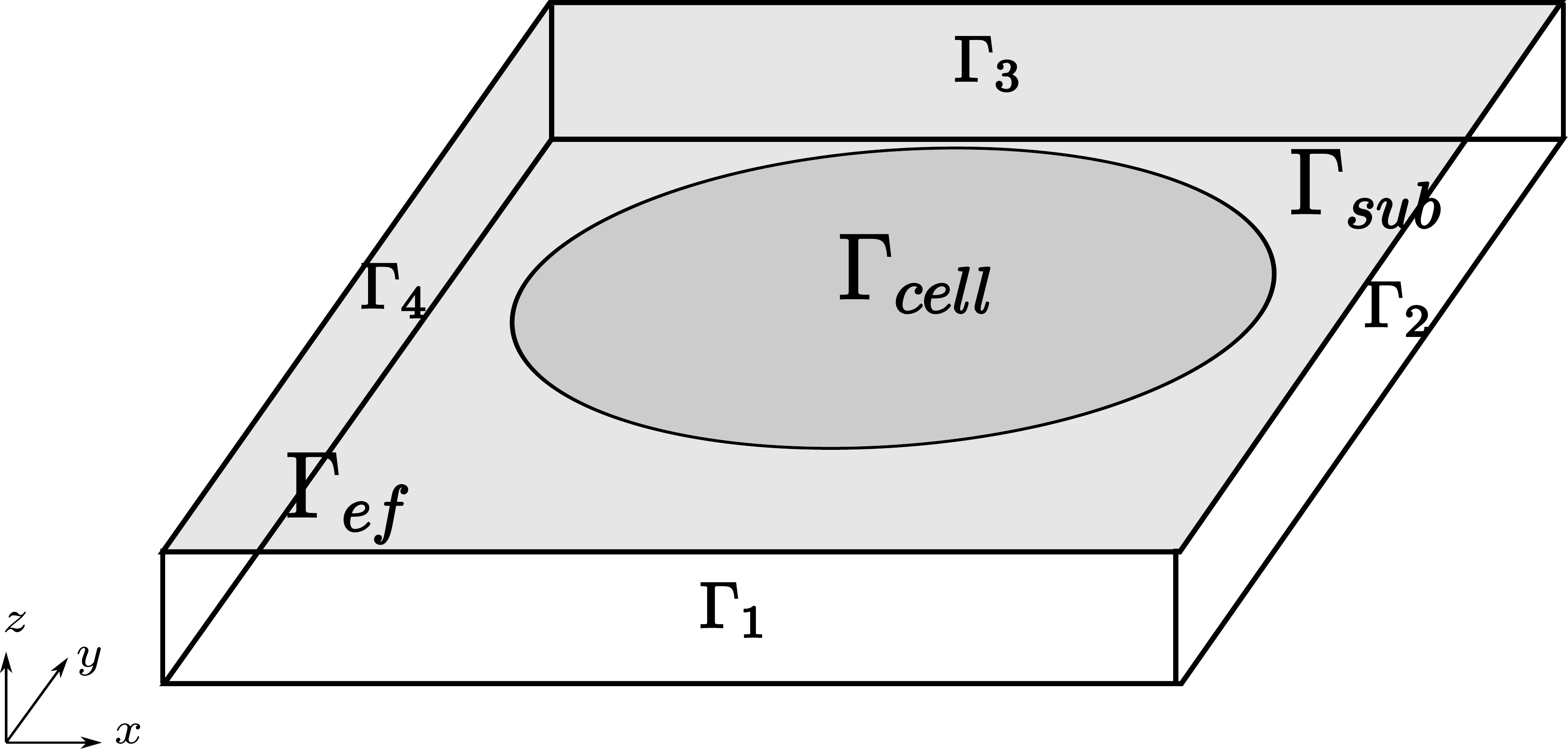}
\label{fig:omega3D}}
\caption{Left: geometrical model of the whole bio-hybrid system. 
A cell surrounded by an electrolyte bath is attached to an electronic device.
Right: the computational domain $\Omega_{el}$ is the 
thin layer of electrolyte between the cell and the substrate. 
Seven different boundary regions are distinguished: the upper surface
is divided in $\Gamma_{cell}$ (the cell attachment area) and $\Gamma_{ef}$
(surface separating the electrolyte cleft from the extracellular fluid).}
\label{fig:3Ddomains}
\end{figure}

Fig.~\ref{fig:Physical-model} shows a 3D schematic picture 
of a cell-to-substrate interface.
The cell shape is represented as a 
rotational solid separated by the planar substrate 
by a thin electrolyte domain $\Omega_{el}$.
Despite the extracellular fluid is surrounding all the cell, 
most of the effects resulting from cell stimulation occur 
in the adhesion region as is demonstrated in the numerical
experiments reported in Section~\ref{sec:voltage_clamp_validation}. 
For this reason, we restrict the geometrical space where to solve 
the mathematical problem to the 3D region $\,\Omega_{el}$, see Fig.~\ref{fig:omega3D}, 
which includes the cleft between the attached cell membrane
and the device, but also the part of electrolyte in the neighborhood of the cell
close to the substrate.
The layer thickness $\delta_{j}$ is of the order of 50\,$\div$\,\SI{100}{\nano\meter}
while the cell radius is around \SI{10}{\micro\meter} in the considered 
applications~\cite{Brittinger2005,Pabst2007}.

\subsection{The Poisson-Nernst-Planck system}
\label{sub:The-Poisson-Nernst}
The Poisson-Nernst-Planck system (PNP) for ion electrodiffusion 
reads~\cite{rubinstein1990electro}:
\begin{subequations}\label{eq:PNP_flux}
\begin{align}
\dfrac{\partial c_{i}}{\partial t}+\text{div}\,\mathbf{f}_{i}
\left(c_{i},\varphi\right) & = 0
&& i=1,\ldots,M
\label{eq:continuity-1-1}
\\
\mathbf{f}_{i}\left(c_{i},\varphi\right) & 
= -D_{i}\nabla c_{i} + \mu_{i} \Frac{z_{i}}{|z_i|} c_{i}\mathbf{E}
&&
i = 1,\ldots,M
\label{eq:nernstplanck-1}
\\
D_{i} & 
= \dfrac{\mu_{i}V_{th}}{\left|z_{i}\right|}
&&
i=1,\ldots,M
\label{eq:einstein_rel-1}
\\
\text{div}\,\mathbf{E} & 
= \dfrac{1}{\varepsilon} \rho 
\label{eq:poisson-1-1}
\\
\mathbf{E} &
=  -\nabla\varphi
\label{eq:cel_2}
\\
\rho &
= q \sum_{i=1}^M z_{i}c_{i}.
\label{eq:rho}
\end{align}
Eq.~\eqref{eq:continuity-1-1} is the continuity equation
describing mass conservation for each ion whose concentration
is denoted by $c_i$ (\si{\metre^{-3}}), $i=1, \ldots, M$, $M \geq 1$ being the number
of ions flowing in the electrolyte fluid.
Each ion flux density $\mathbf{f}_{i}$ (\si{\metre^{-2} s^{-1}}) 
is defined by the Nernst-Planck relation~\eqref{eq:nernstplanck-1}
in which 
it is possible to recognize a chemical contribution and an electric contribution,
in such a way that the model be regarded as an extension of Fick's law of diffusion 
to the case where the diffusing particles are also moved by electrostatic forces with
respect to the fluid. The quantity $z_i$ is the valence of the $i$-th ion, while 
$\mu_{i}$ and $D_{i}$ are the mobility and diffusivity of the chemical species, respectively,
related by the Einstein relation \eqref{eq:einstein_rel-1} 
where $V_{th}=k_{B}T/q$ (\si{V}) is the thermal potential
($k_{B}$ is the Boltzmann constant, $T$ is the absolute temperature
and $q$ is the elementary charge).
The electric field $\mathbf{E}$ (\si{V m^{-1}}) 
due to space charge distribution $\rho$ in the electrolyte is determined 
by the Poisson equation~\eqref{eq:poisson-1-1} which represents Gauss' 
law in differential form, $\varepsilon$ being uthe dielectric permittivity of 
the fluid medium.
For further analysis, it is useful  to introduce 
the electrical current density 
$\mathbf{j}_{i}$ (\si{A \metre^{-2}}), equal to the 
number of ion charges flowing through a given surface area per unit time
and defined as
\begin{equation}
\mathbf{j}_{i}:=q z_i \mathbf{f}_{i} = 
- q z_i D_i \nabla c_i + q \mu_i |z_i| c_i \mathbf{E} \qquad 
i=1, \ldots, M.
\label{eq:current_density}
\end{equation}

\begin{remark}
The PNP system~\eqref{eq:PNP_flux} has the same format and structure
as the Drift-Diffusion equations for 
semiconductors (see, e.g.,~\cite{jerome1996analysis}),
but it is applied to a different medium (water instead of a semiconductor
crystal lattice) and includes, in general, 
more charge carriers than just holes and electrons, as in the case
of semiconductor device theory.
\end{remark}
\end{subequations}

\subsection{Boundary and initial conditions}
\label{sub:Boundary-and-initial}
Let $t$ and $\mathbf{x}$ denote the time variable and 
the spatial coordinate,
respectively. We denote also by $\Gamma:= \partial \Omega_{el}$
the boundary of the computational domain in Fig.~\ref{fig:omega3D} 
and by $\mathbf{n}$ the unit outward normal vector on $\Gamma$.

The initial conditions $c_{i}^{0}(\mathbf{x})=c_{i}\left(0,\mathbf{x}\right)$
and $\varphi^{0}(\mathbf{x})=\varphi(0,\mathbf{x})$ are determined
by solving the static version of the PNP system~\eqref{eq:PNP_flux}
in the domain $\Omega_{el}$, which corresponds to setting 
$\frac{\partial c_{i}}{\partial t}=0$
in~\eqref{eq:continuity-1-1} for each ion $i=1,\ldots, M$.

\begin{subequations}\label{eq:bcs}
The boundary conditions deserve a deeper discussion because they need
to mathematically express the coupling between the electrolyte cleft 
and the surrounding environment, comprising the extracellular fluid and the
active parts of the bio-hybrid system, namely, the cell, the
membrane and the electronic substrate (see Fig.~\ref{fig:Physical-model}).
Referring to Fig.~\ref{fig:omega3D} for the notation, 
we distinguish among seven different regions in the boundary 
$\Gamma$: four lateral sides
$\Gamma_i$, $i=1, \ldots, 4$, a lower face $\Gamma_{sub}$ in contact
with the substrate, and an upper face
divided into two parts, the area attached to the cell
$\Gamma_{cell}$, and the free region $\Gamma_{ef}$, covered 
on top by the surrounding volume of extracellular fluid. 
Accordingly, the following conditions
are enforced on the electric field and the particle fluxes:
\begin{align}
\varphi 
& 
= V_{bath}
&& \text{on} \, \Gamma_{1} \cup \Gamma_{2}
\cup \Gamma_{3} \cup \Gamma_{4}
\label{eq:gamma1_phi}
\\
\llbracket\mathbf{D\cdot n}\rrbracket_{\Gamma_{ef}} 
& 
= 0
&&
\text{on} \, \Gamma_{ef}
\label{eq:Edotn_electr}
\\
\llbracket\mathbf{D\cdot n}\rrbracket_{\Gamma_{cell}} 
& 
= 0 
&&
\text{on} \, \Gamma_{cell}
\label{eq:jumpDcell}
\\
\llbracket\mathbf{D\cdot n}\rrbracket_{\Gamma_{sub}} 
& 
= 0
&&
\text{on} \, \Gamma_{sub}
\label{eq:jumpDsub}
\\
c_{i} & 
= c_{i}^{bath}
&&
\text{on} \, 
\Gamma_{1}\cup\Gamma_{2}\cup\Gamma_{3}\cup\Gamma_{4}\label{eq:gamma1_ci}
\\
\llbracket\mathbf{f}_{i}\cdot\mathbf{n}\rrbracket_{\Gamma_{ef}} & 
= 0
&&
\text{on} \, \Gamma_{ef}
\label{eq:Jdotn_electr}
\\
\llbracket\mathbf{f}_{i}\cdot\mathbf{n}\rrbracket_{\Gamma_{cell}} & 
= 0
&&
\text{on} \, \Gamma_{cell}
\label{eq:Jjump_cell}
\\
\mathbf{f}_{i}\cdot\mathbf{n} & 
= 0
&&
\text{on} \, \Gamma_{sub}
\label{eq:JdotnSub}
\end{align}
having denoted by
$\llbracket\cdot\rrbracket_{\zeta}$ the jump operator
restricted to the interface $\zeta$.

Eqns.~\eqref{eq:gamma1_phi} and~\eqref{eq:gamma1_ci} are Dirichlet
boundary conditions that can be interpreted as 
\textquotedblleft{}far field conditions\textquotedblright{}, meaning that sufficiently far from the surface where the cell is attached
to the substrate, we can assume the electric potential $\varphi$ to be fixed at a
given constant reference value and 
each ion concentration to be fixed at a given 
constant value $c_{i}^{bath}$.
The quantities $c_{i}^{bath}$, $i=1, \ldots, M$, are 
physiologically given values in such a way that 
the electrolyte is a neutral solution, i.e.,
\begin{equation}
\rho = 
\rho_{bath}=q\sum_{i=1}^{M}z_{i}c_{i}^{bath}=0.\label{eq:electroneutr_bath}
\end{equation}
\end{subequations}

\subsubsection*{Cleft-cell coupling}
\label{sub:Cleft-cell-coupling}
The most relevant physiological phenomena occurring in the
bio-hybrid interface depend on the properties of the 
cell membrane and can be modeled as the sum of two contributions,
one from the lipid portion of the membrane and one from the 
ionic channels.
\begin{figure}[t!]
\centering
\subfigure[Membrane]
{\includegraphics[height=0.3\textwidth]{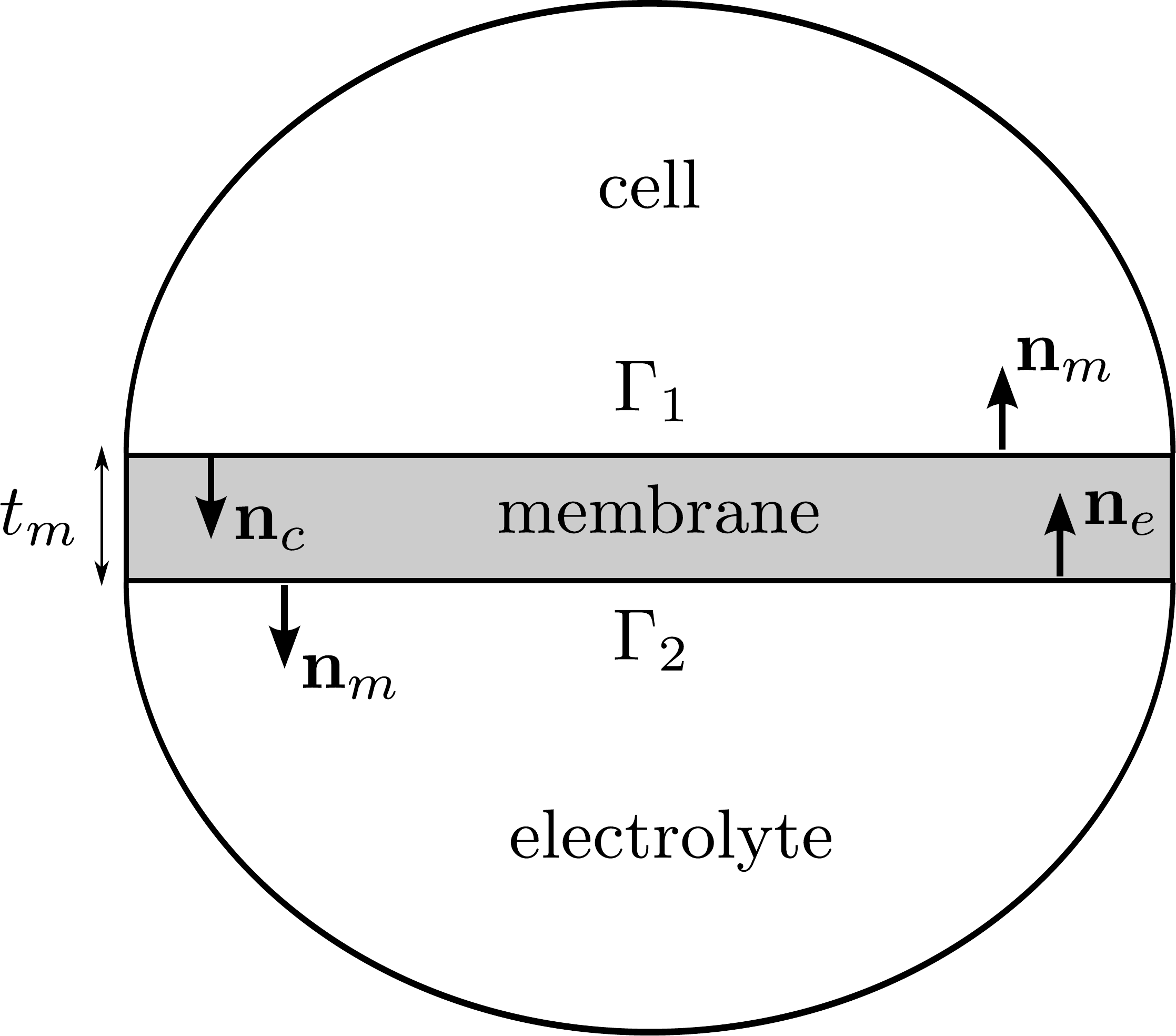}
\label{fig:membrane}}
\qquad
\subfigure[Lumping in $\Gamma_{cell}$]
{\includegraphics[height=0.3\textwidth]{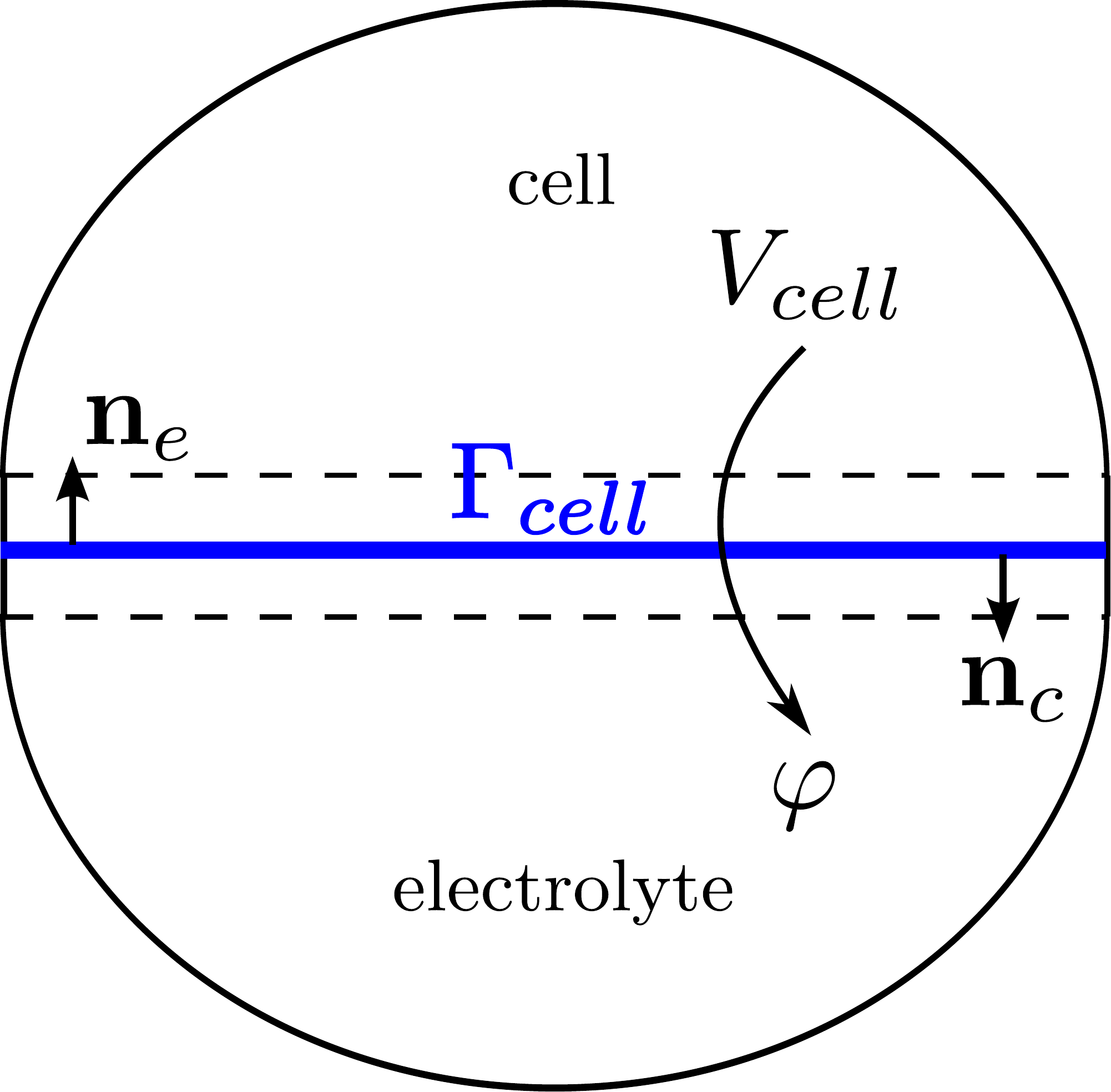}
\label{fig:lumpingmembrane}}
\caption{On the left: cell
and electrolyte separated by the membrane with its physical thickness
t$_{M}$. On the right: cell and electrolyte separated by an
interface $\Gamma_{cell}$ with zero thickness, result of the lumping
of the original boundaries $\Gamma_{1}$ and $\Gamma_{2}$ of the
membrane region.}
\label{fig:cleft_cell_coupling}
\end{figure}

The membrane subdomain (shown in Fig.~\ref{fig:membrane})
has a thickness t$_{M}$ of the order of 5\,$\div$\,\SI{10}{\nano\meter}, 
which is much smaller than the characteristic size of the domain 
$\Omega_{el}$. Therefore, a geometrical discretization of this 
small region may give rise to a huge number of degrees of freedom
for the numerical method. To reduce computational complexity, 
we apply the membrane model proposed in~\cite{mori2006three,brera2010conservative}.
This amounts to assuming that $\varphi$ varies linearly across the
membrane thickness so that, upon  introducing 
the two dimensional manifold $\Gamma_{cell}$ corresponding
to the middle cross-section of the membrane volume, 
the transmission condition~\eqref{eq:jumpDcell}
across the two dimensional manifold $\Gamma_{cell}$ becomes
\begin{subequations}\label{eq:membrane_models}
\begin{equation}
\mathbf{D}_{c}\cdot\mathbf{n}_{c}=-\mathbf{D}_{e}
\cdot\mathbf{n}_{e}=-\varepsilon_{M}\dfrac{\varphi_{m2}-\varphi_{m1}}
{\text{t}_{M}}
\simeq-C_{M}\left(\varphi_{m2}-V_{cell}\right),\label{eq:CM_condition}
\end{equation}
where $\varepsilon_{M}$ is the membrane permittivity, 
$V_{cell}$ is the intracellular potential while
$\varphi_{m1}$ and $\varphi_{m2}$ are the traces of $\varphi$
at both sides of $\Gamma_{cell}$ (cell and electrolyte, respectively).
Condition~\eqref{eq:CM_condition} expresses a capacitive coupling
between cell and cleft through the 
membrane specific capacitance $C_{M}:= \varepsilon_{M}/\text{t}_{M}$ (\si{\farad\per\square\metre}). 

To account for membrane channels, the ionic current in the case 
of active cells
is described by
the following 
generalized Hodgkin-Huxley (HH) model (for a detailed description 
see~\cite{hodgkin1952currents,hodgkin1952quantitative,hille2001ion,keener2009mathematical})
\begin{equation}
j_{i}^{tm}=j_{i}^{tm}\left(t,\mathbf{x},\mathbf{s},V_{cell},
\varphi,\mathbf{c}^{cell},\mathbf{c}\right)\label{eq:hodgkin_hux}
\end{equation}
where $\mathbf{s}$ is a vector collecting the gating variables $n, m$ and
$h$ responsible of channel opening probabilities, while
$\mathbf{c}^{cell}$ and $\mathbf{c}$ are arrays of size $M$
containing all the ion concentrations inside and outside the cell.
If instead only passive cells are considered, as in 
most simulations performed in this work, 
the adopted model is the Goldman-Hodgkin-Katz
(GHK) equation for the current density~\cite{hille2001ion}
\begin{equation}
j_{i}^{tm}=p_{i}z_{i}q
\left[Be\left(-\dfrac{z_{i}\left(V_{cell}-\varphi\right)}
{V_{th}}\right)c_{i}^{cell}-
Be\left(\dfrac{z_{i}\left(V_{cell}-\varphi\right)}
{V_{th}}\right)c_{i}\right],\label{eq:ghk_model}
\end{equation}
where $p_{i}$ is the permeability constant of the specific ion 
(\si{\metre\per\second}) and $Be(x):=x/(e^{x}-1)$ is the inverse of the Bernoulli function.
With these descriptions of the transmembrane current density $j_{i}^{tm}$
for each ion species,
condition~\eqref{eq:Jjump_cell} becomes
\begin{equation}
\mathbf{f}_{e}^{i}\cdot\mathbf{n}_{e}=-\mathbf{f}_{c}^{i}\cdot\mathbf{n}_{c}=-\dfrac{j_{i}^{tm}}{qz_{i}}.\label{eq:flux_through_membr}
\end{equation}
\end{subequations}

\subsubsection*{Cleft-substrate coupling}
In the present work, the action
of the electro-chemical bounding of ions at the interface between 
electrolyte cleft and substrate is neglected, and
the semiconductor device is assumed to behave
as a MOS capacitor mathematically described through a lumped equivalent
model as in~\eqref{eq:CM_condition}. 
Referring to Fig.~\ref{fig:cleft_bulk_coupling}
the capacitive coupling on $\Gamma_{sub}$ is
\begin{subequations}\label{eq:cleft-sub-coupling}
\begin{equation}
\mathbf{D}_{s}\cdot\mathbf{n}_{s}=-\mathbf{D}_{e}\cdot\mathbf{n}_{e}
\simeq-\epsilon_{s}\dfrac{\varphi_{s2}-\varphi_{s1}}{\text{t}_{S}}=-C_{S}\left(\varphi_{s2}-V_{G}\right),\label{eq:CS_condition}
\end{equation} 
where $\varphi_{s1}$ and $\varphi_{s2}$ are the traces of $\varphi$
on both sides of $\,\Gamma_{sub}$. The function $V_{G}=V_{G}\left(t\right)$ 
denotes the value of
the potential on the gate contact, taken to be spatially constant
according to the hypothesis of ideal metallic behavior of the gate.
Regarding the particle fluxes, condition~\eqref{eq:JdotnSub} 
already states that there is no current injection in the electronic device.
\begin{figure}[h!]
\centering 
\includegraphics[width=0.5\textwidth]{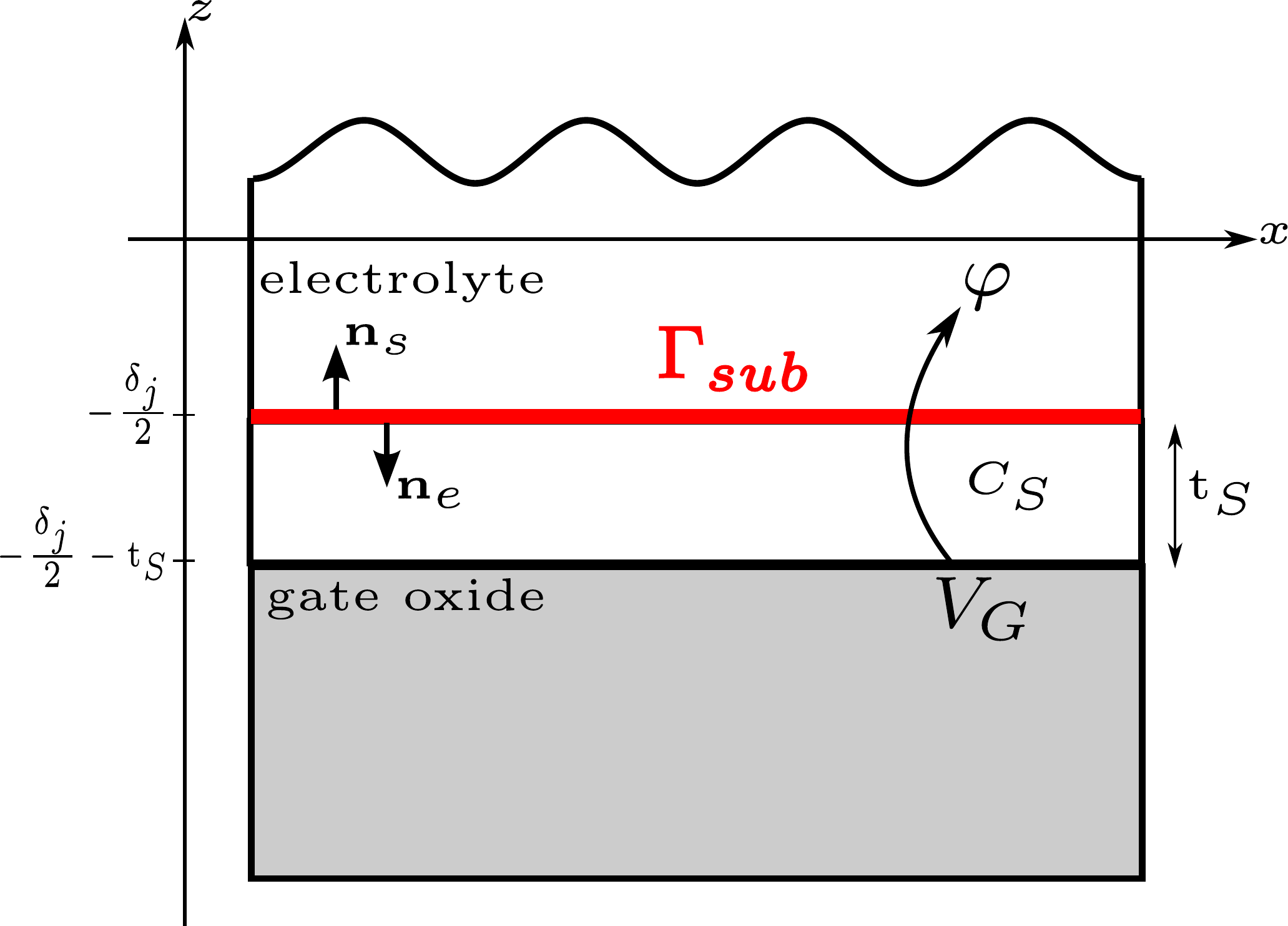}
\caption{Electronic substrate model and coupling with the electrolyte.}
\label{fig:cleft_bulk_coupling}
\end{figure}

\subsubsection*{Electrolyte-electrolyte artificial coupling}
As shown in Fig.~\ref{fig:Physical-model}, 
the electrolyte domain is restricted to a thin sheet 
of amplitude $\delta_{j}$. This approximation leads to 
the boundary conditions
\eqref{eq:Edotn_electr} and~\eqref{eq:Jdotn_electr} on the fictitious boundary
$\Gamma_{ef}$. 
Assuming again that the potential is a linear function of $z$, 
we can rewrite~\eqref{eq:Edotn_electr} as
\begin{equation}
\mathbf{D}_{ext}\cdot\mathbf{n}_{ext}=-\mathbf{D}_{int}\cdot\mathbf{n}_{int}
\simeq C^{*}\left(\varphi_{int}-\varphi_{ext}\right)\simeq C^{*}\left(\varphi_{int}-V_{bath}\right), \label{eq:cstar_first}
\end{equation}
where $\varphi_{int}$ and $\varphi_{ext}$ are the traces of $\varphi$
respectively on the two sides of $\,\Gamma_{ef}$ (see Fig.~\ref{fig:electroelectro_coupling}).
The reduced model~\eqref{eq:cstar_first} consists of assuming
 that far away from the boundary $\,\Gamma_{ef}$ the potential is at the reference value $V_{bath}\,$, $C^{*}$ being 
 a fictitious capacitance introduced to relate the value of the potential in the electrolyte, inside and outside of the computational domain $\Omega_{el}$. 
 A possible modeling approach to estimate the value of $C^*$ 
consists in taking a fraction $1/\kappa$ of the value of $C_{M}$.
Computational experiments indicate that $\kappa = 5$ 
is an appropriate choice.
\begin{figure}[h!]
\centering 
\includegraphics[width=0.55\textwidth]{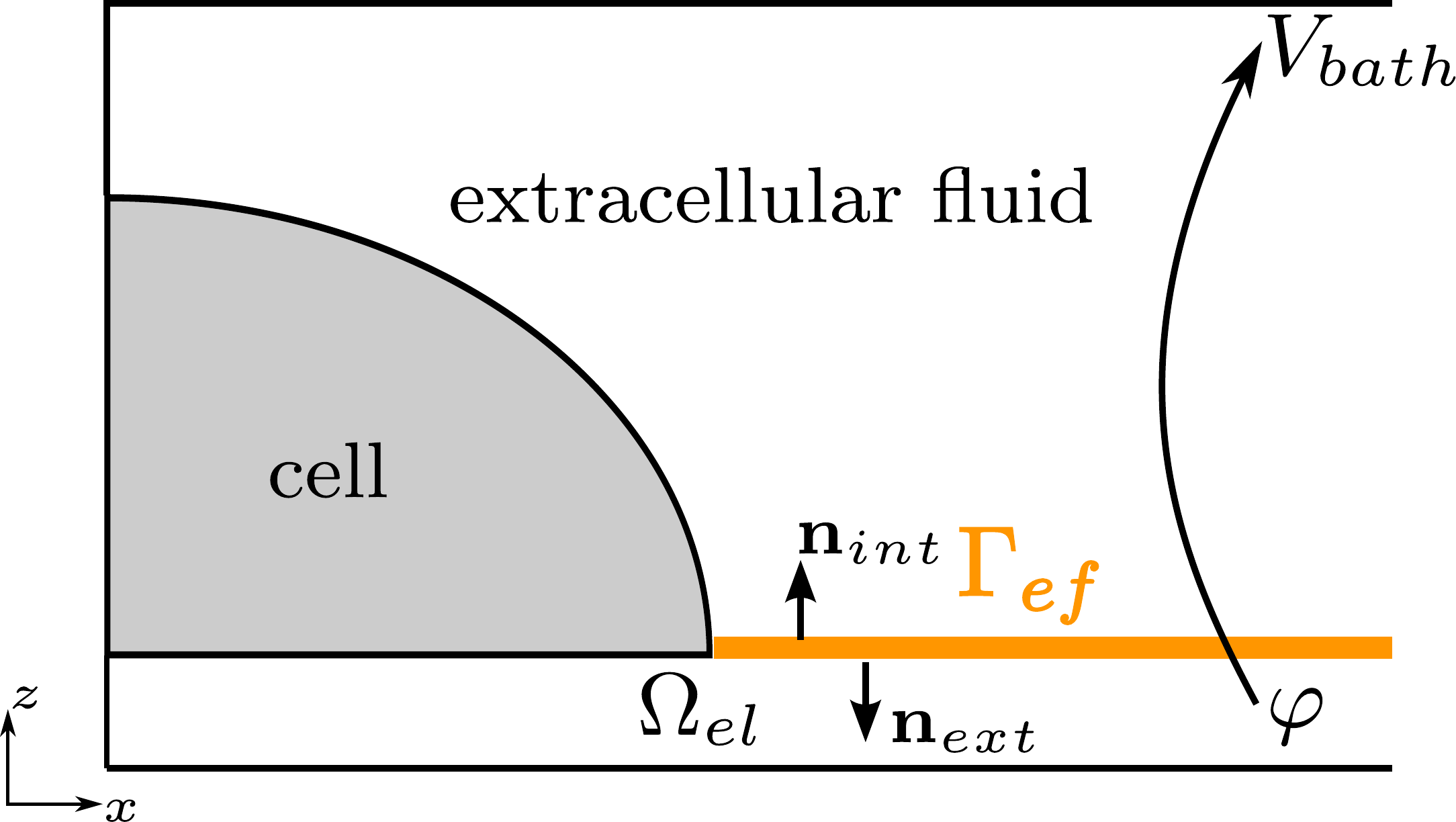}
\caption{Cross section
in the $x$-$z$ plane of the electrolyte bath, illustrating the
coupling condition between $\Omega_{el}$ and the external remaining
electrolyte enforced on $\Gamma_{ef}$.}
\label{fig:electroelectro_coupling} 
\end{figure}

Regarding particle fluxes, we assume that, far away from $\Gamma_{ef}$,
ion concentrations can be considered to be equal to their bath value 
$c_{i}^{bath}$ and the electrolyte to be electroneutral. 
Then condition~\eqref{eq:Jdotn_electr} can be rewritten as
\begin{equation}
\mathbf{f}_{ext}^{i}\cdot\mathbf{n}_{ext}=-\mathbf{f}_{int}^{i}\cdot\mathbf{n}_{int}=-v_{i}^{*}\left(c_{i}^{int}-c_{i}^{ext}\right)\simeq v_{i}^{*}\left(c_{i}^{int}-c_{i}^{bath}\right) 
\quad i=1, \ldots, M.\label{eq:flux_through_el}
\end{equation}
The above relation is a Robin condition for the particle flux density
which physically expresses the fact that ions are allowed to 
cross the fictitious interface $\Gamma_{ef}$, as it should be in 
the non-truncated electrolyte domain.
Mathematically, the quantity $v_{i}^{*}$ is 
an effective permeability (\si{\metre\per\second}) 
whose value can be estimated by equating flux~\eqref{eq:flux_through_el}
to a fraction $1/\kappa^*$ of the flux through the membrane~\eqref{eq:flux_through_membr}.
Computational experiments indicate that $\kappa^* = 20$ 
is an appropriate choice.
\end{subequations}

\section{Hierarchical models}\label{sec:hierarchy}
The analysis of~\cite{Brittinger2005,Pabst2007} 
shows that
the ion current density entering the cleft through the membrane 
 mainly flows parallel to
the $z$-axis. Once inside the cleft, the direction of the
current density changes into the radial one. The time needed
by the ions to flow across the cleft thickness is of the order of \SI{e-7}{\second}. Particles move many times up and down along the 
$z$-direction because the ratio between the cleft thickness 
and the radius of the attached area is of the order of $10^{-3}$
and the resulting contribution of this random motion to the vertical 
current is equal to zero. Thus,  
it appears to be reasonable to 
derive a family of two-dimensional models in the $x-y$ plane from 
the 3D PNP system of Sect.~\ref{sec:models}. 
This is the object of the discussion below.

\subsection{Model reduction: from fully 3D to 2.5D ion electrodiffusion \label{sub:Model-reduction}}

\begin{figure}[h!]
\centering 
\includegraphics[width=0.7\textwidth]{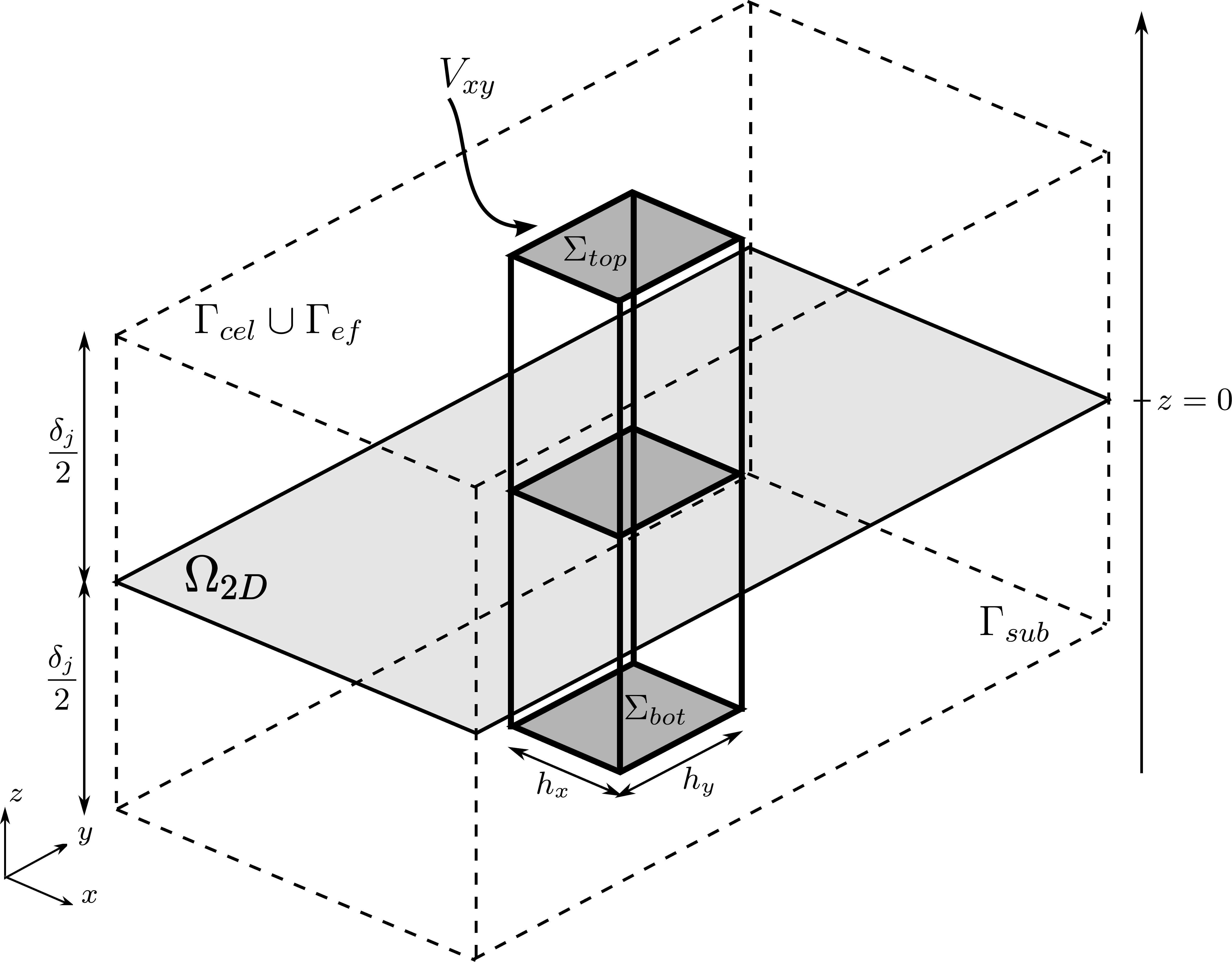}
\caption[Geometrical reduction in the $x$-$y$ plane]
{\label{fig:middleplan_2d}Schematics for the geometrical
reduction in the $x$-$y$ plane: the middle plan of the cleft
becomes the two-dimensional domain $\Omega_{2D}$. $V_{xy}$ is a
control volume used to compute the integrals and the fluxes.}
\end{figure}

We place a coordinate system with the origin
in the middle of $\Omega_{el}$. The plane in the
middle of the cell-chip junction, depicted in Fig.~\ref{fig:middleplan_2d},
is going to be the new two dimensional domain $\,\Omega_{2D}$: it
is equidistant from $\Gamma_{sub}\,$ and from $\Gamma_{ef}\cup\Gamma_{cell}$,
which are respectively placed at $z=-\delta_{j}/2$ and $z=+\delta_{j}/2$. 

We follow a procedure based on the integration of the three-dimensional
Eqs.~\eqref{eq:PNP_flux} on the test volume of Fig.~\ref{fig:middleplan_2d}.
This latter is a parallelepiped of volume $V_{xy}=\delta_{j}h_{x}h_{y}$,
where $h_{x}$ and $h_{y}$ are infinitesimally small.
Introducing 
the following integral means:
\begin{subequations}\label{eq:reduction-3d-2d}
\begin{equation}
\overline{\varphi}\left(t;x,y\right)=
\dfrac{1}{\delta_{j}}
\int_{-\frac{\delta_j}{2}}^{\frac{\delta_j}{2}}\varphi\left(t;x,y,z\right)dz
\qquad \overline{c_i}\left(t;x,y\right)=
\dfrac{1}{\delta_{j}}
\int_{-\frac{\delta_j}{2}}^{\frac{\delta_j}{2}}c_i\left(t;x,y,z\right)dz,
\label{eq:integral_mean}
\end{equation}
and integrating the continuity Eq.~\eqref{eq:continuity-1-1}, we obtain
\begin{align}
\dfrac{\partial\overline{c}_{i}}{\partial t}\delta_{j}h_{x}h_{y}\,+\: & \,\text{f}_{x}^{i}\left(t;x+\dfrac{h_{x}}{2},y,0\right)\delta_{j}h_{y}-\,\text{f}_{x}^{i}\left(t;x-\dfrac{h_{x}}{2},y,0\right)\delta_{j}h_{y}\nonumber \\
+\: & \,\text{f}_{y}^{i}\left(t;x,y+\dfrac{h_{y}}{2},0\right)\delta_{j}h_{x}-\,\text{f}_{y}^{i}\left(t;x,y-\dfrac{h_{y}}{2},0\right)\delta_{j}h_{x}\nonumber \\
+\: & \,\text{f}_{z}^{i}\left(t;x,y,\dfrac{\delta_{j}}{2}\right)h_{x}h_{y}-\,\text{f}_{z}^{i}\left(t;x,y,-\dfrac{\delta_{j}}{2}\right)h_{x}h_{y}=0,\label{eq:mediated_int}
\end{align}
$\text{f}_{x}$, $\text{f}_{y}$ and $\text{f}_{z}$ being the
components of the fluxes in the three directions. 
An analogous result is obtained for the Poisson 
equation~\eqref{eq:poisson-1-1}.
 
The values of the fluxes $\mathbf{f}_{i}$  and of the electric displacement
 $\mathbf{D}$ on $\Sigma_{top}$ and $\Sigma_{bot}$
are unknown and need be computed according to the boundary
conditions applied on these surfaces. The boundary conditions
on $\Sigma_{top}$ and $\Sigma_{bot}$ can be written as:
\begin{align}
\mathbf{f}_{i}\cdot\mathbf{n} & = f_{i}^{top}\left(t;c_{i}^{top},\overline{c}_{i},\varphi_{top},\overline{\varphi}\right) & & \text{on }
\Sigma_{top}
\\
\mathbf{f}_{i}\cdot\mathbf{n}& =  f_{i}^{bot}\left(t;c_{i}^{bot},\overline{c}_{i},\varphi_{bot},\overline{\varphi}\right)
& & \text{on }
\Sigma_{bot}
\\
\mathbf{D}\cdot\mathbf{n}& =  g_{top}\left(t;c_{i}^{top},\overline{c}_{i},\varphi_{top},\overline{\varphi}\right)
& & \text{on } \Sigma_{top}
\\
\mathbf{D}\cdot\mathbf{n}& =  g_{bot}\left(t;c_{i}^{bot},\overline{c}_{i},\varphi_{bot},\overline{\varphi}\right)
& & \text{on }\Sigma_{bot},
\end{align}
where the functions $f_{i}^{top}$,
$f_{i}^{bot}$, $g_{top}$ and $g_{bot}$ depend on the
``averaged'' quantities defined in~\eqref{eq:integral_mean}, 
but also on the quantities evaluated on the surfaces $\Sigma_{top}$ and $\Sigma_{bot}$, defined as:
\begin{eqnarray}
c_{i}^{top}:=\left.c_{i}\right|_{\Sigma_{top}} & \qquad\qquad & c_{i}^{bot}:=\left.c_{i}\right|_{\Sigma_{bot}}\label{eq:ctopbot}\\
\varphi_{top}:=\left.\varphi\right|_{\Sigma_{top}} &  & \varphi_{bot}:=\left.\varphi\right|_{\Sigma_{bot}}.\label{eq:phitopbot}
\end{eqnarray}
\end{subequations}
Dividing~\eqref{eq:mediated_int} and the analogue for the Poisson equation
by $\left|V_{xy}\right|=\delta_{j}h_{x}h_{y}$,  
and taking the limit as $h_{x},h_{y}\rightarrow0$, we obtain 
the following averaged 2.5D PNP model in
$\Omega_{2D}$:
\begin{subequations}\label{eq:ridotto_2D}
\begin{align}
\dfrac{\partial\overline{c}_{i}}{\partial t}+\text{div}_{xy}\overline{\mathbf{f}}_{i}+\dfrac{1}{\delta_{j}}f_{i}^{top}+\dfrac{1}{\delta_{j}}f_{i}^{bot}\,=\,0\label{eq:continuity_2d}\\
\overline{\mathbf{f}}_{i}\,=\,-D_{i}\left(\nabla_{xy}\overline{c}_{i}+\dfrac{z_{i}}{V_{th}}\overline{c}_{i}\nabla_{xy}\overline{\varphi}\right)\label{eq:flux_2d}\\
\text{div}_{xy}\overline{\mathbf{D}}+\dfrac{1}{\delta_{j}}g_{top}+\dfrac{1}{\delta_{j}}g_{bot}\,= \: q\sum_{i}z_{i}\overline{c}_{i}\label{eq:poiss_2d}\\
\overline{\mathbf{D}}\:=\:-\epsilon\nabla_{xy}\overline{\varphi}.\qquad\label{eq:displacement_2d}
\end{align}
The boundary conditions applied on $\partial\Omega_{2D}$ simply reduce to:
\begin{align}
\overline{c}_{i}=c_{i}^{bath}\label{eq:reduction_bc_ci}\\
\overline{\varphi}=\varphi_{bath}.\label{eq:reduction_bc_phi}
\end{align}
\end{subequations}
The reason for qualifying the novel electrodiffusive 
model~\eqref{eq:ridotto_2D} as a 2(+1/2)D=2.5D formulation is related
to the definition of the source flux terms 
$f_{i}^{top},f_{i}^{bot}$
and $g_{top},g_{bot}$, object of the next section.

\subsubsection{The reason for the "2.5D": boundary layer model 
approximation}
Let $\chi_{\zeta}$ denote
the characteristic function of a domain $\zeta\subset\mathbb{R}^{2}$.
Then, since the upper surface of $\Omega_{el}$ 
is the union of two different parts, $\Gamma_{cell}$
and $\Gamma_{ef}$ (cf. Fig.~\ref{fig:middleplan_2d}), 
we decompose $f_{i}^{top}$ and $g_{top}$ into the sum of two contributions as:
\begin{subequations}\label{eq:top_with_chi}
\begin{eqnarray}
f_{i}^{top} & = & f_{i,top}^{cell}\,\chi_{\Gamma_{cell}}+f_{i,top}^{ef}\,\chi_{\Gamma_{ef}}\label{eq:ftop_chi}\\
g_{top} & = & g_{top}^{cell}\,\chi_{\Gamma_{cell}}+g_{top}^{ef}\,\chi_{\Gamma_{ef}}.
\label{eq:gtop_chi}
\end{eqnarray}
\end{subequations}

\begin{figure}[tb]
\centering 
\includegraphics[width=0.7\textwidth]
{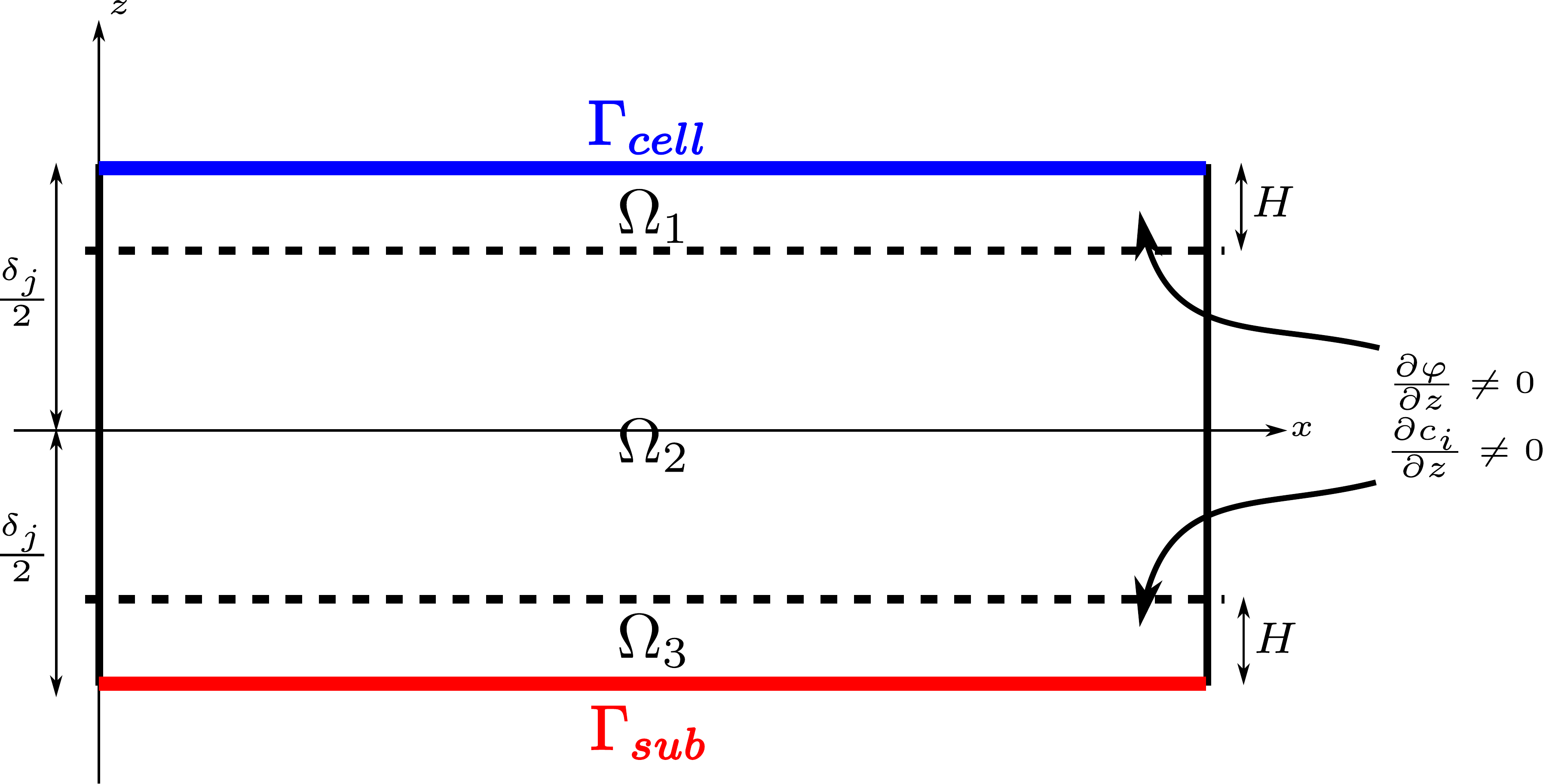}
\caption[Schematics for the boundary layers hypotheses.]
{\label{fig:boundary_layer}Cross section in the $x$-$z$
plane of the three-dimensional cleft, showing the scheme for the
model of the boundary layers near the surfaces $\Gamma_{cell}$
and $\Gamma_{sub}$ ($H$ is the layer amplitude). Decomposition into 
three subdomains:
$\Omega_{1}=\left\{ (x,z)\text{ s.t. }z\in\left[\delta_{j}/2-H;\:\delta_{j}/2\right]\right\} $,
$\Omega_{2}=\left\{ (x,z)\text{ s.t. }z\in\left[-\delta_{j}/2+H;\:\delta_{j}/2-H\right]\right\} $
and $\Omega_{3}=\left\{ (x,z)\text{ s.t. }z\in\left[-\delta_{j}/2;\:-\delta_{j}/2+H\right]\right\} $ }
\end{figure}

\begin{figure}[tb]
\centering
\subfigure[Potential]
{\includegraphics[width=0.47\textwidth]{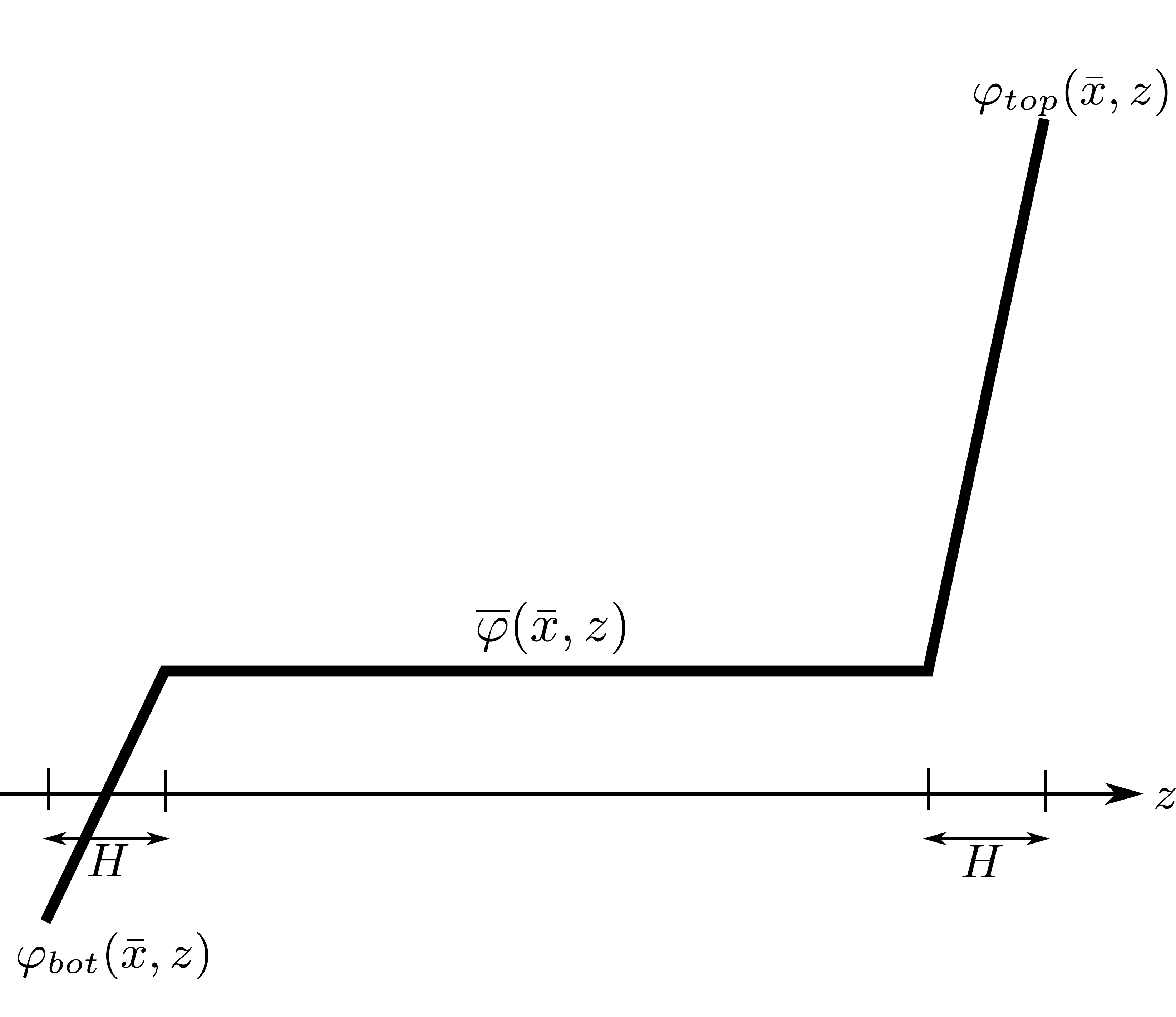}
\label{fig:pot_z}}
\subfigure[Concentration]
{\includegraphics[width=0.47\textwidth]{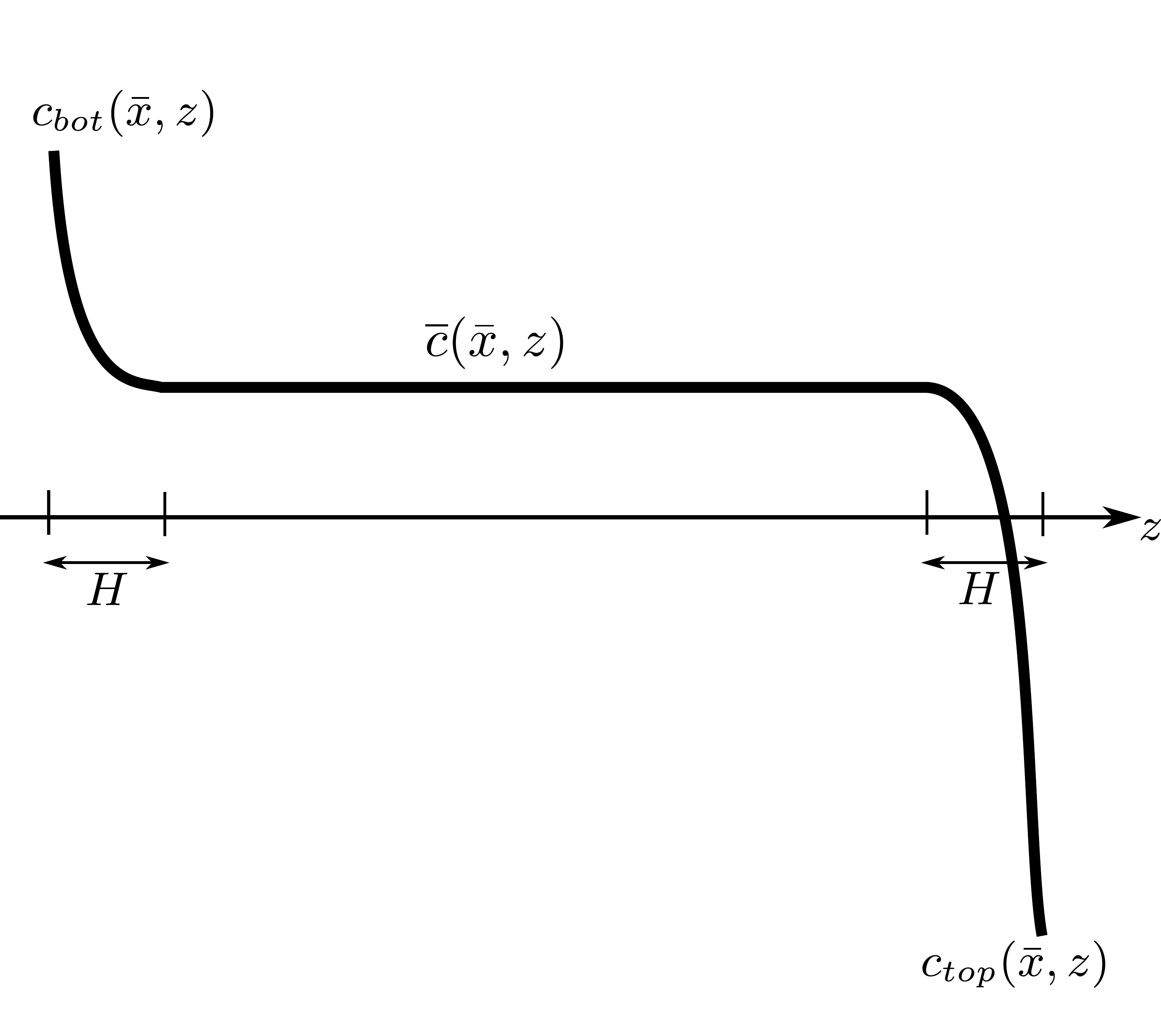}
\label{fig:conc_z}}
\caption{Scheme of the distributions in the $z$-direction 
at a fixed $\bar{x}$
for $\varphi(\bar{x},z)$ and for $c_{i}(\bar{x},z)$ (where the considered
ion is positively charged) for the cross section depicted in 
Fig.~\ref{fig:boundary_layer}.}
\label{fig:distributions_z}
\end{figure}

The coupling conditions on $\,\Gamma_{cell}$ and $\,\Gamma_{sub}$
are between different environments and typically 
give rise to the occurrence of boundary layers~\cite{brera2009multiphysics,gardner2013simulation,mori2009numerical}.
These latter are in the form of electrical double layers, of which
we only account for the diffuse layer, neglecting the ions attached to the
surfaces, as in the Gouy-Chapman approximation~\cite{grahame1947electrical}.
In Fig.~\ref{fig:boundary_layer}, we focus our attention on
a $x$-$z$ cross section of the whole three-dimensional electrolyte cleft 
at $y=\bar{y}$, denoted $\Omega_{xz}$. This latter is partitioned
into three distinct subdomains as $\Omega_{xz}=\Omega_{1}\cup\Omega_{2}\cup\Omega_{3}$, 
$H$ being the amplitude of the two boundary layer regions 
$\Omega_{2}$ and $\Omega_{3}$. 
According to physical
evidence, for every fixed point $\bar{x}$ of the $x$ axis, we assume that 
\[
\dfrac{\partial\varphi(\bar{x},z)}{\partial z}=\dfrac{\partial c_{i}(\bar{x},z)}{\partial z}=0\quad\text{in}\:\Omega_{2},
\]
and we set $\varphi(\bar{x},z)=\overline{\varphi}(\bar{x},z)\text{ and }c_{i}(\bar{x},z)=\overline{c}_{i}(\bar{x},z)\text{ for all }z\in\Omega_{2}$.
These two definitions amount to extending along the $z$-direction
(in the sole interval $\Omega_{2}$) the averaged values determined
by the model described in Sect.~\ref{sub:Model-reduction}. We also
introduce further assumptions on the electric potential and the particle
fluxes in the two boundary layer subdomains:
\begin{enumerate}
\item $\varphi$ is linear in $\Omega_{1}$ and $\Omega_{3}$ and continuous
at $z=\delta_{j}/2-H$ and at $z=-\delta_{j}/2+H$;
\item $\mathbf{f}_{i}$ is constant in $\Omega_{1}$ and $\Omega_{3}$.
\end{enumerate}
The spatial distribution of $\varphi(\bar{x},z)$ for a fixed point
$\bar{x}$ is schematically depicted in Fig.~\ref{fig:pot_z}.
Assumption 1. indicates that the electric field is piecewise constant
over $\Omega_{xz}$ (and equal to zero in $\Omega_{2}$). Also the
particle fluxes are piecewise constant over $\Omega_{xz}$ (and equal
to zero in $\Omega_{2}$ because both drift and diffusion terms are
null there). 
In order to determine the concentration $c_{i}(\bar{x},z)$, we integrate
the Nernst-Planck transport equation~\eqref{eq:nernstplanck-1} in
$\Omega_{1}$ and $\Omega_{3}$. The resulting distribution of ions
is piecewise exponential over $\Omega_{xz}$, continuous at $z=\delta_{j}/2-H$
and at $z=-\delta_{j}/2+H$, and constant in $\Omega_{2}$, as depicted
in Fig.~\ref{fig:conc_z}. The corresponding mathematical
expressions for the boundary fluxes $f_{i,top}^{cell}$ and $f_{bot}$
are:
\begin{subequations}\label{eq:scharf_gumm_fluxes}
\begin{eqnarray}
f_{i,top}^{cell} =  -\dfrac{D_{i}}{H}\left(Be\left(-\dfrac{z_{i}(\varphi_{top}-\overline{\varphi})}{V_{th}}\right)c_{i}^{top}-Be\left(\dfrac{z_{i}(\varphi_{top}-\overline{\varphi})}{V_{th}}\right)\overline{c}_{i}\right)\label{eq:fi_top_sg}\\
f_{i}^{bot} =  -\dfrac{D_{i}}{H}\left(Be\left(\dfrac{z_{i}(\overline{\varphi}-\varphi_{bot})}{V_{th}}\right)c_{i}^{bot}-Be\left(-\dfrac{z_{i}(\overline{\varphi}-\varphi_{bot})}{V_{th}}\right)\overline{c}_{i}\right).
\label{eq:fi_bot_sg}
\end{eqnarray}
\end{subequations}
The above described modeling reduction procedure is equivalent to
applying the Scharfetter-Gummel (SG) exponentially fitted approximation
in $\Omega_{1}$ and $\Omega_{3}$~\cite{scharfetter1969large}. 
Regarding the electric displacement, with the above approximation
we obtain:
\begin{subequations}\label{eq:electr_displ_fluxes}
\begin{eqnarray}
g_{top}^{cell}   -\epsilon\dfrac{\varphi_{top}-\overline{\varphi}}{H}\label{eq:gi_top_sg}\\
g_{bot} =  -\epsilon\dfrac{\varphi_{bot}-\overline{\varphi}}{H}.\label{eq:gi_bot_sg}
\end{eqnarray}
\end{subequations}
To determine the values $c_{i}^{top}$, $c_{i}^{bot}$, $\varphi_{top}$
and $\varphi_{bot}$ we use \eqref{eq:scharf_gumm_fluxes}
into the coupling boundary conditions~\eqref{eq:flux_through_membr}
and~\eqref{eq:JdotnSub}, and \eqref{eq:electr_displ_fluxes}
into the coupling boundary conditions~\eqref{eq:CM_condition}
and~\eqref{eq:CS_condition}, to obtain:
\begin{subequations}\label{eq:top_bot_final} 
\begin{gather}
\left.c_{i}^{top}\right|_{\Gamma_{cell}}=\frac{1}{Be\left(-z_{i}(\varphi_{top}-\overline{\varphi})/V_{th}\right)}\left(\overline{c}_{i}Be\left(z_{i}(\varphi_{top}-\overline{\varphi})/V_{th}\right)+\dfrac{j_{i}^{top}H}{qz_{i}D_{i}}\right)\\
\left.c_{i}^{bot}\right|_{\Gamma_{sub}}=\frac{Be\left(-z_{i}(\overline{\varphi}-\varphi_{bot})/V_{th}\right)}{Be\left(z_{i}(\overline{\varphi}-\varphi_{bot})/V_{th}\right)}\overline{c}_{i}\\
\left.\varphi_{top}\right|_{\Gamma_{cell}}=\dfrac{1}{C_{M}+\epsilon/H}\left(C_{M}V_{cell}+\dfrac{\epsilon}{H}\overline{\varphi}\right)\\
\left.\varphi_{bot}\right|_{\Gamma_{sub}}=\dfrac{1}{C_{S}+\epsilon/H}\left(C_{S}V_{G}+\dfrac{\epsilon}{H}\overline{\varphi}\right).
\end{gather} 
\end{subequations}
Let us now consider the artificial surface $\Gamma_{ef}$. 
No boundary layer is expected to occur there, 
so that we simply set
\begin{equation}\label{eq:cond_top_ef}
\left.\varphi_{top}\right|_{\Gamma_{ef}} = \overline{\varphi}, \qquad
\left.c_{i}^{top}\right|_{\Gamma_{ef}} = \overline{c}_{i}
\end{equation}
and use these values into the
electrolyte-electrolyte coupling conditions~\eqref{eq:flux_through_el} and~\eqref{eq:cstar_first}
to compute the functions $g_{top}^{ef}$ and $f_{i,top}^{ef}$ 
introduced in~\eqref{eq:top_with_chi}.
The final 2.5D ion electrodiffusion model in $\Omega_{2D}$ then reads:
\begin{subequations}\label{eq:final_2d_layers}
\begin{gather}
\dfrac{\partial\overline{c}_{i}}{\partial t}+\text{div}_{xy}\overline{\mathbf{f}}_{i}+\dfrac{1}{\delta_{j}}f_{i,cell}^{top}\left.\chi\right|_{\Gamma_{cell}}+\dfrac{1}{\delta_{j}}f_{i}^{bot}+  \dfrac{1}{\delta_{j}}v^{*}\left(\overline{c}_{i}-c_{i}^{bath}\right)\left.\chi\right|_{\Gamma_{ef}}=0\\
\overline{\mathbf{f}}_{i}=-D_{i}(\nabla_{xy}\overline{c}_{i}+\dfrac{z_{i}}{V_{th}}  \overline{c}_{i}\nabla_{xy}\overline{\varphi})\\
\text{div}_{xy}\overline{\mathbf{D}}+\dfrac{1}{\delta_{j}}g_{top}^{cell}\left.\chi\right|_{\Gamma_{cell}}+\dfrac{1}{\delta_{j}}g_{bot}+\dfrac{1}{\delta_{j}}C^{*}  \left(\overline{\varphi}-V_{bath}\right)\left.\chi\right|_{\Gamma_{ef}}=q\sum_{i}z_{i}\overline{c}_{i}\\
\overline{\mathbf{D}}=-\epsilon\nabla_{xy}\overline{\varphi}.
\end{gather}
\end{subequations}
System~\eqref{eq:final_2d_layers} is completed by the same kind of
initial conditions as for the 3D PNP system~\eqref{eq:PNP_flux}, 
by the boundary conditions~\eqref{eq:reduction_bc_ci}-\eqref{eq:reduction_bc_phi}
and by the set of 
relations~\eqref{eq:scharf_gumm_fluxes}-\eqref{eq:cond_top_ef}
 
\subsection{A 2D electrical model for ion transport:
the Area Contact formulation\label{sub:A-simplified-version}}
In the same spirit as done for single cells by Fromherz et al. in~\cite{Brittinger2005,waser2012nanoelectronics}
and for multi-electrode arrays (MEAs) in ~\cite{joye2008electrical}, 
it is possible to obtain from the 2.5D equation 
system~\eqref{eq:ridotto_2D}
a genuine 2D ion transport model, denoted Area Contact model. 
To this purpose, following the approach of~\cite{Brittinger2005}, 
we consider only the attached area as computational domain $\Omega_{2D}$, 
we neglect the variation of $\varphi$ and $c_i$ in the 
$z$-direction and 
we also assume that ion densities are spatially homogeneous. 
In this manner, the functions $f_{i}^{top}$, $f_{i}^{bot}$, $g_{top}$
and $g_{bot}$ are still computed using the 3D boundary conditions of Sect.~\ref{sub:Boundary-and-initial} but setting $\varphi_{top}=\varphi_{bot}=\overline{\varphi}$ and $c_{i}^{top}=c_{i}^{bot}=\overline{c}_{i}$. 

Omitting from now the notation $\,\overline{(\cdot)}\,$, we sum
the $M$ continuity equations~\eqref{eq:continuity_2d}
and use the Poisson equation and the boundary coupling conditions defined
in~\eqref{eq:CM_condition} and~\eqref{eq:CS_condition}
to express the time derivative of $\rho$ as
\begin{equation*}
\dfrac{\partial}{\partial t}\left(\rho\delta_{j}\right)=\dfrac{\partial}{\partial t}\left(\text{div}_{xy}\mathbf{D}\right)+\dfrac{\partial}{\partial t}\left(C_{M}\varphi+C_{S}\varphi\right)-\dfrac{\partial}{\partial t}\left(C_{M}V_{cell}+C_{S}V_{G}\right).\label{eq:rho_areacont}
\end{equation*}
Replacing the previous relation into the sum of the continuity equations 
we end up with the 2D Area Contact model for ion transport:
\begin{subequations}\label{eq:BF_model_phi}
\begin{gather}
\left(C_{M}+C_{S}\right)\dfrac{\partial\varphi}{\partial t}+\text{div}_{xy}\left(\mathbf{j}_{tot}^{cond}-\dfrac{\partial\left(\epsilon\nabla\varphi\right)}{\partial t}\right)\:=\: j_{tot}^{tm}  +\dfrac{\partial}{\partial t}\left(C_{M}V_{cell}+C_{S}V_{G}\right)\label{eq:model_BF_phi}\\
\mathbf{j}_{tot}^{cond}\,=\,-q\sum\left|z_{i}\right|\mu_{i}c_{i}  \delta_{j}\nabla\varphi.\label{eq:model_BF_phi2}
\end{gather}
For given ion concentrations $c_i$,
system~\eqref{eq:BF_model_phi} is a parabolic boundary value problem for 
the dependent variable $\varphi=\varphi(x,y, t)$. 
Following the terminology adopted in~\cite{Brittinger2005}, we refer
to system~\eqref{eq:BF_model_phi} as the \lq\lq 2D electrical model\rq\rq.
Ion dynamics should also be accounted for, therefore
at each time level we first solve~\eqref{eq:BF_model_phi} and 
compute the integral mean of $\varphi$ over $\Omega_{2D}$
\begin{equation}
V_{J}(t):=\dfrac{\int_{\Omega_{2D}}
\varphi\left(x,y, t\right)dxdy}{\left|\Omega_{2D}\right|},\label{eq:VJ_def}
\end{equation}
\end{subequations}
$\left|\Omega_{2D}\right|$ denoting the area of the contact area.
Then, we use $V_{J}(t)$ as an input voltage to
solve the system of ODEs corresponding to the electrical equivalent circuit
proposed in~\cite{Brittinger2005}, 
to determine the concentrations $c_i$ that must be passed 
to~\eqref{eq:model_BF_phi2} to advance to the next time level.
Following again the terminology adopted in~\cite{Brittinger2005}, 
we refer to this latter PDE-ODE coupled system  as
the \lq\lq 2D electrodiffusion model\rq\rq. The 
physical accuracy of the 2D models introduced above
are investigated in Sect.~\ref{sec:area_contact_simul}.

\section{Functional Iterations and Numerical Discretization}
\label{sec:methods}
In this section we describe the techniques used to
numerically solve the mathematical models introduced in 
Sections~\ref{sec:models} and~\ref{sec:hierarchy}.
The adopted strategy is composed of three steps: 
\begin{enumerate}
\item Temporal discretization
\item Linearization
\item Spatial discretization
\end{enumerate}
\subsection*{Step 1 \label{sub:time_disc}}
For the temporal discretization
we adopt the Backward-Euler scheme
to approximate all time derivatives. 
Since in most applications considered in this work the
input signal (usually the intracellular potential $V_{cell}(t)$)
is expressed as a combination of Heaviside functions, the time step 
of temporal advancement $\Delta t$ is a-priori appropriately chosen 
in numerical simulations according to the following strategy: 
in correspondence of the on/off switching time
of the signal, $\Delta t$ is set equal to a small value, 
in the order of \SI{1e-8}{\second}; then, once transients are exhausted, 
$\Delta t$ is suitably increased up to a value in the order of \SI{1e-3}{\second}.

\subsection*{Step 2 \label{sub:linear}}
In order to handle the intrinsic nonlinearity of the models,
we apply a functional iteration procedure that is widely used in 
the decoupled solution of the DD semiconductor device equations. 
The method is the well known Gummel Map
\cite{jerome1996analysis,kerkhoven1988spectral,gummel1964self,bosisio2000discretization},
a staggered algorithm where each variable of the problem and its corresponding
equation are treated in sequence until convergence.
Gummel's map in the steady state regime 
is theoretically investigated in the seminal 
book~\cite{jerome1996analysis} to which we refer for details and further bibliography. Under proper conditions on problem geometry and 
boundary data, the Gummel decoupled iteration is proved to admit 
a fixed point and also to be a contraction. This implies uniqueness 
of the solution of the nonlinear PDE equation system under investigation.
The convergence rate predicted by the analysis 
of~\cite{jerome1996analysis} is linear, but computational experience 
reveals that the Gummel algorithm in stationary conditions 
is exceptionally rapid and robust with respect to the choice of the 
initial guess, which makes it much preferable than Newton's method
despite its theoretically predicted lower order of convergence. 

\subsection*{Step 3}
\label{sub:finite_el}
For the spatial discretization of the linearized PDEs
we adopt the piecewise linear conforming 
Galerkin-finite element method (G-FEM) stabilized by means of an
exponential fitting technique (Edge Averaged Finite Element method (EAFE))~\cite{bank1998finite,de2006quantum,gatti1998new,xu1999monotone}, in order to deal with possibly dominating drift terms and avoid the onset of spurious oscillations in the computed solutions.

\subsection{The EAFE method in axisymmetric geometries}
In this section we illustrate the extension of
the EAFE method proposed in~\cite{xu1999monotone} to treat the case of
three-dimensional electrodiffusive problems in 
axisymmetric geometries as in the
numerical simulations reported in Sect.~\ref{sec:results}.
\begin{figure}[h!]
\centering 
\includegraphics[width=0.5\textwidth]{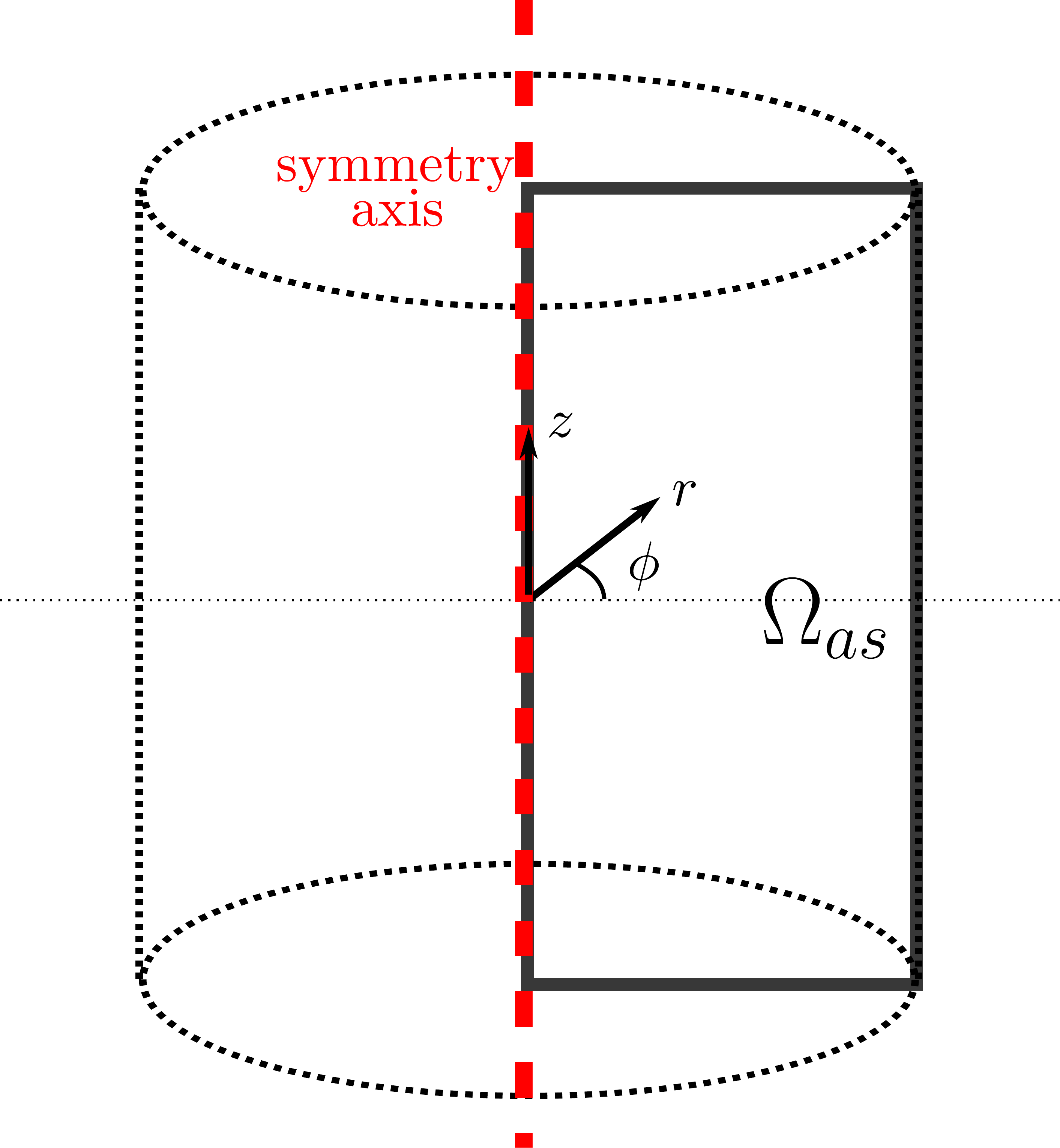}
\caption[Schematics of an axisymmetric domain]
{\label{fig:axialsymm}Schematics of an axisymmetric configuration.}
\end{figure}

\subsubsection{The electrodiffusion model problem}
Let $\Omega_{as} = (0, R) \times (0, Z)$ 
be the computational domain shown in 
Fig.~\ref{fig:axialsymm} with boundary $\Gamma:=\partial \Omega_{as}$
divided into two disjoint subsets $\Gamma_D$ and $\Gamma_N$ and 
with outward unit normal vector $\mathbf{n}=[n_r, \, n_z]^T$.
Let $u$ be the dependent variable (having the physical meaning of 
ion density), related
to the ancillary dependent variable $n$ through the constitutive 
equation
\begin{subequations}\label{eq:eafe_cyl}
\begin{equation}
u = n e^{\psi}
\label{eq:slotboom}
\end{equation}
where $\psi$ is a given function representing a normalized electric
potential. The change of variable~\eqref{eq:slotboom} allows
to write the DD flux density~\eqref{eq:nernstplanck-1} in an equivalent diffusive form and $n$ undergoes the name of \emph{Slotboom variable}~\cite{slotboom}. The model of 
electrodiffusion of $n$ that we consider in the present section is 
the following boundary value problem in self-adjoint form:
\begin{equation}
\begin{cases}
\dfrac{1}{r}\dfrac{\partial}{\partial r}\left(rJ_{r}\right)+\dfrac{\partial}{\partial z}J_{z}+ce^{\psi}n=f & \text{in}\;\Omega_{as}\\
n=0 & \text{on}\;\Gamma_{D}\\
-\mathbf{J}\left(n\right)\cdot\mathbf{n}=j_{N} & \text{on}\;\Gamma_{N},
\end{cases}\label{eq:model_cyl_axial}
\end{equation}
where $\mathbf{J}(n) = [J_{r}(n) \, J_{z}(n) ]^T =
\mu e^{\psi} [\partial n/\partial r, \, 
\partial n/\partial z]^T$, 
$\mu\in C^{0}(\overline{\Omega}_{as})$ being the ion mobility such that $\mu=\mu\left(x\right)\geq\mu_{0}>0\;\forall x\in\overline{\Omega}_{as}$,
and where $\psi$ is a continuous piecewise linear function over $\,\overline{\Omega}_{as}$ such that the drift field is
$\mathbf{b}:=\nabla\psi$ (a normalized electric field). We also assume
the reaction coefficient $c\in L^{\infty}(\Omega_{as})$, 
with $c\geq0\;\text{a.e.\ in }\Omega_{as}$, 
the source term $f\in L^{2}\left(\Omega_{as}\right)$
and the boundary flux $j_{N} \in L^2\left(\Gamma_N\right)$.
For any given functions $\phi$ and $\theta$ belonging to 
$L^2(\Omega_{as})$, 
let us endow $L^{2}(\Omega_{as})$ with the weighted scalar product 
$\left\langle \cdot,\cdot\right\rangle _{\omega}$
\begin{equation}
\left\langle \phi, \theta \right\rangle _{\omega}: = \int_{0}^{Z}\int_{0}^{R}\widetilde{\phi}(r,z) \widetilde{\theta}(r,z)drdz=\left\langle \widetilde{\phi},\widetilde{\theta}\right\rangle,  \label{eq:scalar_prod} 
\end{equation}
where the weight function is $\omega=\sqrt{r}$, 
$\widetilde{\phi}:=\omega\phi$, 
$\widetilde{\theta}:=\omega\theta$ and 
$\left\langle \cdot,\cdot\right\rangle$ is the usual scalar product
in $L^2(\Omega_{as})$.
Finally, let 
\begin{equation}
V:=H_{\Gamma_{D}}^{1}\left(\Omega_{as}\right)=\left\{ v\in H^{1}\left(\Omega_{as}\right):\:\left.v\right|_{\Gamma_{D}}=0\right\} ,\label{eq:space_V}
\end{equation}
endowed with the equivalent norm, using Poincar\'e-Friedrichs 
inequality (see~\cite{quarteroni2008numerical}, Chapt.~1)
\begin{equation}
\left\Vert w\right\Vert _{V}:=\left\Vert \nabla w\right\Vert _{\omega}
= \left( \left\langle \nabla w,\nabla w\right\rangle_{\omega}\right)^{1/2} 
\qquad \forall w \in V.\label{eq:V_norm_def}
\end{equation}
Then, the weak formulation of problem~\eqref{eq:model_cyl_axial} reads:\\
find $n\in V$ such that:
\begin{align}
&
a_{\omega}\left(n,v\right) = F_{\omega}\left(v\right) &\qquad\forall v\in V\label{eq:weak_form}
\end{align}
where:
\begin{align}
& a_{\omega}\left(n,v\right) = -\left\langle \mu e^{\psi}\nabla n,\nabla v\right\rangle _{\omega}+\left\langle ce^{\psi}n,v\right\rangle _{\omega} 
\label{eq:bilinear} \\
&
F_{\omega}\left(v\right) = \left\langle f,v\right\rangle _{\omega}+\int_{\Gamma_{N}}j_{N}vds_{\omega}  \label{eq:functional_rhs}
\end{align}
and where $ds_{\omega}=rds$
is the curvilinear abscissa in radial coordinates, 
$ds$ being the usual curvilinear abscissa.
Using the Lax-Milgram lemma (cf.~\cite{quarteroni2008numerical}, Chapt.~5) 
we can prove existence and uniqueness of the solution 
$n \in V$ of~\eqref{eq:weak_form} and the following stability
estimate~\cite{abbate2014}
\[
\left\Vert n\right\Vert _{V}\leq\dfrac{C_{P}\left\Vert f\right\Vert _{L^{2}(\Omega_{as})}+C_{N}\left\Vert j_{N}\right\Vert _{L^{2}(\Gamma_{N})}}{e^{\psi_{m}}\mu_{0}}
\]
where $C_{P}$ and $C_N$ are the Poincar\'e's and trace constants, respectively,
while $\psi_m$ is the minimum of $\psi$ in $\overline{\Omega}_{as}$.
Existence and uniqueness of $n$ obviously imply
the existence and uniqueness of the function $u$ that solves 
the boundary value problem corresponding to~\eqref{eq:model_cyl_axial} 
but written as the continuity equation~\eqref{eq:continuity-1-1}.

\subsubsection{Extension of the EAFE method to axisymmetric geometries}
Let $\mathcal{T}_{h}$ be a regular triangulation of the
domain $\Omega_{as}$ made by triangles $K$ 
(cf.~\cite{quarteroni2008numerical}, Chapt.~3).
On $\mathcal{T}_{h}$ we introduce the finite dimensional subspace 
$V_h \subset V$ made of piecewise linear conforming finite elements
vanishing on $\Gamma_D$.
Denoting by $n_h \in V_h$ the finite element approximation of $n$, 
the application of the EAFE method to problem~\eqref{eq:weak_form}
consists of finding $n_h \in V_h$ such that: 
\begin{align}
&
a_{\omega,h}\left(n_h,v_h\right) = F_{\omega}\left(v_h\right) &\qquad\forall v_h\in V_h\label{eq:weak_form_h}
\end{align}
where $a_{\omega,h}(\cdot, \cdot)$ is an approximate bilinear form 
constructed in such a way that the approximation 
$\mathbf{J}_h(n_h)$ of the flux $\mathbf{J}(n)$ over each 
triangle $K$ is a constant vector whose tangential 
component over each edge $e \in \partial K$ is computed 
by replacing the diffusion coefficient $\mu e^{\psi}\big|_e$ with 
its harmonic average along $e$. 
After computing the local stiffness matrix $A^{K}$ as in the case of 
Cartesian orthogonal coordinates, 
$A^{K}$ in still multiplied by the weighting factor $\int_{K}rdrdz$ (see~\cite{abbate2014,de2006quantum,porro2009solar}).
The following formula is used for the approximate 
evaluation of the local stiffness matrix 
of the EAFE method in axisymmetric geometries
\begin{align}
& A^{K}_{\omega}= 
\begin{bmatrix}\vspace{0.2cm}\overline{r}_{1} & 0 & 0\\
\vspace{0.2cm}0 &  \overline{r}_{2} & 0\\
0 & 0 & \overline{r}_{3}
\end{bmatrix} A^{K} & \forall K \in \mathcal{T}_h \label{eq:local_A_K}
\end{align}
$\overline{r}_{i}$ being the radial distance of 
the midpoint of edge $e_i$ from the origin of the coordinate system, 
$i=1,2,3$. The above approach proves 
to be numerically stable and accurate as demonstrated in
all the computational experiments reported in Sect.~\ref{sec:results}.
For the discretization of the local reaction and source terms, 
we adopt the same trapezoidal quadrature used in the Cartesian case 
and obtain
\begin{align}
M^{K}_{\omega}=\begin{bmatrix}\vspace{0.2cm}c_{1} e^{\psi_1}r_{1} & 0 & 0\\
\vspace{0.2cm}0 & c_{2}e^{\psi_2}r_{2} & 0\\
0 & 0 & c_{3}e^{\psi_3}r_{3}
\end{bmatrix}\dfrac{|K|}{3}, \label{eq:mass_K}\\
\mathbf{F}^{K}_{\omega}=
\begin{bmatrix}
f_{1}r_{1}\\
f_{2}r_{2}\\ 
f_{3}r_{3}
\end{bmatrix}\dfrac{|K|}{3} +
\dfrac{1}{2}
\begin{bmatrix}
(j_{N,2}+ j_{N,3})\overline{r}_{1} |e_1| \delta_{e_1}
\\
(j_{N,3}+ j_{N,1})\overline{r}_{2} |e_2| \delta_{e_2}\\
(j_{N,1}+ j_{N,2})\overline{r}_{3} |e_3| \delta_{e_3}\\
\end{bmatrix} & \qquad \forall K \in \mathcal{T}_h
\label{eq:rhs_K}
\end{align}
where ${r}_{i}$ is the radial distance of 
node $i$ from the origin of the coordinate system, $i=1,2,3$,
$|K|$ is the area of $K$, $|e_i|$ is the length of
edge $e_i$, $i=1,2,3$, while $\delta_{e_i}$ is equal to 
1 if $e_i \in \Gamma_N$ and 0 otherwise, and, finally,
$c_i$, $f_i$, $j_{N,i}$ and $\psi_i$ are the values of the $\mathbb{P}_1$-interpolants 
of $c$, $f$, $\psi$ and $j_N$ at each node $i$ of $K$, $i=1,2,3$.
Upon assembling~\eqref{eq:local_A_K},~\eqref{eq:mass_K} 
and~\eqref{eq:rhs_K} over 
the grid and applying the inverse of~\eqref{eq:slotboom} at each node of $\mathcal{T}_h$~\cite{Brezzi1989_1,brezzi_1989_2}, we end up with the following linear algebraic
system 
\begin{equation}
\Sigma_{\omega}\mathbf{u}\,=\,\mathbf{F}_{\omega},\label{eq:global_EAFE}
\end{equation}
where $\mathbf{u}$ is the vector of nodal values of the dependent variable
$u$ while $\Sigma_{\omega}$ and $\mathbf{F}_{\omega}$ are the global
stiffness matrix and load vector of the EAFE method, respectively.
As proved in~\cite{abbate2014}, $\Sigma_{\omega}$ 
is an irreducible M-matrix with respect to its columns.
This implies that~\eqref{eq:global_EAFE} admits a unique solution 
and that the discrete maximum principle 
holds under the same conditions
as in the Cartesian case studied in~\cite{xu1999monotone}. 
In particular, if $\,\mathbf{F}_{\omega}\geq \mathbf{0}$
(in the componentwise sense), then 
the solution of~\eqref{eq:global_EAFE} is such that 
$\mathbf{u} > \mathbf{0}$.
\end{subequations}

\section{Numerical Results}
\label{sec:results}

In this section, we carry out an extensive validation of all the mathematical
models discussed in Sects.~\ref{sec:models} and~\ref{sec:hierarchy}.
To this purpose, we divide the conducted simulations into two categories:
\begin{itemize}
\item validation of the PNP model of Sect.~\ref{sec:models} 
in 3D axisymmetric geometries;
\item validation of the model reduction of Sect.~\ref{sec:hierarchy}.
\end{itemize}

\subsection{Convergence analysis}\label{sec:convergence_analysis}
The numerical schemes of Sect.~\ref{sec:methods} have been implemented 
in \textsf{Octave} using
the Octave-Forge package \textsf{bim}~\cite{abbate2014} for matrix assembly.
The need to resort to radial and cylindrical
coordinates is intrinsic in most of the geometries considered in the
description of bio-hybrid devices, therefore we extended the 
package \textsf{bim} for these configurations and accurately validated the code
with a broad range of test cases.
Here we discuss the results of a convergence analysis carried out
on the two dimensional advection-diffusion problem~\eqref{eq:model_cyl_axial},
solved on the square domain $\Omega_{as} = \left[1,2\right] \times \left[0,1\right]$
in the $r$-$z$ plane (the symmetry axis is $r=0$).
We consider the case with $\mu=1$, 
$\mathbf{b} = \nabla\psi=\left[1,1\right]^T$, $c=1$, 
and where $f$ and the boundary conditions enforced on $\Gamma_{as}$ 
are chosen in such a way that the exact solution is 
\begin{equation*}
u\left(r,z\right)=z^{2}\ln r.
\end{equation*}

\begin{figure}[h!]
\centering 
\includegraphics[width=0.55\textwidth]{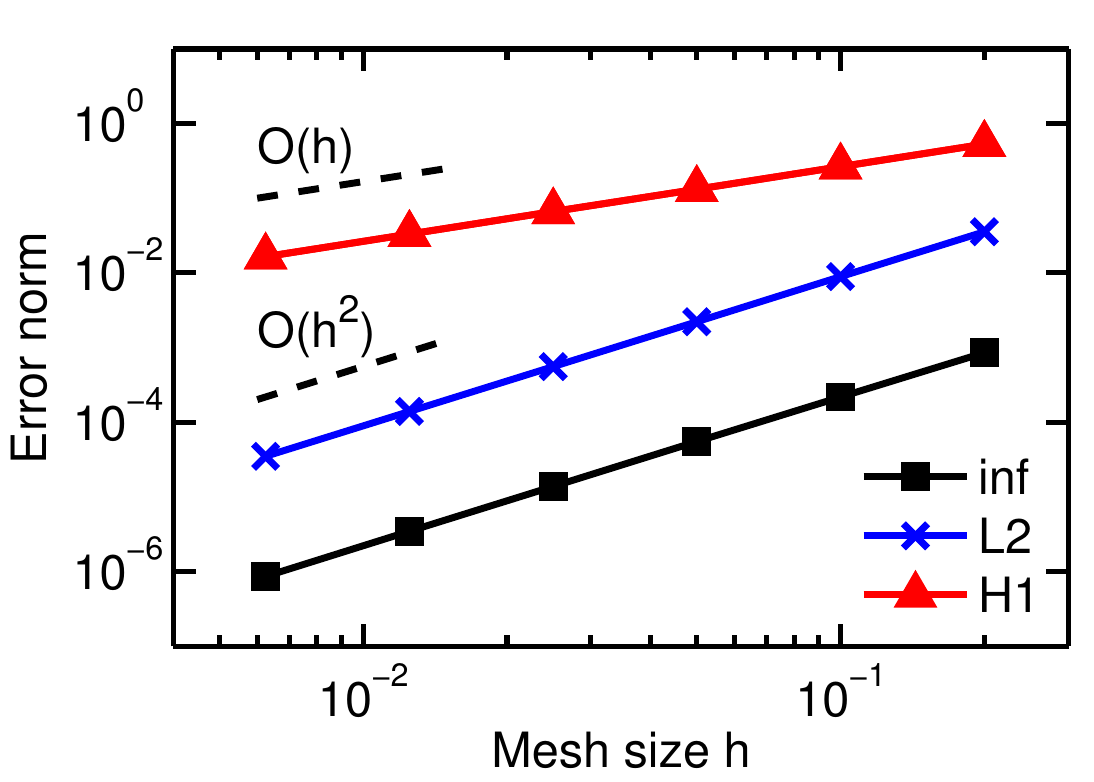}
\caption[Convergence analysis for the test case]
{Convergence analysis: 
$\left\Vert u-u_{h}\right\Vert _{L^\infty\left(\Omega_{as}\right)}$,
$\left\Vert u-u_{h}\right\Vert _{L^{2}\left(\Omega_{as}\right)}$
and
$\left\Vert u-u_{h}\right\Vert _{H^{1}\left(\Omega_{as}\right)}$
as a function of the mesh size $h$.}
\label{fig:testcase_conv}
\end{figure}

Fig.~\ref{fig:testcase_conv} reports the values of the $L^\infty$, $L^2$ and
$H^1$-norm of the error $u - u_{h}$ as a function of the mesh size $h$.
It is remarkable to notice that the numerical solution $u_{h}$ 
shows the same convergence properties in the energy norm
proved for the EAFE method in Cartesian coordinates
in~\cite{xu1999monotone,gatti1998new,lazarov2012exponential}.
Moreover, results clearly indicate superconvergence of the scheme 
in the $L^2$ and $L^\infty$ norms, as with usual piecewise linear
finite elements (see~\cite{quarteroni2008numerical}, Chapt.~6).

\subsection{Voltage-clamp stimulation: validation of the domain reduction}
\label{sec:voltage_clamp_validation}

\begin{figure}[h]
\centering 
\subfigure[Geometry]
{\includegraphics[width=0.35\textwidth]{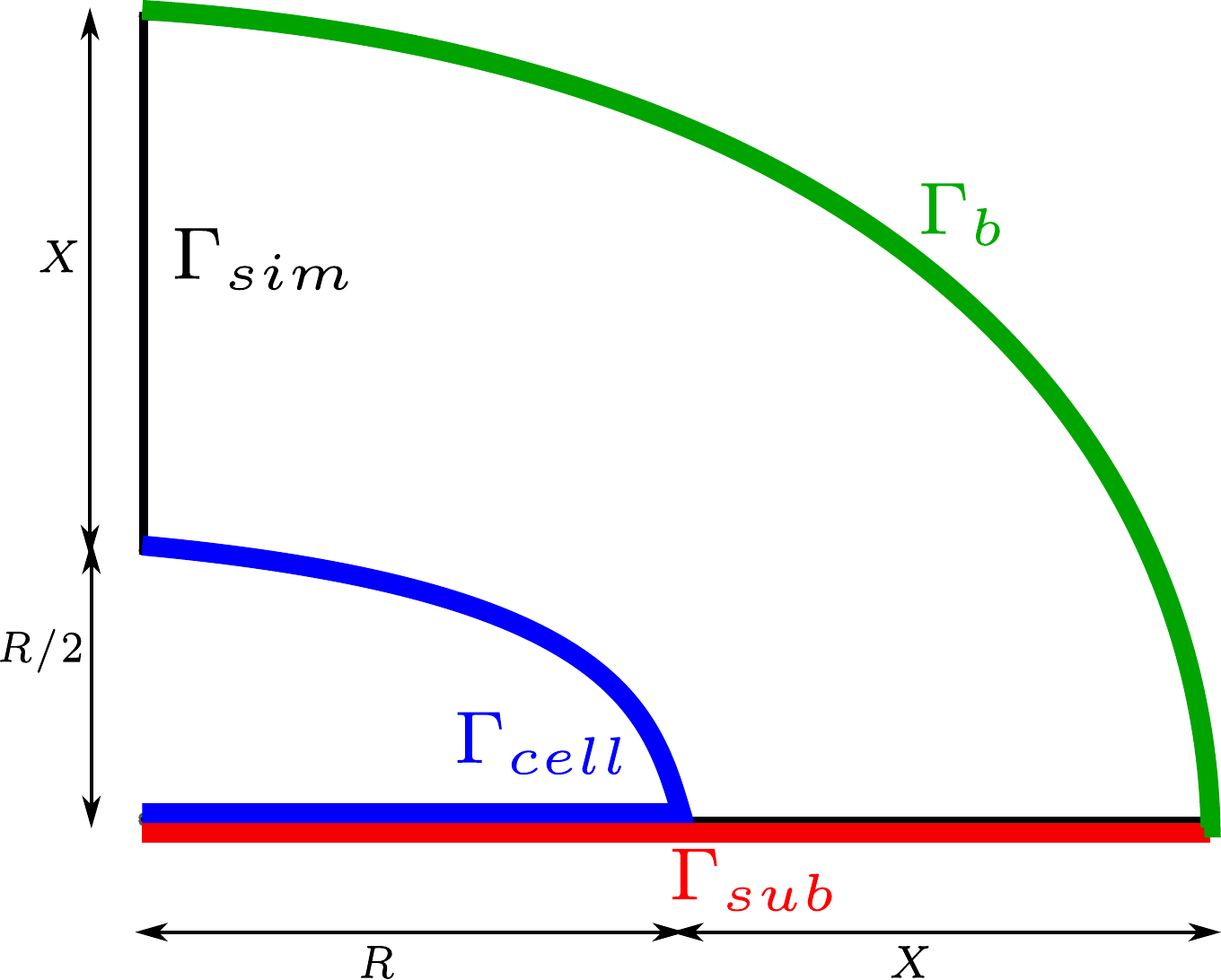}
\label{fig:geom_pabst2d}}
\hfill
\subfigure[Mesh]
{\includegraphics[width=0.58\textwidth]{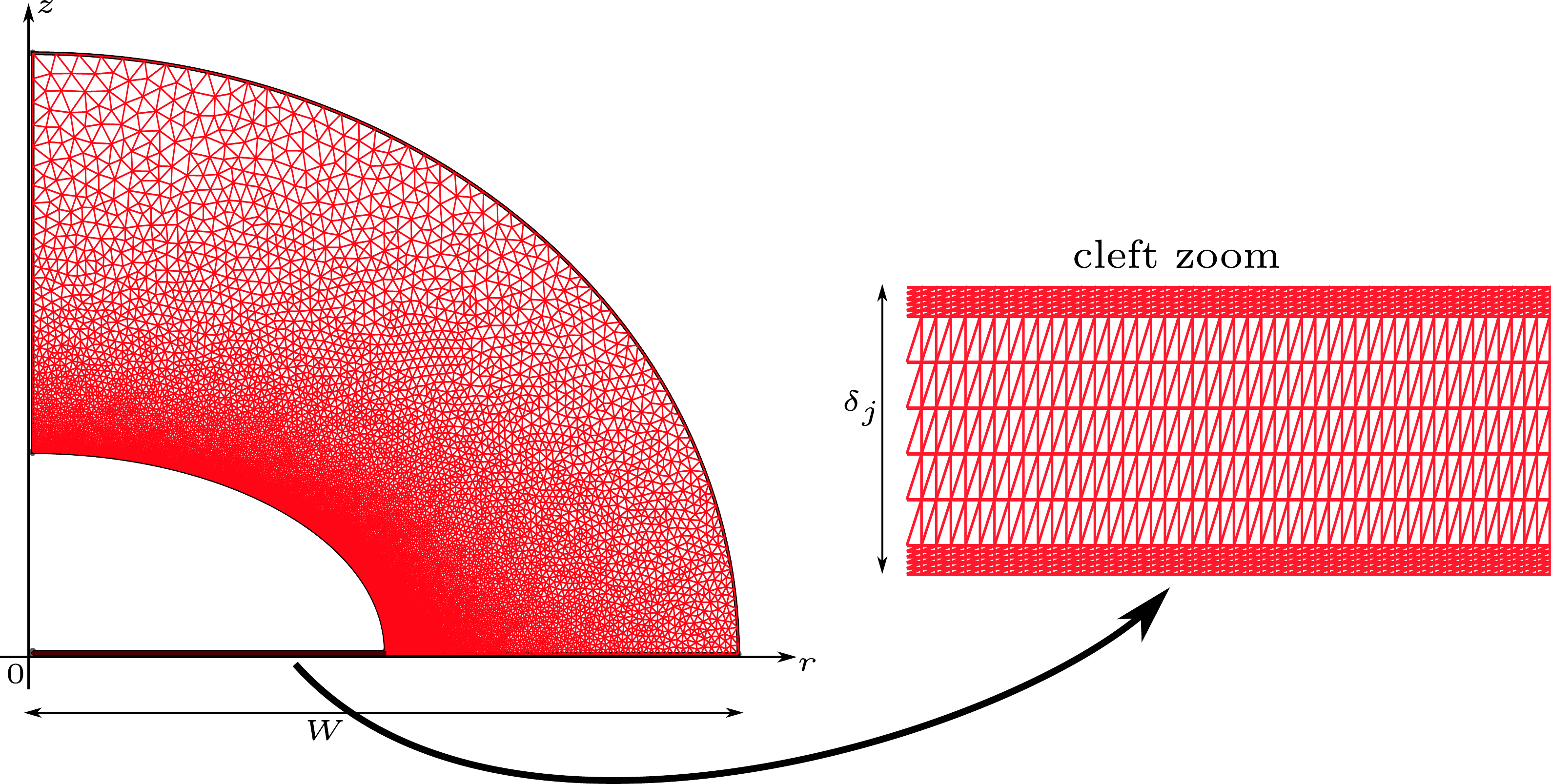}
\label{fig:mesh_pabst2d}}
\caption{Left: 2D geometry of the electrolyte surrounding the cell (cross
 section in the $r$-$z$ plane). Dimensions: $R=$\SI{10}{\micro\meter}
and $X=R$. The cell is approximated
as an ellipsoid with the major semiaxis equal to $R$
and the other one equal to $R/2$. The cleft
is the line between $\Gamma_{cell}$ and $\Gamma_{sub}$: its height
is $\delta_{j}=$ \SI{100}{\nano\meter}. Right: computational mesh, refined all
around the cell (in the zoom of a part of the cleft zone:
the mesh is structured and refined at the boundaries).}
\label{fig:pabst2d}
\end{figure}

As a first
numerical experiment, we consider a configuration comprising
the entire electrolyte bath
surrounding the cell, 
in order to demonstrate that
the main
phenomena occur in the electrolyte cleft chosen as computational domain
in Sect.~\ref{sec:models}. The geometry, represented
in Fig.~\ref{fig:geom_pabst2d}, is a cross section in the $r$-$z$ plane 
of the three-dimensional computational domain.
Boundary conditions are enforced
according to the framework of Sect.~\ref{sub:Boundary-and-initial} for the
bath and for the coupling conditions, while on $\Gamma_{sim}$ 
homogeneous Neumann conditions for both potential and concentrations
are considered, as required by the axial symmetry.
The mesh used in the numerical computations is shown in 
Fig.~\ref{fig:mesh_pabst2d} and is characterized by a local refinement
in the regions close to the cell membrane. 
In the cleft between the cell and the chip (thickness 
$\delta_{j}=$ \SI{100}{\nano\meter}) 
we use a structured mesh in order to independently
control the mesh characteristic dimension 
in the $r$ and $z$ directions. 
This allows the use of very stretched triangles
to achieve a more detailed description of this area as required
by the geometrical multiscale nature of the problem in which
the ratio between cell radius and $\delta_j$ is $10^3$.
\begin{figure}[h!]
\centering
\subfigure[$\varphi$]
{\includegraphics[width=0.45\textwidth]{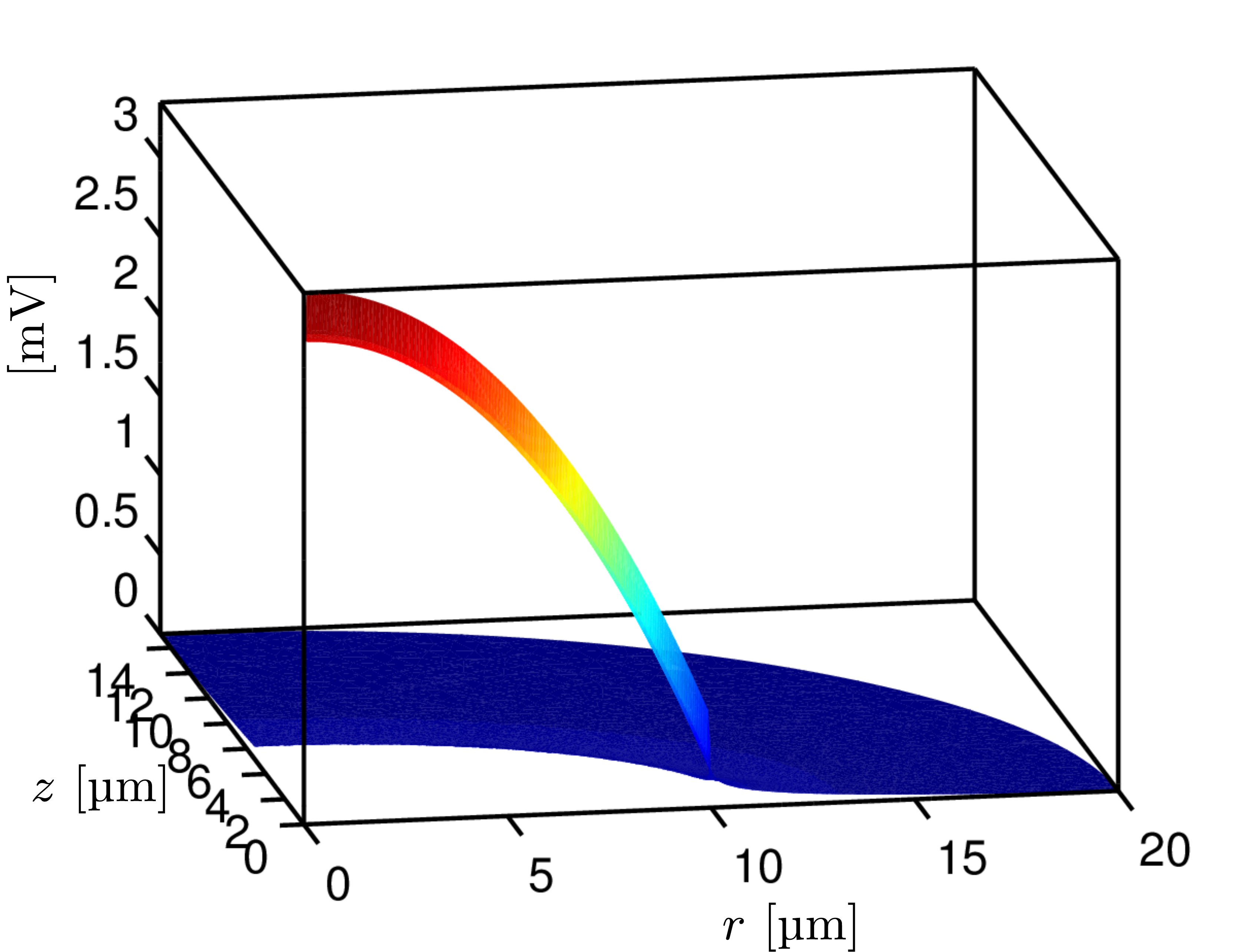}}
\hfill
\subfigure[$c_{Cl}$]
{\includegraphics[width=0.45\textwidth]{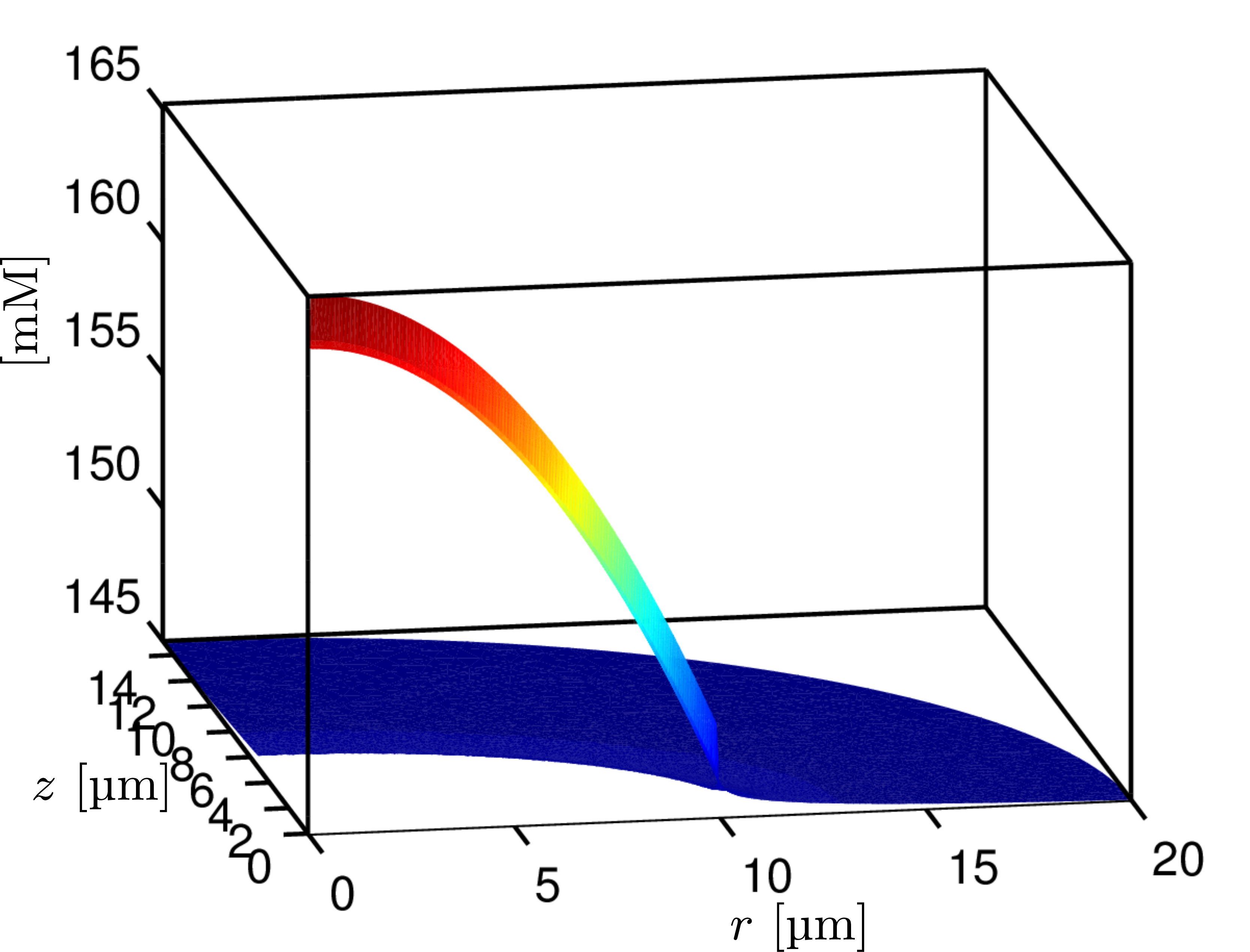}}
\subfigure[$c_{K}$]
{\includegraphics[width=0.45\textwidth]{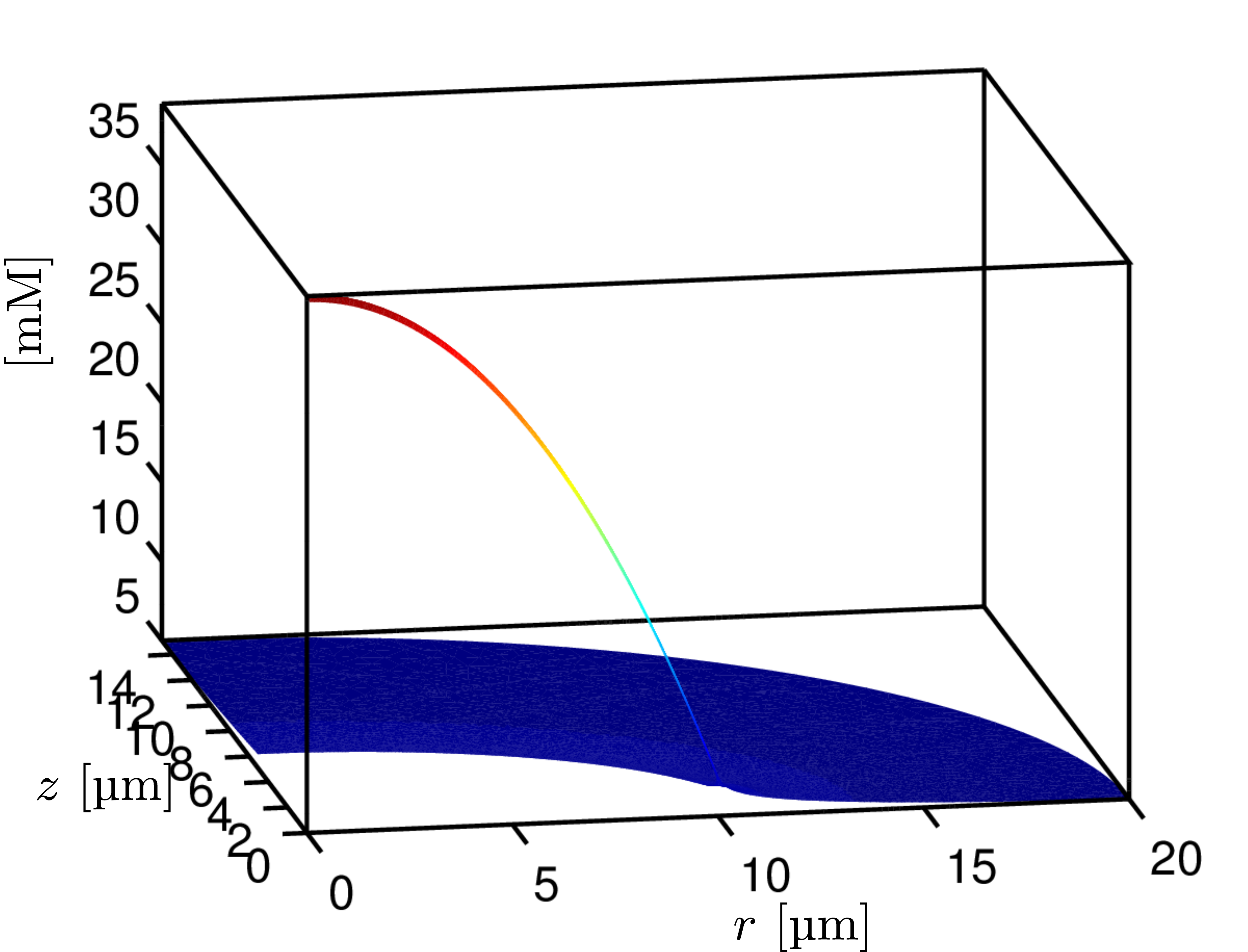}}
\hfill
\subfigure[$c_{Na}$]
{\includegraphics[width=0.45\textwidth]{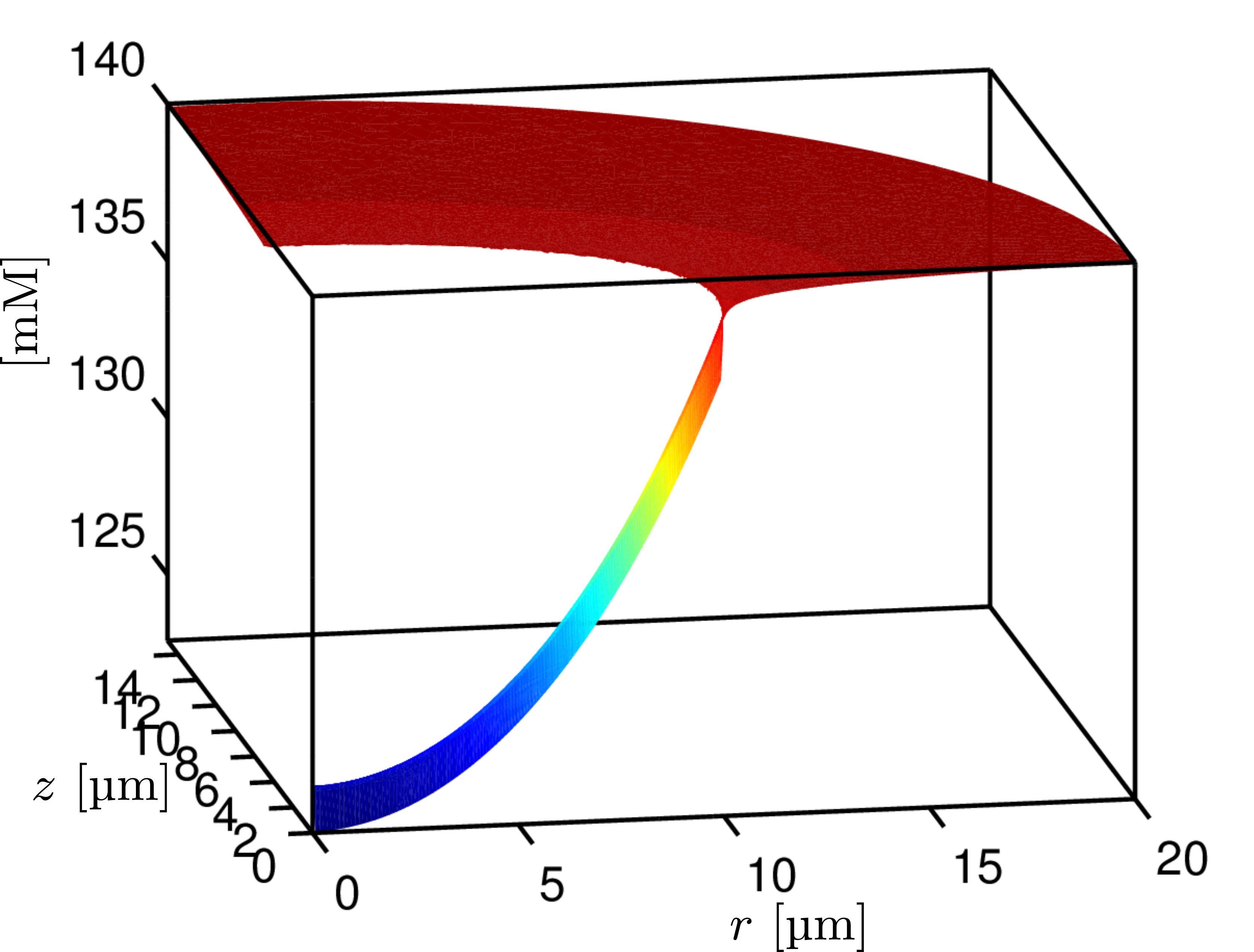}}
\caption{Spatial distribution
of the electric potential $\varphi$, and of the ion concentrations $c_{i}$, 
at the end of the transient resulting from a voltage-clamp depolarization at
$V_{cell}=$ \SI{50}{\milli\volt}.}
\label{fig:pabst_withcouplings}
\end{figure}
\begin{table}[h!]
\centering
\begin{tabular}{lll}\hline
\textbf{Parameter}		& \textbf{Symbol} 	& \textbf{Value}\\ [.25ex]
\hline
Intracellular potassium concentration    & $c_{K}^{int}$		& $140$ mM \\[.25ex]
Intracellular sodium concentration 	    & $c_{Na}^{int}$		& $4$ mM  \\[.25ex]
Intracellular chloride concentration 	& $c_{Cl}^{int}$		& $144$ mM  \\[.25ex]
Extracellular bath potassium concentration    & $c_{K}^{bath}$	& $5$ mM  \\[.25ex]
Extracellular bath sodium concentration 	    & $c_{Na}^{bath}$	& $140$ mM  \\[.25ex]
Extracellular bath chlorine concentration 	& $c_{Cl}^{bath}$	& $145$ mM \\[.25ex]
Potassium conductance   & $g_M^K$			& $250$ \si{\siemens\per\square\meter} \\[.25ex]
Membrane specific capacitance   & $C_M$      & 1 \si{\micro\farad\per\square\centi\meter} \\[.25ex]
Substrate specific capacitance  & $C_S$    & 0.3 \si{\micro\farad\per\square\centi\meter} \\[.25ex]
Initial transmembrane potential  & $V_{cell}-\left.\varphi\right|_{\Gamma_{cell}}$    & -85 \si{\milli\volt}  \\[.5ex]
\hline
\end{tabular} 
\caption{Model parameter values considered in the simulations involving passive HEK cells.}
\label{tab:Parameters_3D}
\end{table}
\begin{figure}[tb]
\centering
\subfigure[$\varphi$]
{\includegraphics[width=0.4\textwidth]{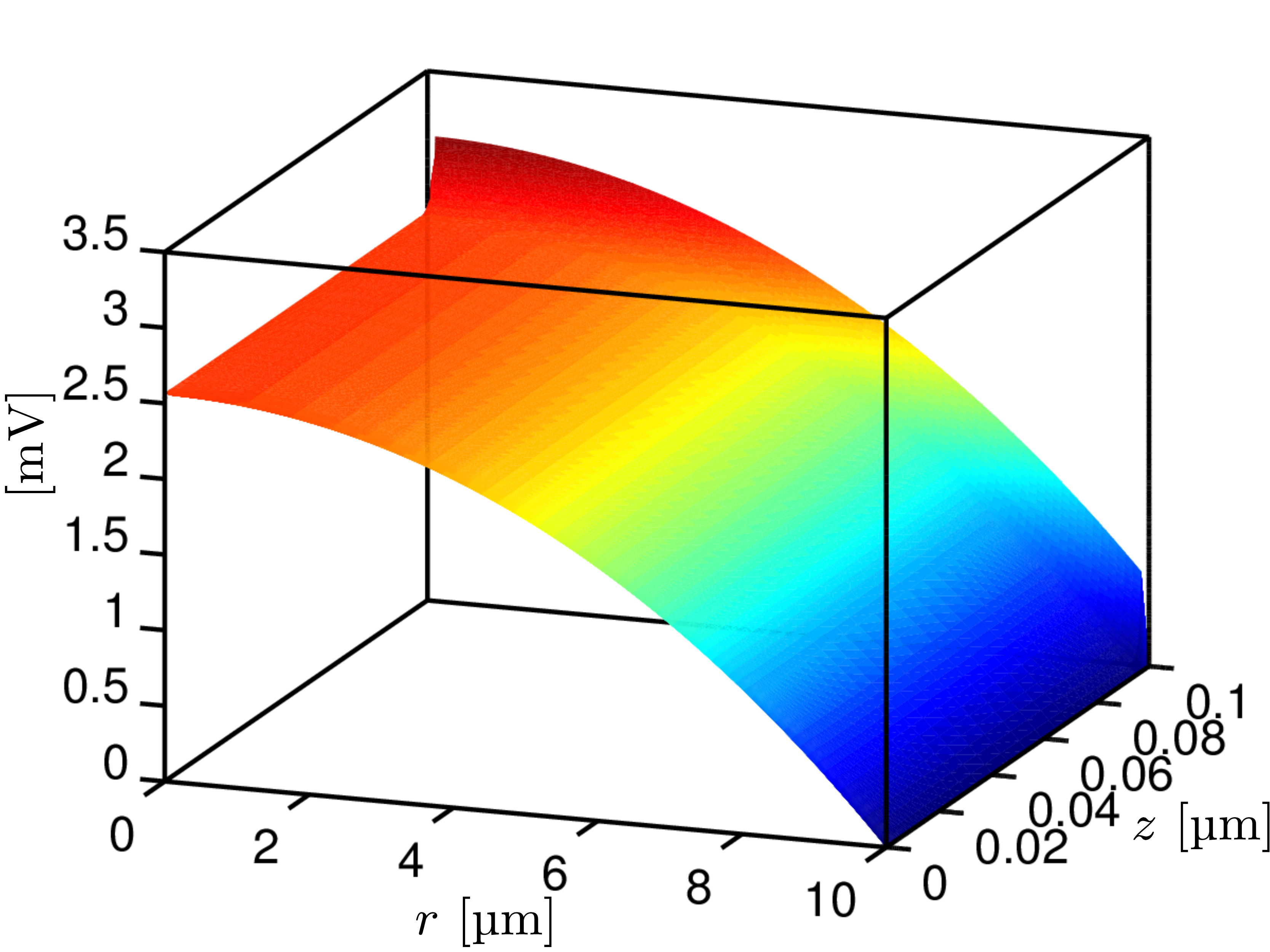}}
\hfill
\subfigure[$c_{Cl}$]
{\includegraphics[width=0.4\textwidth]{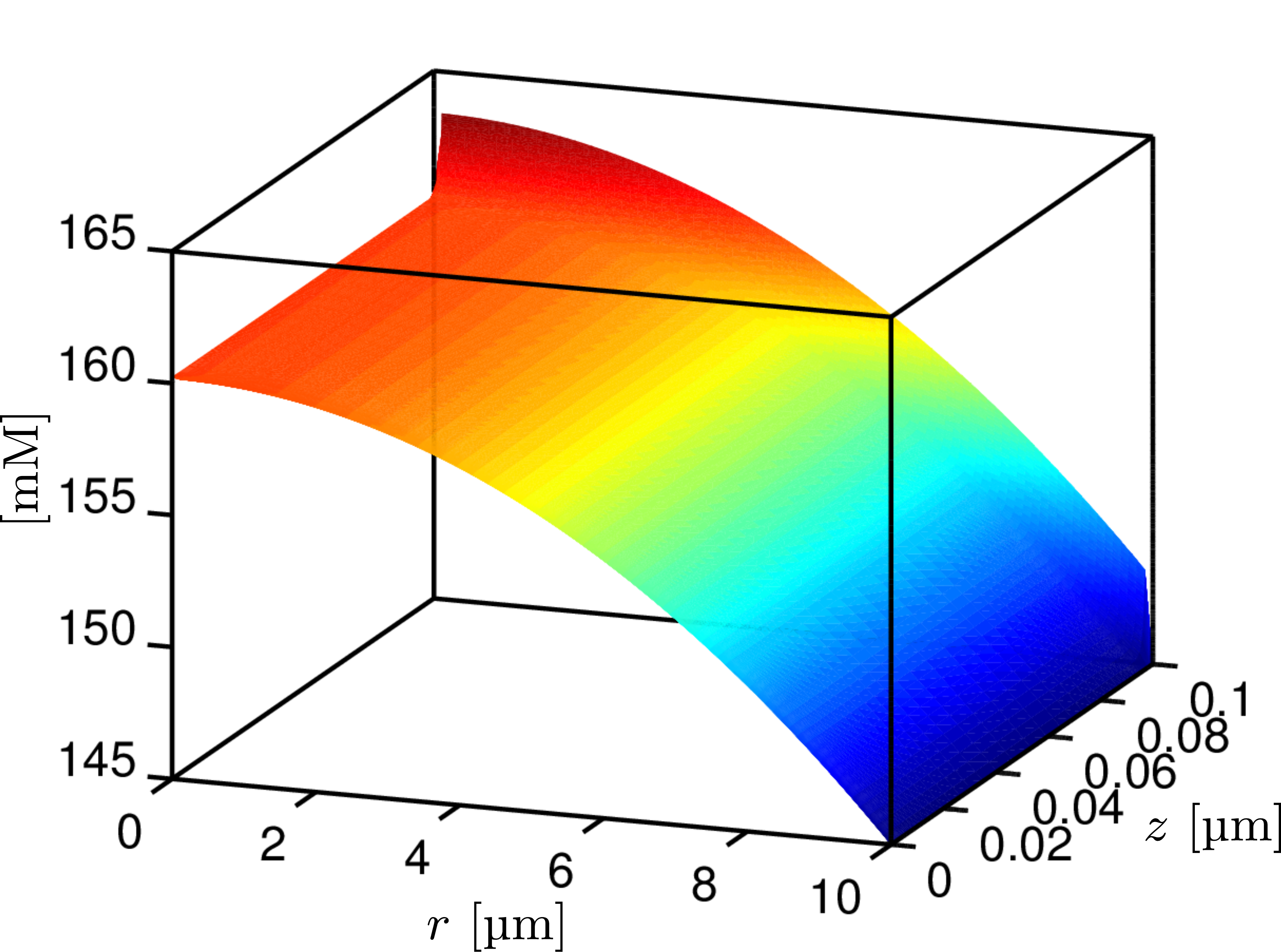}}
\subfigure[$c_{K}$]
{\includegraphics[width=0.4\textwidth]{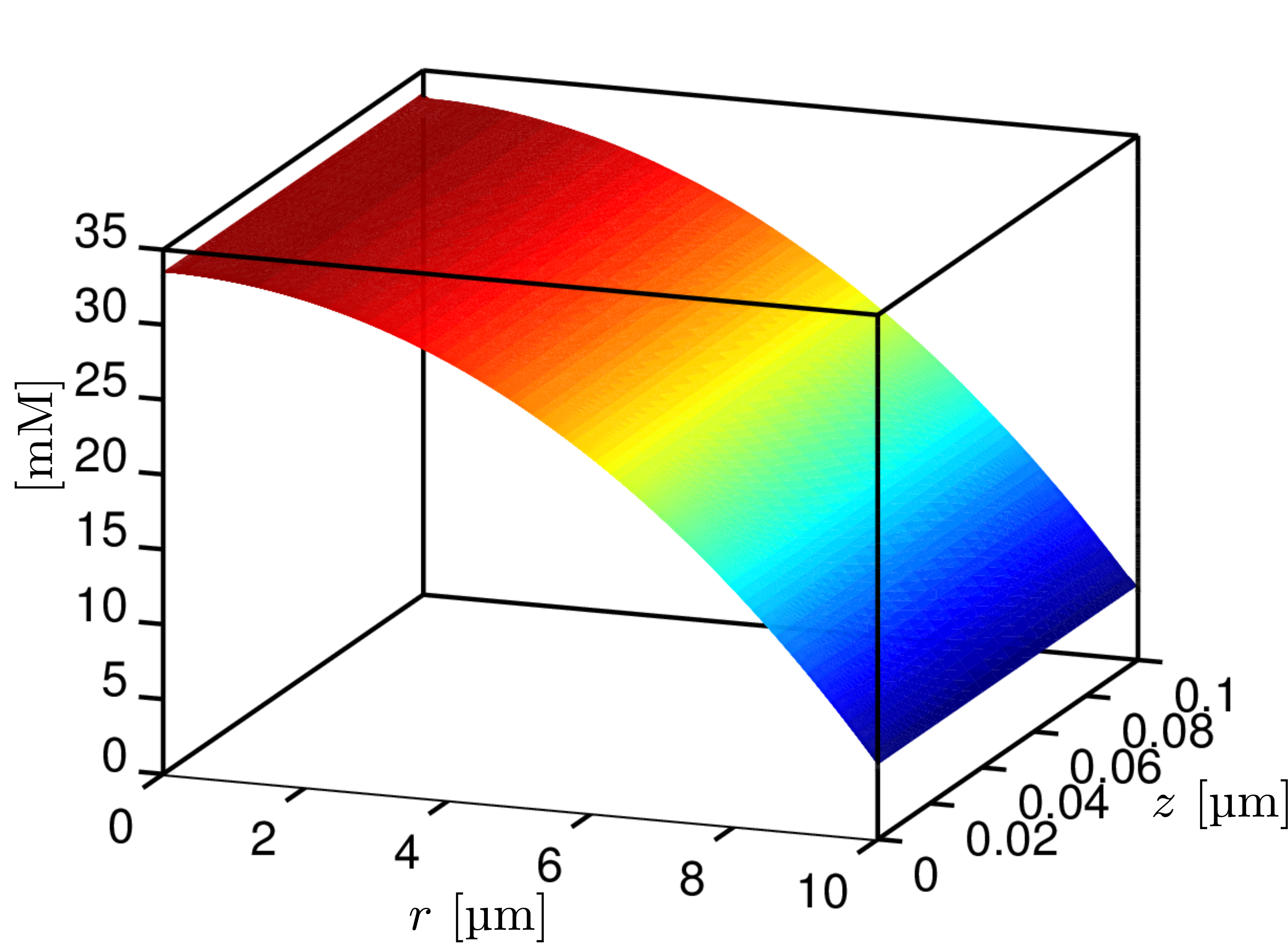}}
\hfill
\subfigure[$c_{Na}$]
{\includegraphics[width=0.4\textwidth]{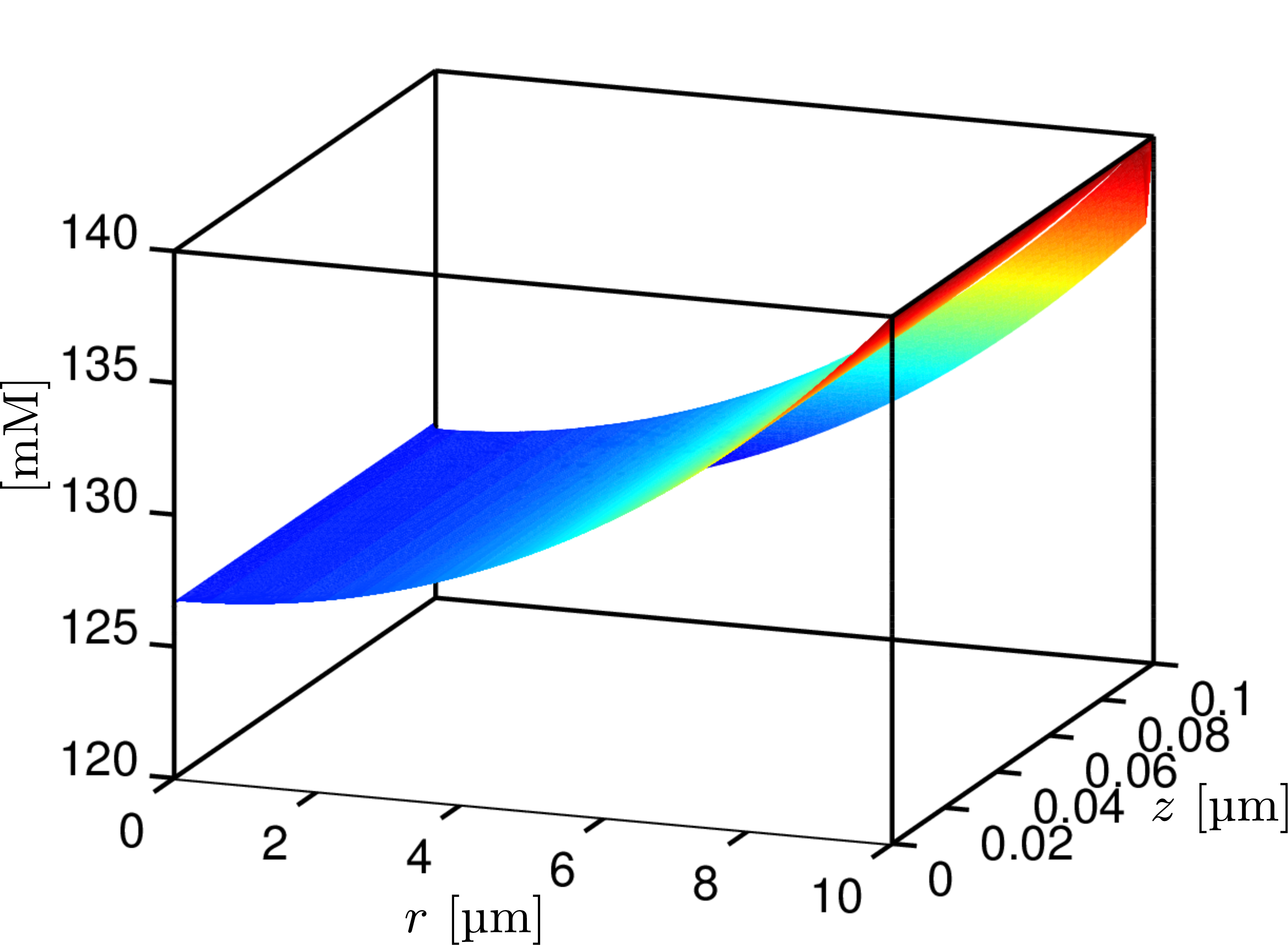}}
\caption{Spatial distribution of $\varphi$ and
$c_{i}$ in the portion
of electrolyte between cell and substrate, 
at the end of the transient
resulting from a voltage-clamp depolarization at
$V_{cell}=$ \SI{50}{\milli\volt}.
}
\label{fig:caso_brera}
\end{figure}

We consider a vol\-ta\-ge-clamp configuration in which the cell is
stimulated by varying the intracellular 
potential from \SI{-85}{\milli\volt}
to \SI{50}{\milli\volt},
and we describe the transmembrane currents
with the Goldman-Hodgkin-Katz model 
illustrated in Sect.~\ref{sub:Boundary-and-initial}.
We consider a human embryonic kidney cell HEK293, as in~\cite{Brittinger2005,Pabst2007},
expressing just potassium channels,
and the adopted parameter values are
reported in Table~\ref{tab:Parameters_3D}.

Fig.~\ref{fig:pabst_withcouplings} shows the obtained results 
in terms of potential and ion concentrations.
As a consequence of the change of the 
intracellular potential, 
the \ce{K+} channels open and \ce{K+} ions flow in the 
extracellular space.
The variation of the potassium concentration is not very 
large outside the cleft region because ions
can diffuse to the bulk of the solution and electroneutrality is 
restored in a few nanometers. In the cleft region, instead, 
ion diffusion is limited by the presence of the substrate and 
ions are forced to follow a radial pathway towards the bulk of the solution.
As a consequence, a much higher potassium concentration is determined
and consequently chloride and sodium ions are respectively attracted 
and repelled from the cleft region. The determined electric potential variation
shows a parabolic profile with a peak value of approximately \SI{3}{\milli\volt}
at the center of the cell.
In view of these results, we can therefore state that the approximations introduced in 
Sect.~\ref{sec:3Dgeometrical_model} are sound. 
Fig.~\ref{fig:caso_brera} reports the results obtained with a simulation conducted only underneath the cell where
it is possible to spot the steep layers at 
the membrane due to capacitive coupling and charge screening effects.

\bigskip

Having started our analysis from a single cell
on an electronic substrate,
we now investigate configurations with more than one
cell and/or more than one electrode, focusing on 
the mutual influence between neigbouring devices.

\subsection{Voltage-clamp stimulation: signal on a neighboring electrode}
\label{sub:Cell-to-chip}
\begin{figure}[h!]
\centering
\includegraphics[width=0.8\textwidth]{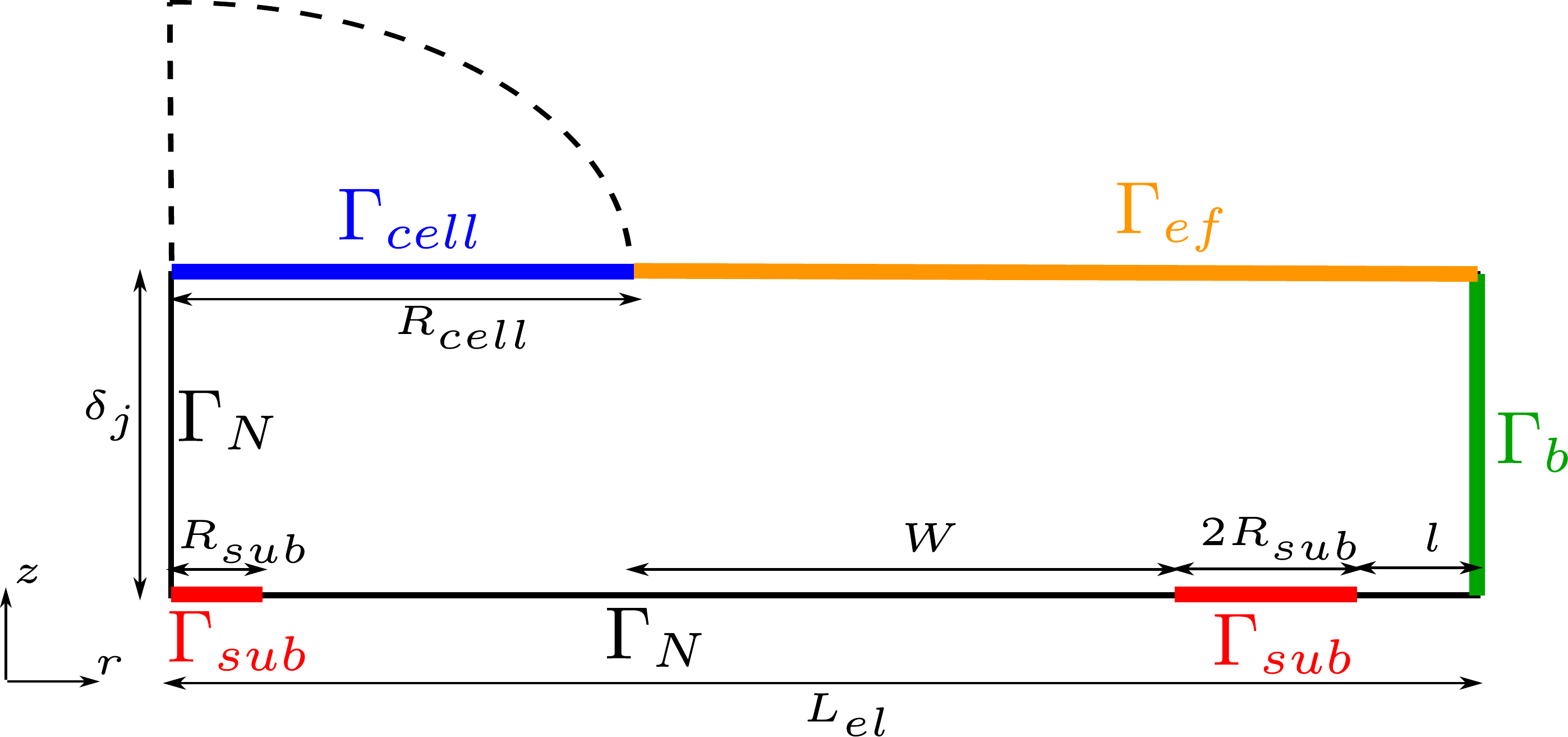}
\caption
{Computational domain with two electrodes
for a voltage-clamp stimulation recording. 
Figure not in scale: $\delta_{j} = \SI{100}{\nano\meter}, 
R_{cell} = \SI{10}{\micro\meter}, 
R_{sub} =  \SI{1.5}{\micro\meter}$ and
$l = \SI{2}{\micro\meter}$.}
\label{fig:celltochip_geom}
\end{figure}

In an ideal device for biological signal recording only the closest 
electrode to the excited cells should be recording a signal,
but this actually does not occur in realistic experiments.
In order to investigate this non ideality effect,
we consider a configuration
in which a cell is stimulated and the resulting signal is probed by two 
electrodes, one just below the cell and the other facing the bath electrolyte
at a variable distance $W$. The representation of the computational
domain with the two electrodes $\Gamma_{s1}$ and $\Gamma_{s2}$ 
and the overlying cell is shown in Fig.~\ref{fig:celltochip_geom}.
The boundary conditions on $\Gamma \setminus \Gamma_{ef}$ are the same as in Sect.~\ref{sec:voltage_clamp_validation}
and the parameter values are reported in Table~\ref{tab:Parameters_3D}.
To account for the whole electrolyte surrounding the cell, 
we enforce on $\Gamma_{ef}$ the Robin boundary condition
described in the modeling procedure of 
Sect.~\ref{sub:Boundary-and-initial}, setting
 $C^*=5\cdot10^{-3}$\si{\farad\per\meter} and $v^*=10^{-3}$\si{\meter\per\second}
in all our computations.

\begin{figure}[h!]
\centering
\subfigure[$W=$ \SI{2}{\micro\meter}]
{\includegraphics[width=0.45\textwidth]{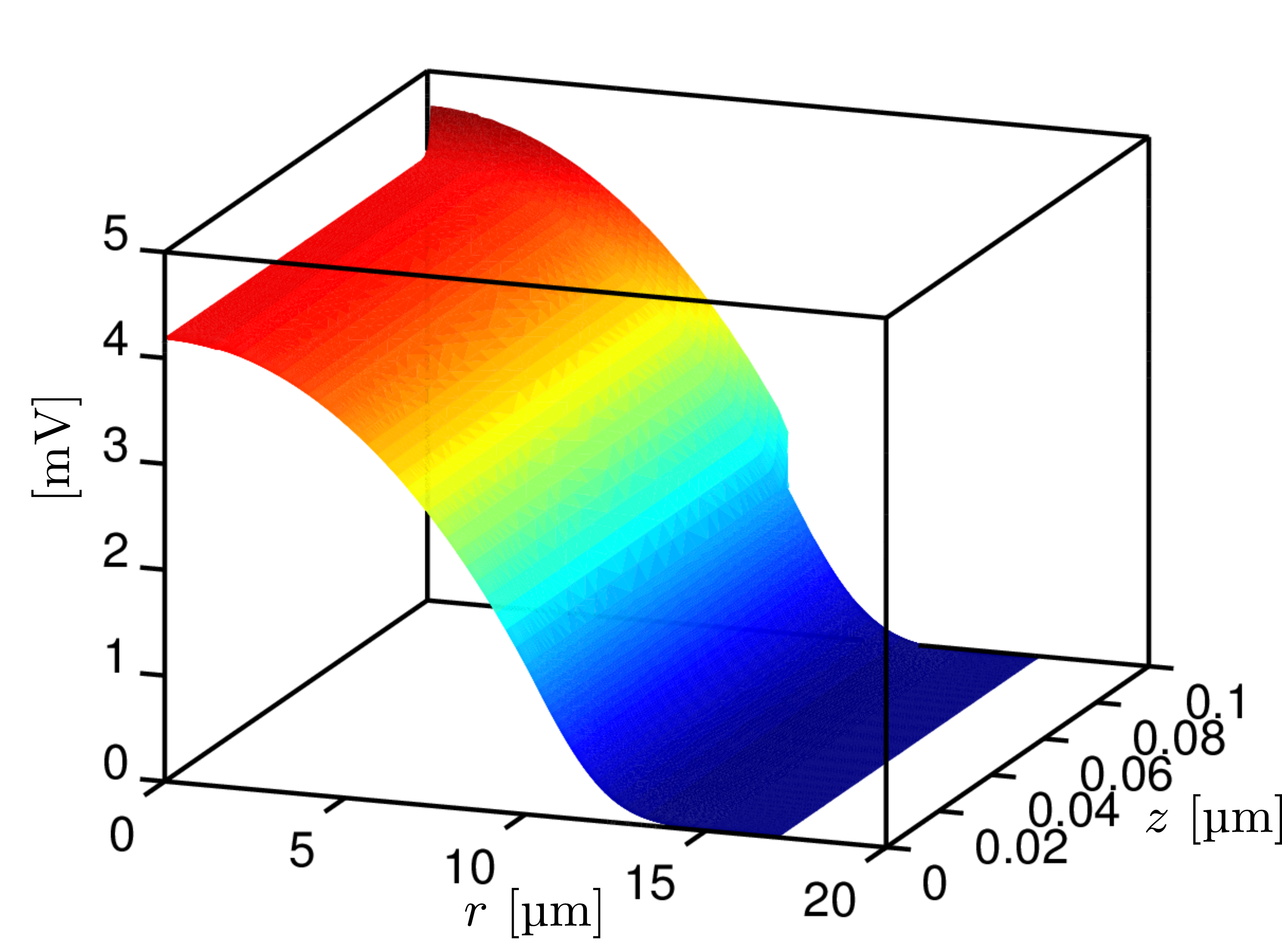}
\label{fig:pot2DW1}}
\hfill
\subfigure[$W=$ \SI{10}{\micro\meter}]
{\includegraphics[width=0.45\textwidth]{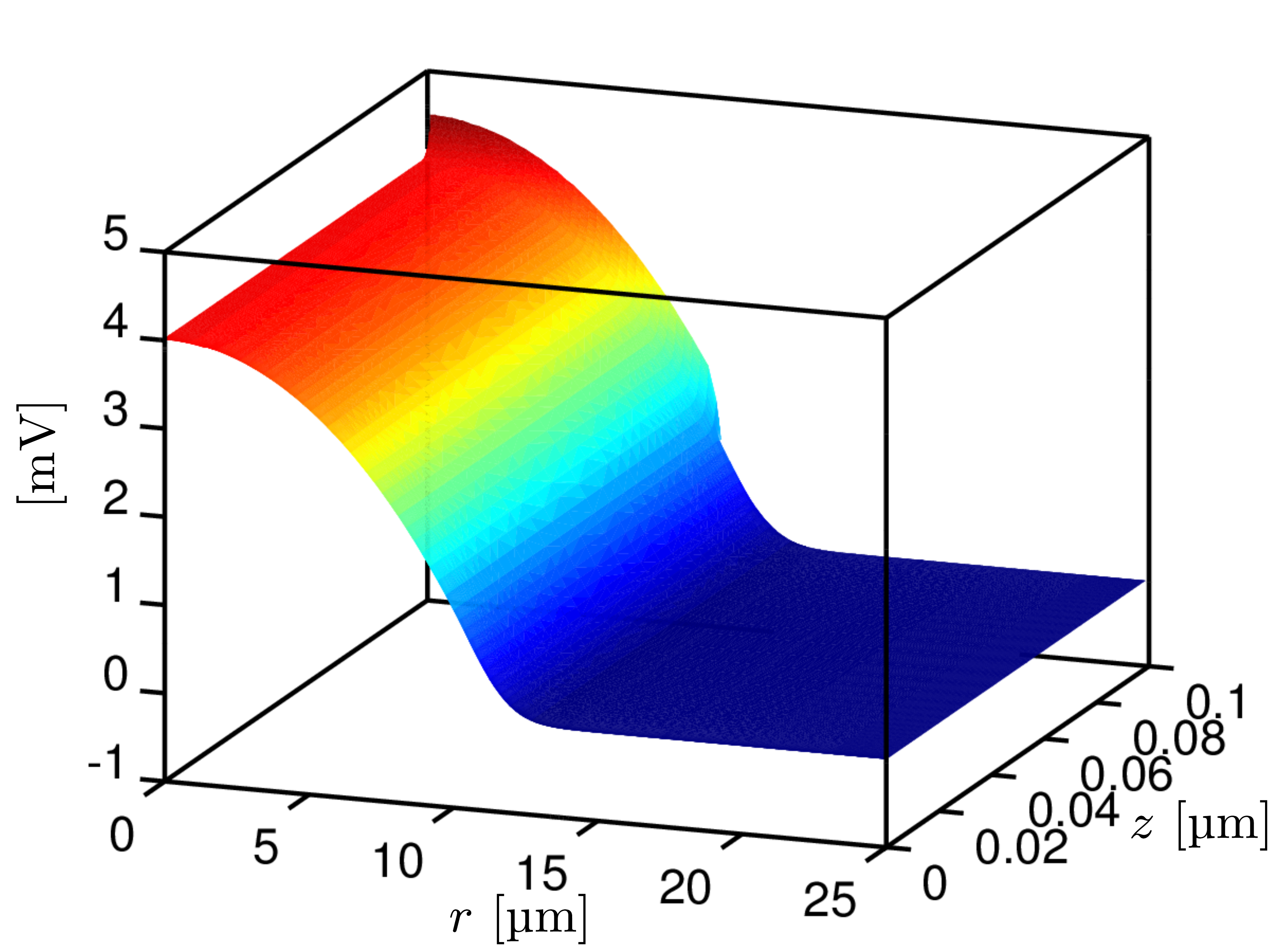}
\label{fig:pot2DW2}}
\caption{Spatial distribution
of the potential $\varphi$ with different domains: the distance between
the cell and the second gate is on the left $W=$ \SI{2}{\micro\meter} and
on the right $W=$ \SI{10}{\micro\meter}. In both cases: $\delta_{j}=$ 
\SI{100}{\nano\meter}.}
\label{fig:celltochip_potential_w2-10}
\end{figure}
The output of the simulations are the potentials measured by the two
electrodes as a consequence of the voltage-clamp stimulation.
Fig.~\ref{fig:celltochip_potential_w2-10} reports the computed 
potential profile in the cases $W = \SI{2}{\micro\meter}$ and 
$W = \SI{10}{\micro\meter}$. Far from the cell adhesion region
$r > \SI{10}{\micro\meter}$, the potential decays to the bulk value ($\varphi=0$) with a trend proportional to $1/r$. 
In the first device configuration the second electrode is located in the region
with $r \in [12,15]\,\si{\micro\meter}$
and the probed electric potential 
is significantly different from the bulk value.
In the second configuration, instead, the electrode is placed between 
$[20,23]\,\si{\micro\meter}$, where the
perturbation is almost equilibrated, so the measured signal is very low.

\begin{figure}[h!]
\centering
\subfigure[Cl]{\includegraphics[width=0.45\textwidth]{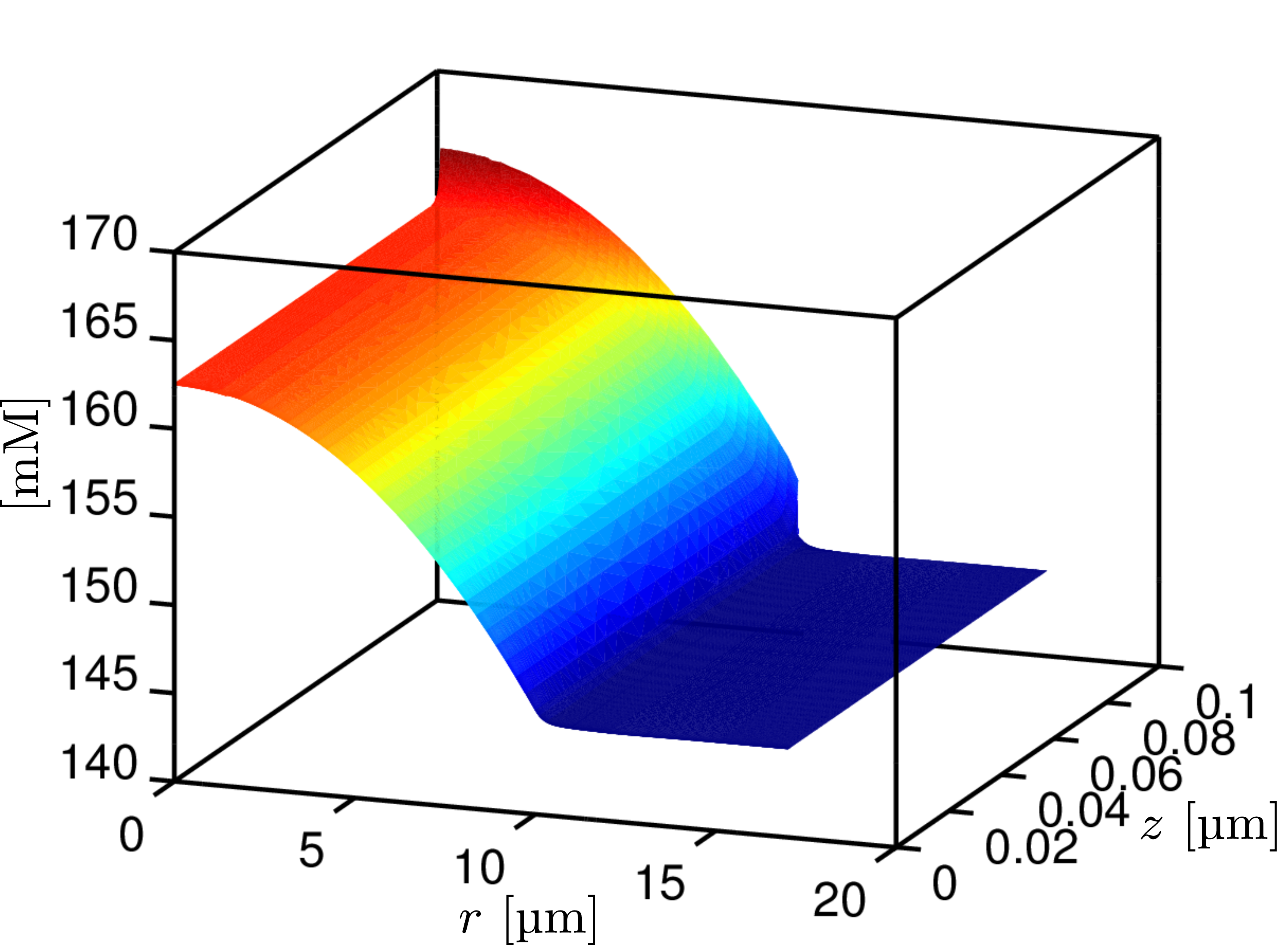}
\label{fig:Cl2DW1}} 
\hfill
\subfigure[K]{\includegraphics[width=0.42\textwidth]{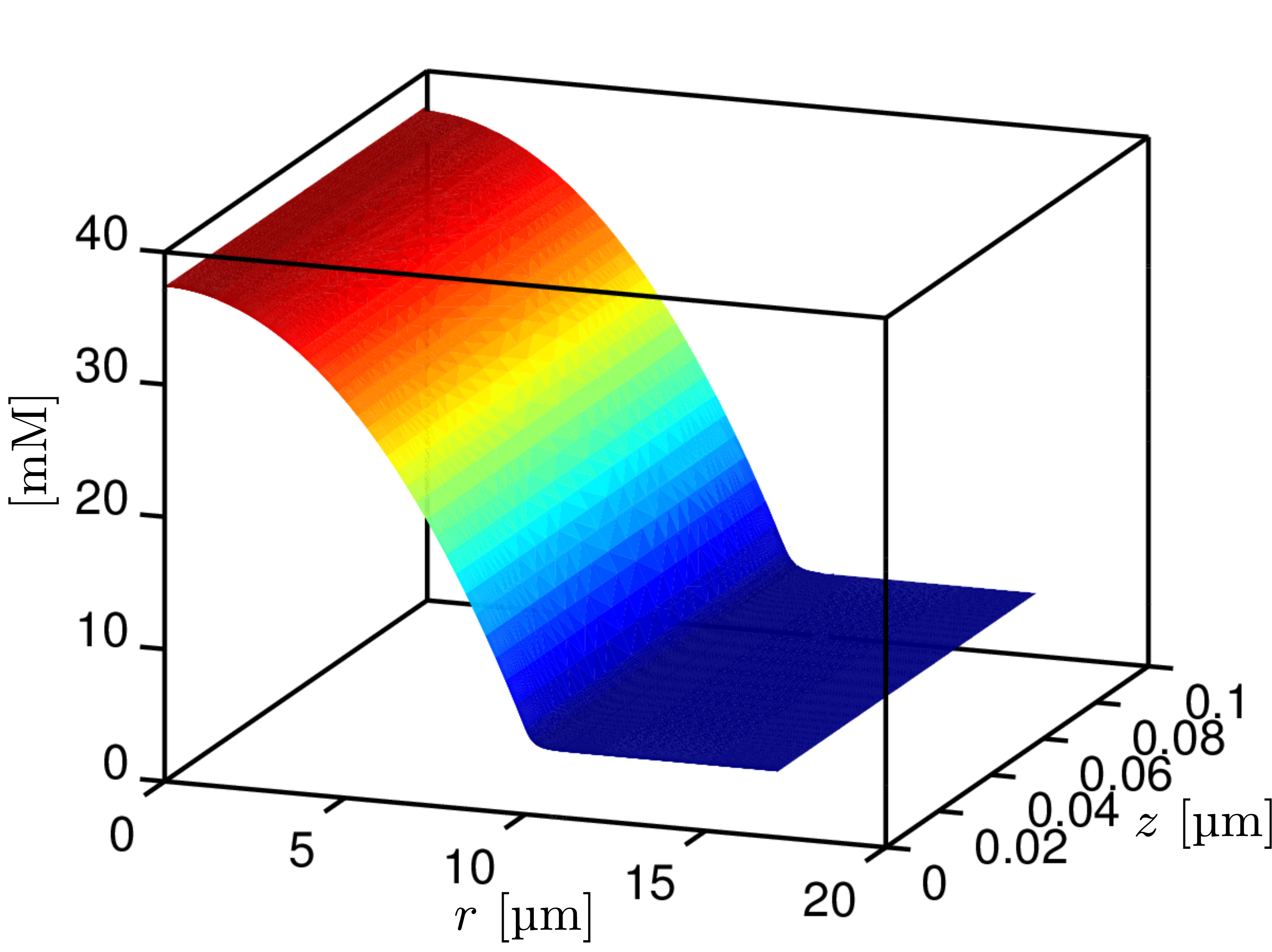}
\label{fig:K2DW1}} 
\hfill
\subfigure[Na]{\includegraphics[width=0.45\textwidth]{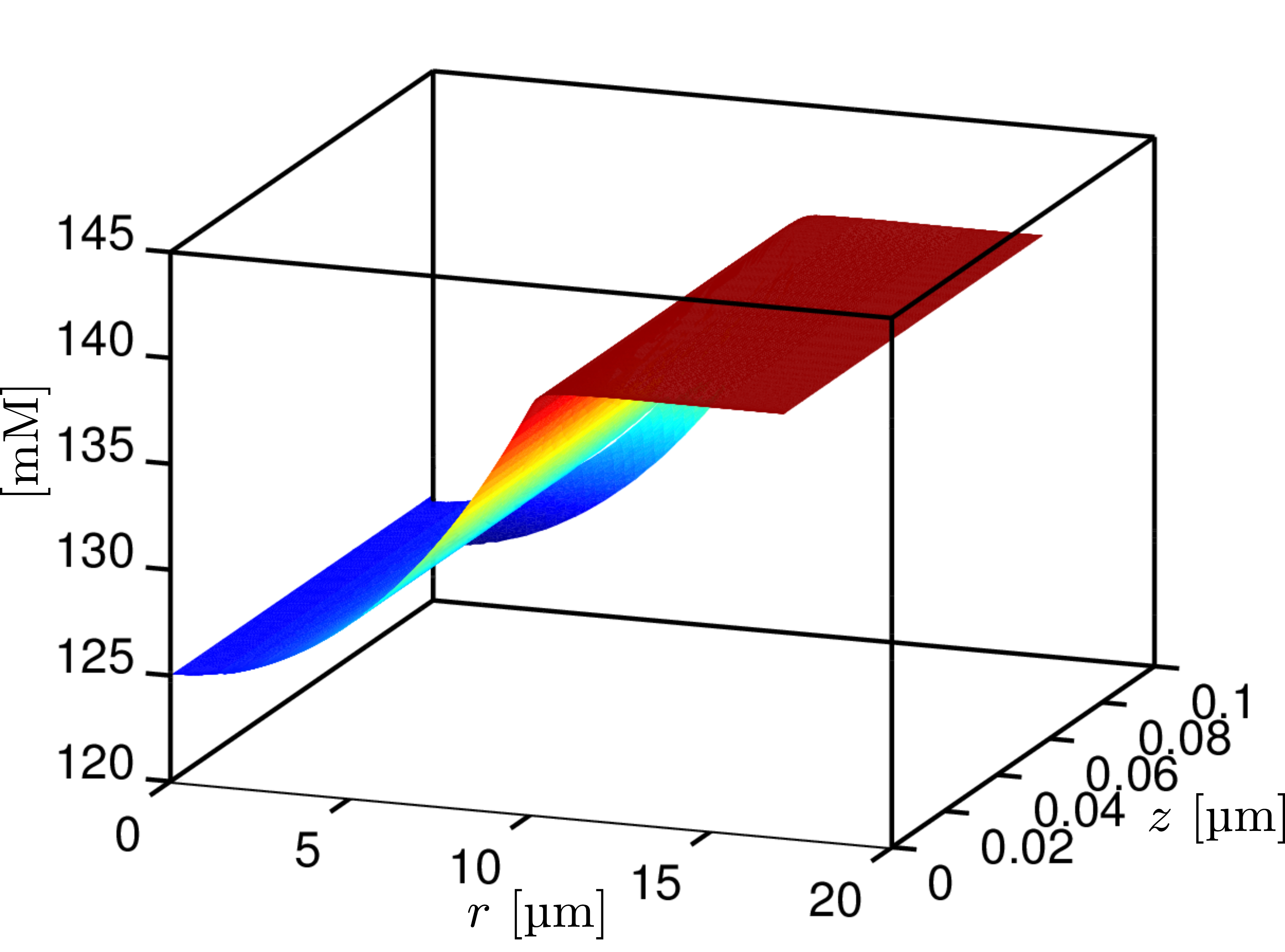}
\label{fig:Na2DW1}}
\caption{Spatial distribution of the computed concentrations 
$c_{i}$ with $W = \SI{2}{\micro\meter}$.}
\label{fig:celltochip_concentr}
\end{figure}

The computed concentration profiles are shown in 
Fig.~\ref{fig:celltochip_concentr}, and 
notably the profiles are in very good agreement with the results of
Sect.~\ref{sec:voltage_clamp_validation}.

\begin{figure}[tb]
\centering
\includegraphics[width=0.5\textwidth]{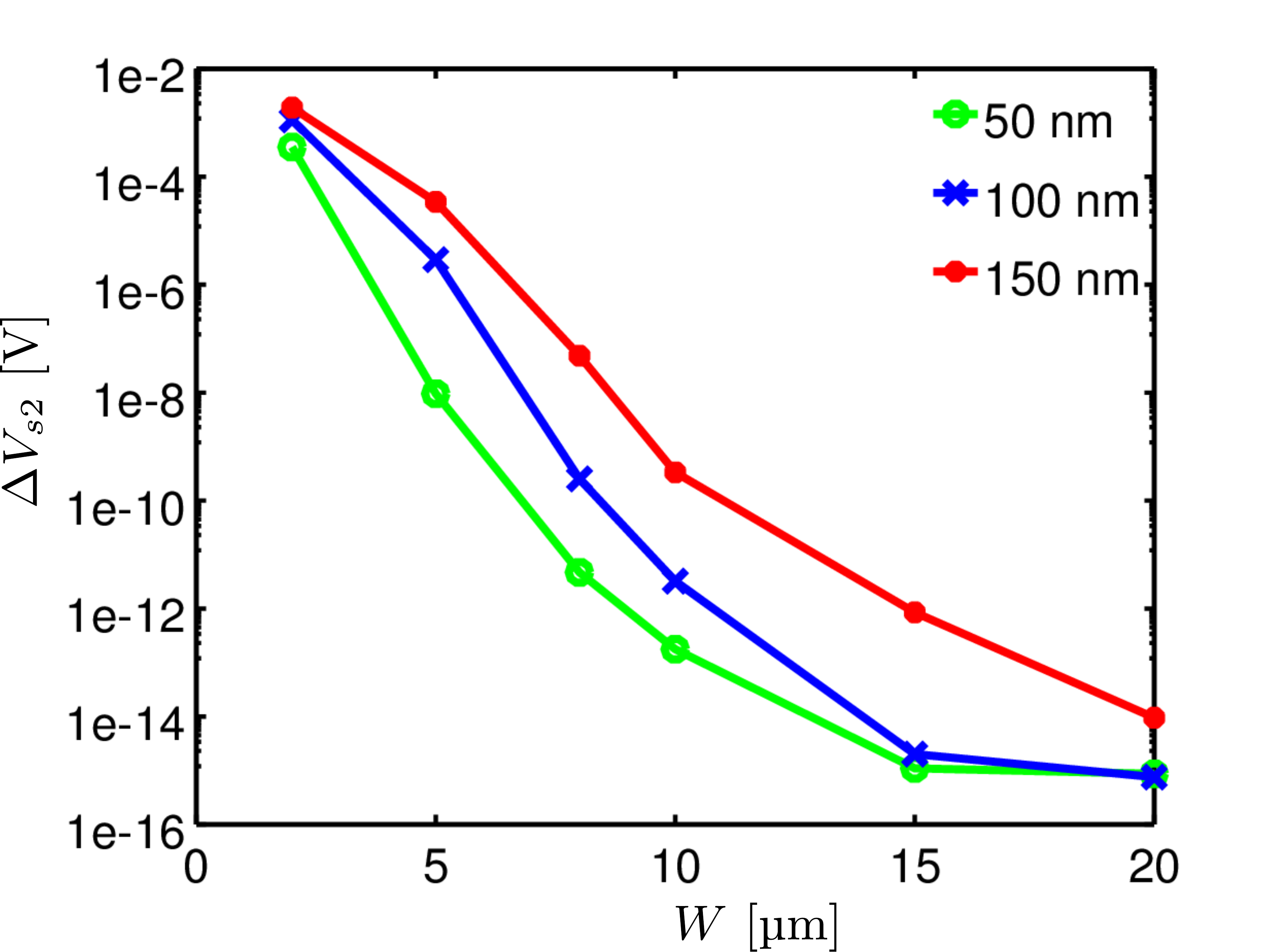}
\caption{Probed voltage $\Delta V_{s2}$ as a function of $W$.
Results for three different values of the cleft thickness $\delta_{j}$,
obtained with a depolarizing pulse with $V_{cell}=$ \SI{50}{\milli\volt}.}
\label{fig:celltochip_param_delta}
\end{figure}

\begin{figure}[tb]
\centering 
\subfigure[$\delta_j=$ \SI{50}{nm}]
{\includegraphics[width=0.4\textwidth]{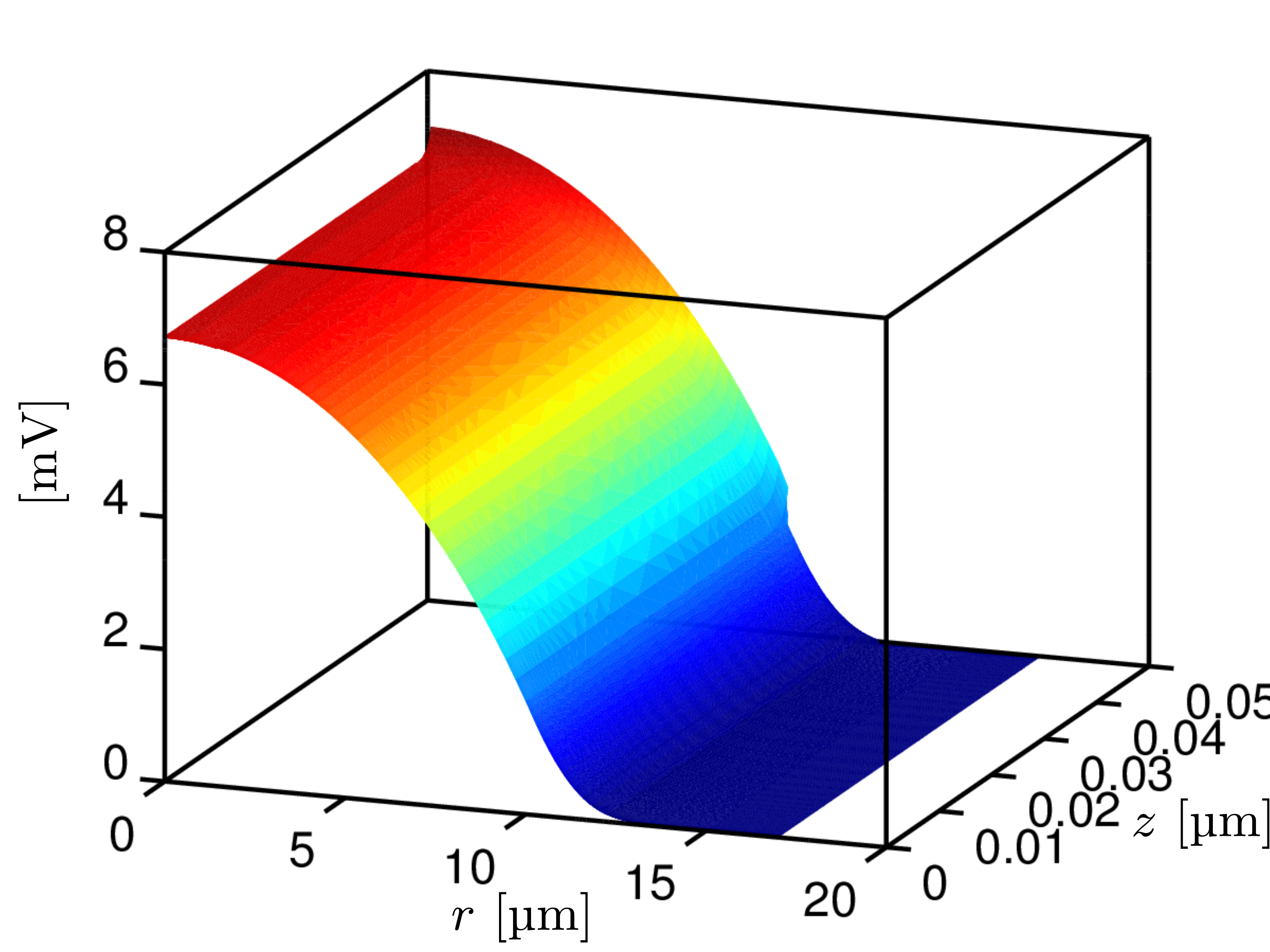}
\label{fig:pot_delta1}} 
\hfill
\subfigure[$\delta_j=$ \SI{150}{nm}]
{\includegraphics[width=0.4\textwidth]{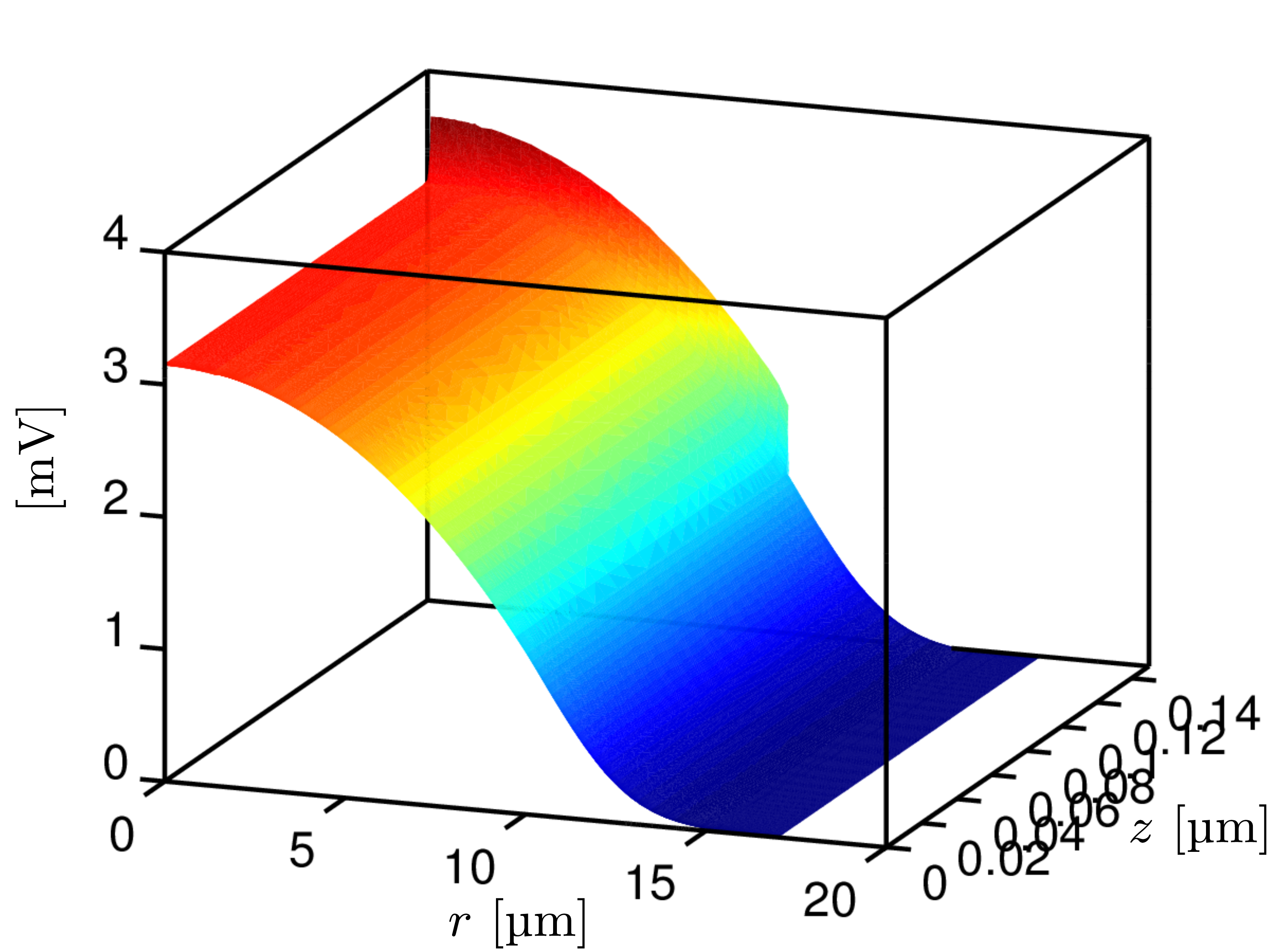}
\label{fig:pot_delta2}} 
\caption
{Spatial distribution
of the potential $\varphi$ with different domains: the distance between
the cell and the second gate is set at $W=$ \SI{2}{\micro\meter}, but on
the left we have a cleft thickness $\delta_{j}=$ \SI{50}{\nano\meter} and
on the right $\delta_{j}=$ \SI{150}{\nano\meter}.}
\label{fig:celltochip_potential_delta}
\end{figure}

Since the potential $\varphi$ in proximity of the second 
electrode does not attain a uniform value, we introduce the local average
\begin{equation*}
\Delta V_{s2} = 
\dfrac{1}{\left|\Gamma_{s2}\right|}
\int_{\Gamma_{s2}}\left(\left.\varphi\right|_{\Gamma_{s2}}-V_{s2}\right)d\gamma
\end{equation*}
to be interpreted as the value probed by the second electrode and returned as
output.
We consider values of $W$ in the range between 2 and \SI{20}{\micro\meter}
and three different values of the cleft thickness $\delta_{j}$, 
namely 50, 100 and \SI{150}{\nano\meter}, 
In Fig.~\ref{fig:celltochip_param_delta} we report the computed
values of $\Delta V_{s2}$ and
we observe an almost exponential decrease 
of the signal when considering an increasingly distant electrode. 
Moreover, we observe in Fig.~\ref{fig:celltochip_potential_delta} that
in configurations characterized by a smaller value of $\delta_{j}$, 
a more intense variation of the 
potential is registered in the portion of electrolyte under 
the cell, due to the fact 
that higher values of the potassium concentration occur there.
From Fig.~\ref{fig:celltochip_param_delta}
we also notice that the value of $\Delta V_{s2}$ decreases with the cleft thickness. 
A physical explanation for
this latter result can be provided by resorting to the definition of electrical
resistance for the cleft, which applies here since
the electrolyte solution is an electrical conductor. We have
\begin{equation*}
R_{el} = \rho_{el}\dfrac{L_{el}}{S_{el}},
\end{equation*}
where $\rho_{el}$ is the electrolyte resistivity, $L_{el}$ 
and $S_{el}$ are the 
length and the cross sectional area of the cleft, 
this latter being
linearly proportional to the thickness $\delta_{j}$. 
A smaller value of $\delta_{j}$ hence results in
a larger
resistance $R_{el}$, 
determining a more pronounced decay of the potential
along the cleft radius, as shown in Fig.~\ref{fig:pot_delta1}. 

\subsection{Voltage-clamp stimulation: effect on a neighboring cell}
\label{sec:cell_cell}
\begin{figure}[h!]
\centering
\includegraphics[width=0.75\textwidth]{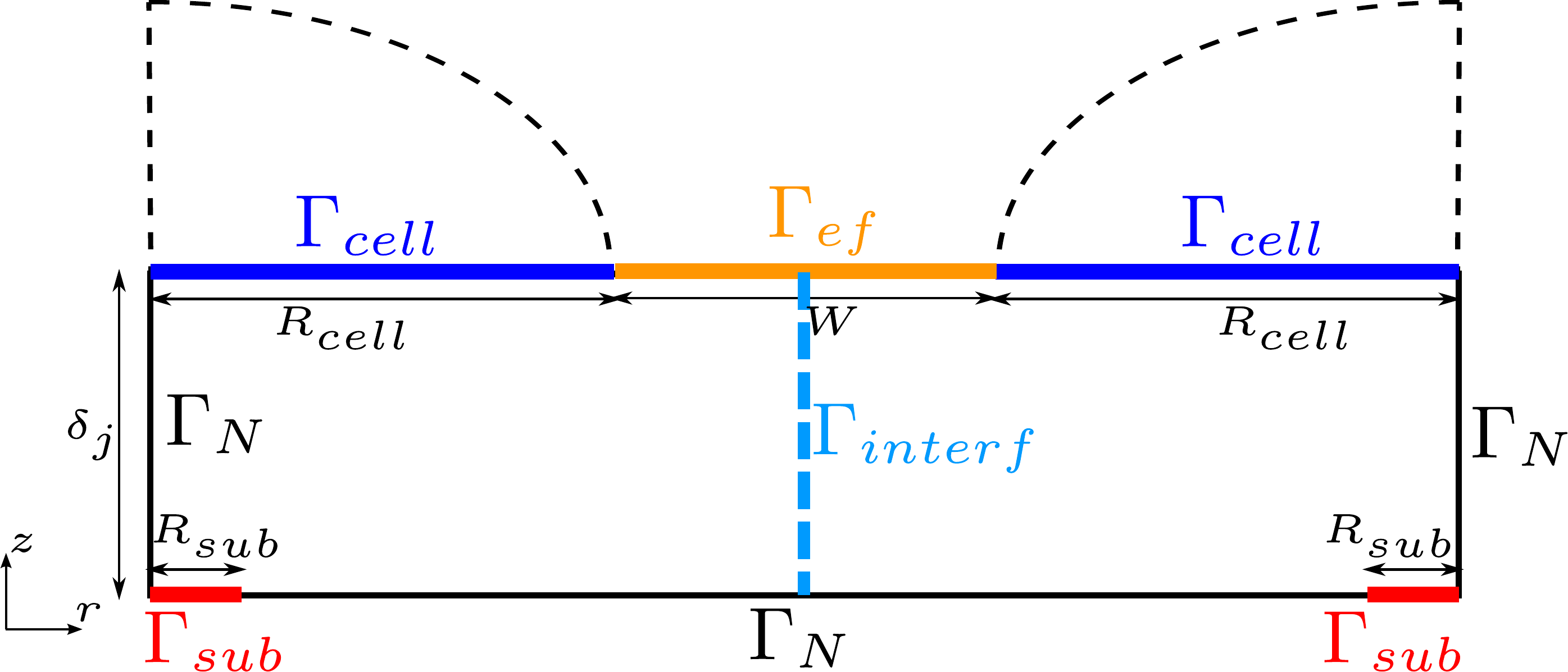}
\caption{Computational domain with two cells
and two electrodes. Figure not in scale:
$\delta_{j}= \SI{100}{nm}, R_{cell}= \SI{10}{\micro\meter}$ and
$R_{sub}= \SI{1.5}{\micro\meter})$.
}
\label{fig:celltocell_geom}
\end{figure}

When more cells are attached to a bioelectronic device, if 
a stimulation is applied to just one cell, the perturbation can be sensed
also by the neighboring ones.
Here we investigate this effect by studying the device configuration 
of Fig.~\ref{fig:celltocell_geom}, where a second cell on the right
(denoted from now on as cell B) is 
located at a distance $W$ from the stimulated one on the left
(denoted from now on as cell A).

The scheme of Fig.~\ref{fig:celltocell_geom} 
is not characterized by any rotational symmetry,
so in principle it cannot be described using the geometrical framework introduced in 
Sect.~\ref{sub:finite_el} and adopted in the previously presented
numerical results. Our approach consists of 
dividing the computational domain of  
Fig.~\ref{fig:celltocell_geom} into two parts along 
the artificial interface $\Gamma_{interf}$ and assuming 
axial symmetry to hold separately on each of the two 
subdomains. This allows us to formulate the mathematical 
problem in each subdomain using cylindrical coordinates as in 
Sect.~\ref{sub:finite_el} so that the solution of the 
problem in the whole domain is obtained through subdomain coupling 
across $\Gamma_{interf}$ using the substructuring techniques described in~\cite{quarteroni1999domain,defalco2012,abbate2014,porro2009solar}.
Since the perturbation due to cellular stimulation 
is not symmetric with respect to the axis of cell B, the above described approach introduces a certain 
level of approximation, which is expected to become less
significant as the distance between the two cells is increased.

We study an electrophysiological 
configuration in which the internal potential of both cells 
is controlled and while cell A is stimulated applying a voltage step
to \SI{50}{\milli\volt}, cell B is kept at the resting value and
the membrane current is measured.
Since the configuration differs from that considered in 
Sects.~\ref{sec:voltage_clamp_validation} and~\ref{sub:Cell-to-chip},
the resting voltage $V_{ceq}$ at which the cells are in equilibrium, 
(i.e., when the overall current flowing through the membrane is zero),
is not expected to coincide 
with the value of \SI{-85}{\milli\volt} attained in the case of
the single cell configuration.
The new value of $V_{ceq}$ can be determined by solving
the PNP system subject to 
the condition $V_{c1}=V_{c2}=V_{ceq}$, $V_{ceq}$ being 
an unknown quantity, and to the following constraints
\begin{gather*}
\dfrac{1}{\left|\Gamma_{c1}\right|}
\int_{\Gamma_{c1}}
j^{tm}(V_{c1}, \varphi\vert_{\Gamma_{c1}})
\, d\gamma
=
\dfrac{1}{\left|\Gamma_{c2}\right|}
\int_{\Gamma_{c2}}
j^{tm}(V_{c2}, \varphi\vert_{\Gamma_{c2}})
\, d\gamma = 0.
\end{gather*}
The obtained value for $V_{ceq}$ is about \SI{-90}{\milli\volt}, 
slightly more negative than the previous value 
of \SI{-85}{\milli\volt}. This is probably due to the fact that 
in the case of the two-cell configuration, potassium concentration
in the cleft region is, on average, larger
than in the case of the sole cell A, because of potassium current injection
also from cell B. As a consequence, 
a larger diffusive flux tends to drive positive ions from the
electrolyte cleft towards the intracellular sites so that
a more negative clamp voltage is needed to increase 
the transmembrane electric field required to counterbalance the
increased diffusive flux and restore the equilibrium 
condition of zero transmembrane current flow. 

\begin{figure}[tb]
\centering 
\subfigure[$\varphi$ with $W = \SI{20}{\micro\meter}$]
{\includegraphics[width=0.45\textwidth]{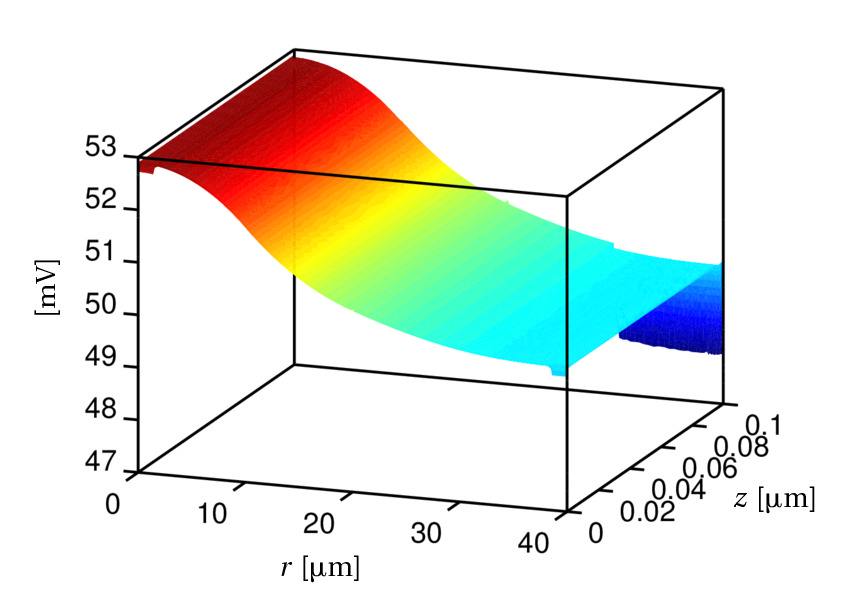}
\label{fig:pot_2cells_W1}} 
\hfill
\subfigure[$c_K\text{ with }W=$ \SI{20}{\micro\meter}]
{\includegraphics[width=0.5\textwidth]{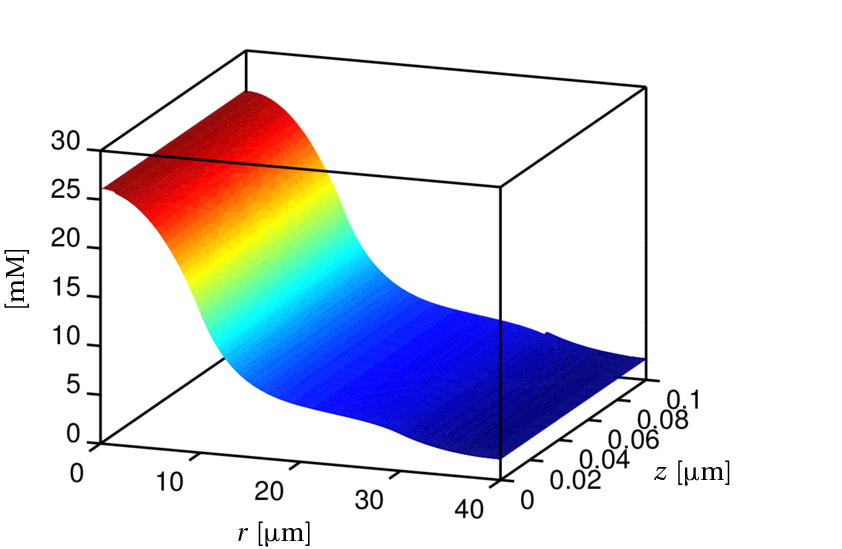}
\label{fig:cK_2cells_W1}} 
\subfigure[$\varphi\text{ with }W=$ \SI{70}{\micro\meter}]
{\includegraphics[width=0.45\textwidth]{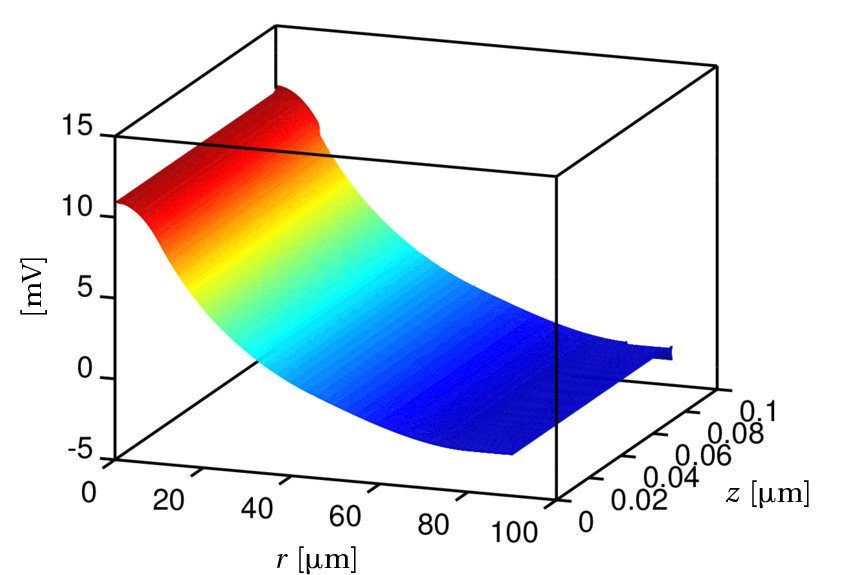}
\label{fig:pot_2cells_W2}} 
\hfill
\subfigure[$c_K\text{ with }W=$ \SI{70}{\micro\meter}]
{\includegraphics[width=0.48\textwidth]{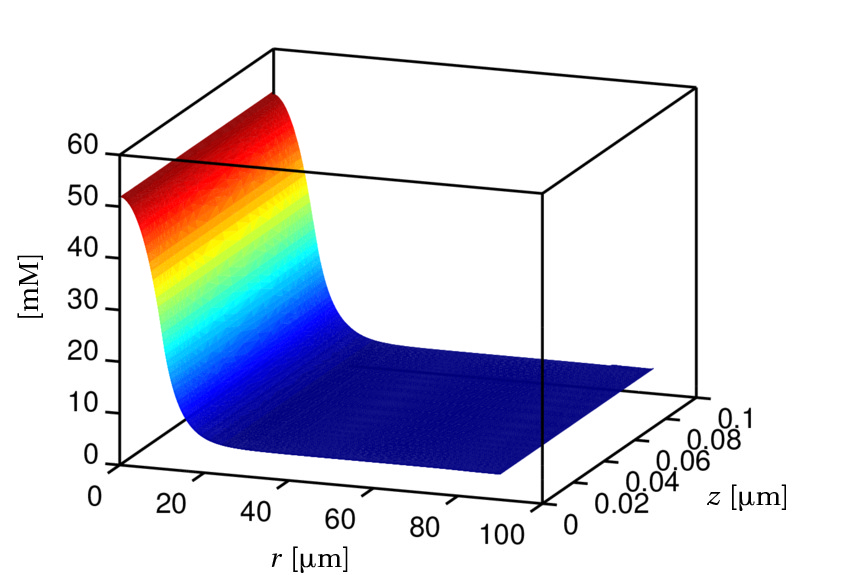}
\label{fig:cK_2cells_W2}} 
\caption{Spatial distributions of $\varphi$ and of $c_{K}$ 
at the end of the transient resulting from a voltage clamp 
stimulation of cell A to $V_{cell}=$ \SI{50}{\micro\volt}, 
keeping cell B at $V_{ceq}$.
Results for two different
distances between the cells: $W = \SI{20}{\micro\meter}$ 
and $W = \SI{70}{\micro\meter}$.}
\label{fig:celltocell_k_phi}
\end{figure}
Fig.~\ref{fig:celltocell_k_phi} shows the spatial distributions of
the potential $\varphi$ and of the potassium concentration for two
different values of the distance $W$ between the two cells.
The channels of cell A are always open and are injecting a
$\text{K}^{+}$ current in the electrolyte, because of depolarization.
This causes an increase of \ce{K+} and of $\varphi$
in the considered domain, which may lead, in turn, to the opening
of the channels of cell B. 
In the case where the two cells are close enough 
(at a distance $W = \SI{20}{\micro\meter}$) the electric potential 
exhibits a significant
boundary layer at $\Gamma_{c2}$ as shown by Fig.~\ref{fig:pot_2cells_W1}, to which corresponds an opening of the 
\ce{K+} channels of cell B. In Fig.~\ref{fig:cK_2cells_W1},
 we observe
an evident depletion in the spatial distribution $c_{K}$
under $\,\Gamma_{c2}$. This is due to the fact that the potassium
current here is entering into cell B: as physically expected
the potassium is injected by one cell (cell A) 
and collected from the other one (cell B), 
in a manner that resembles, using an electronics analogy, 
the working principle of a solid-state transistor 
(see the series of articles~\cite{BobWall,BobStruct,BobAlive}).
In the case of a larger $W$ (Figs.~\ref{fig:pot_2cells_W2} 
and~\ref{fig:cK_2cells_W2}),
the value of the potential in the electrolyte is lower and there is
practically no current entering into cell B, because 
in this configuration ions are free to flow in a larger portion of electrolyte.

\begin{figure}[h!]
\centering
\includegraphics[width=0.55\textwidth]{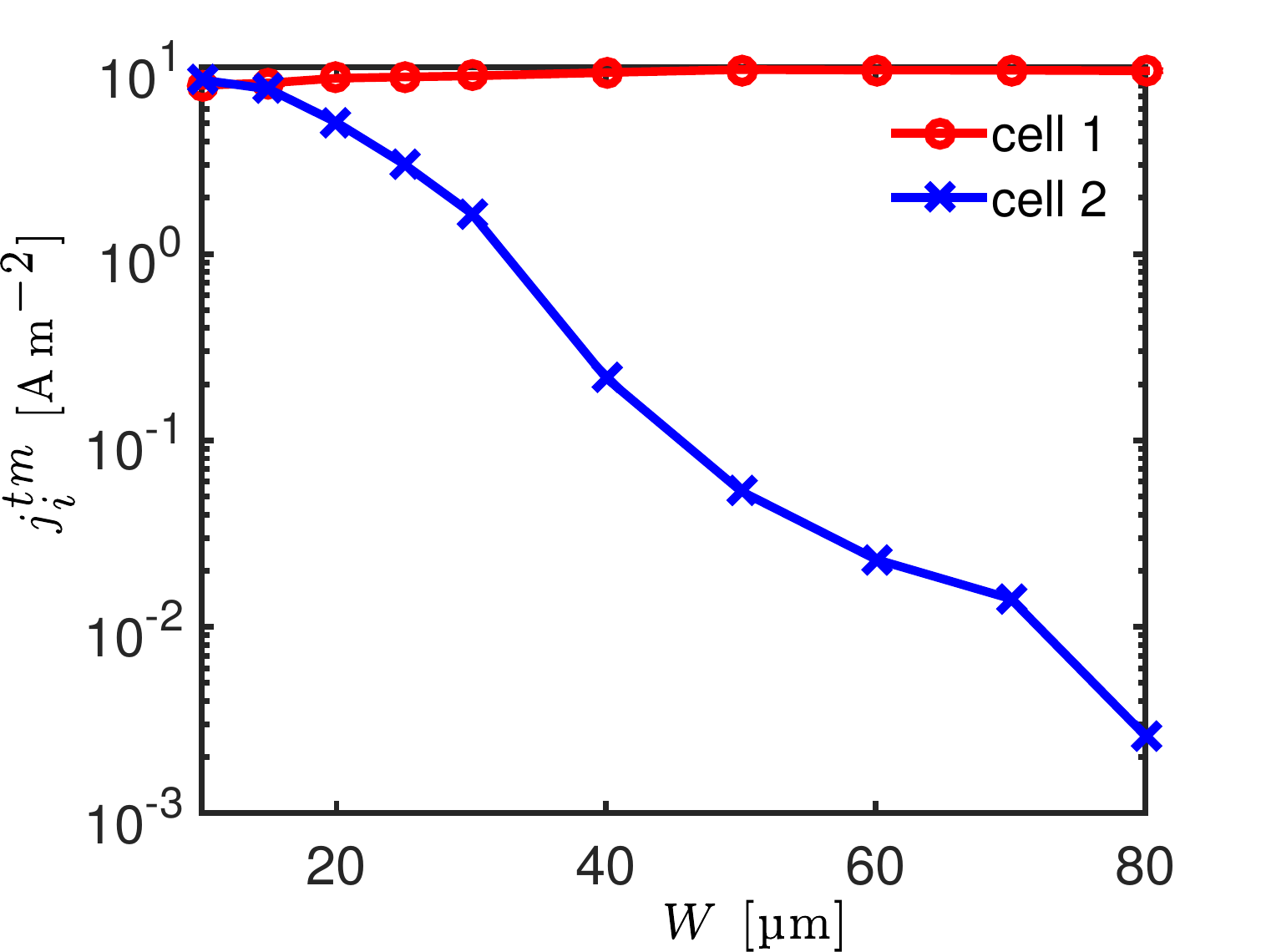}
\caption{Computed transmembrane current densities $j_i^{tm}, \ i=1,2,$ for 
the two cells as a function of $W$. With $j_{1}^{tm}$ we denote the outward current 
flowing from the interior of cell 1 to the electrolyte, 
while with $j_{2}^{tm}$ we denote the inward current flowing into the second cell.
Results obtained with a depolarizing step at $V_{c1}= \SI{50}{\milli\volt}$
and keeping $V_{c2}=V_{ceq}$. 
}
\label{fig:celltocell_param} 
\end{figure}

In Fig.~\ref{fig:celltocell_param} we report the integral averages of 
the transmembrane current densities $j_1^{tm}$ and $j_2^{tm}$, 
flowing out of the first cell and into the second cell, respectively,
as functions of the distance $W$ between 
the two cells. While $j_1^{tm}$ is unaffected by $W$, $j_2^{tm}$ exponentially decays with $W$.

\subsection{Cells with active channels}
\label{sec:active}
In all the numerical experiments described so far,
we have dealt with cells with passive channels. 
We now consider the case of a cardiac cell as the one studied in~\cite{mori2006three},
which expresses voltage-gated ion channels.
In order to describe the behavior of such kind of cells we adopt
the well known Hodgkin-Huxley model (see for the
details~\cite{hodgkin1952currents,hille2001ion,keener2009mathematical}), 
to reproduce the opening/closing dynamics of the channels,
so that the transmembrane ion currents are computed 
at every time $t$ as:
\begin{equation}
\label{eq:tmc_active}
\begin{cases}
j_K =\overline{g}_K n^4 \left(V_{cell}-\left.\varphi\right|_{\Gamma_{cell}}-\left.V_K\right|_{\Gamma_{cell}}\right) 
\\[1ex]
j_{Na} =\overline{g}_{Na} h m^3 \left(V_{cell}-\left.\varphi\right|_{\Gamma_{cell}}-\left.V_{Na}\right|_{\Gamma_{cell}}\right) 
\\[1ex]
j_{leak} =\overline{g}_{Cl} \left(V_{cell}-\left.\varphi\right|_{\Gamma_{cell}}-\left.V_{Cl}\right|_{\Gamma_{cell}}\right), 
\end{cases}
\end{equation}
where $V_K, V_{Na}$ and $V_{Cl}$ are the Nernst potentials
of potassium, sodium and chloride ions at $\Gamma_{cell}$, respectively,
while $\overline{g}_K$ and $\overline{g}_{Na}$ are the maximum values of
potassium and sodium ion conductances.
The value $n^4$ represents the probability of a \ce{K+} channel
to be open,
while the probability that a \ce{Na+} channel
is open is given by the product $m^3 h$.
The gating variables $n, m$ and $h$ need to be consistently computed by 
solving the following system of ODEs at each of the mesh nodes
on the side $\Gamma_{cell}$ representing the cell membrane:
\begin{equation}\label{eq:HH_syst}
\begin{cases}
\dfrac{\mathrm{d}n}{\mathrm{d}t}=\alpha_n(1-n)-\beta_n n 
\\[1.5ex]
\dfrac{\mathrm{d}m}{\mathrm{d}t}=\alpha_m(1-m)-\beta_m m 
\\[1.5ex]
\dfrac{\mathrm{d}h}{\mathrm{d}t}=\alpha_h(1-h)-\beta_h h.
\end{cases}
\end{equation}

Here we study a benchmark
patch-clamp experiment in voltage-clamp configuration,
in which a depolarizing pulse is applied from the resting state
(\SI{-70}{\milli\volt})
holding $V_{cell}$ at \SI{15}{\milli\volt} and
considering the parameter values listed in Table~\ref{tab:Parameters_active}.
\begin{table}[h!]
\centering
\begin{tabular}{lll}\hline
\textbf{Parameter}		& \textbf{Symbol} 	& \textbf{Value} \\[.25ex] \hline
Intracellular potassium concentration    & $c_{K}^{int}$		& $140$ mM \\[.25ex]
Intracellular sodium concentration 	    & $c_{Na}^{int}$		& $10$ mM \\[.25ex]
Intracellular chloride concentration 	& $c_{Cl}^{int}$		& $20$ mM \\[.25ex]
Extracellular bath potassium concentration    & $c_{K}^{bath}$	& $5$ mM \\[.25ex]
Extracellular bath sodium concentration 	    & $c_{Na}^{bath}$	& $145$ mM \\[.25ex]
Extracellular bath chloride concentration 	& $c_{Cl}^{bath}$	& $150$ mM \\[.25ex]
Max potassium conductance   & $\overline{g}_K$			& $180\cdot 10^1$ \si{\siemens\per\square\meter} \\[.25ex]
Max sodium conductance   & $\overline{g}_{Na}$			& $600\cdot 10^1$ \si{\siemens\per\square\meter} \\[.25ex]
Chloride conductance   & $\overline{g}_{Cl}$			& $1.5\cdot 10^1$ \si{\siemens\per\square\meter} \\[.25ex]
Membrane specific capacitance   & $C_M$      & 1 \si{\micro\farad\per\square\centi\meter} \\[.25ex]
Substrate specific capacitance  & $C_S$    & 0.3 \si{\micro\farad\per\square\centi\meter} \\[.25ex]
Initial transmembrane potential  & $V_{cell}-\left.\varphi\right|_{\Gamma_{cell}}$    & -70 \si{\milli\volt}  \\[.5ex]
\hline
\end{tabular} 
\caption{Model parameter values considered in the simulations involving
active cardiac cell.}
\label{tab:Parameters_active}
\end{table}

\begin{figure}[h!]
\centering 
\subfigure[$n^4(r)$]
{\includegraphics[width=0.435\textwidth]
{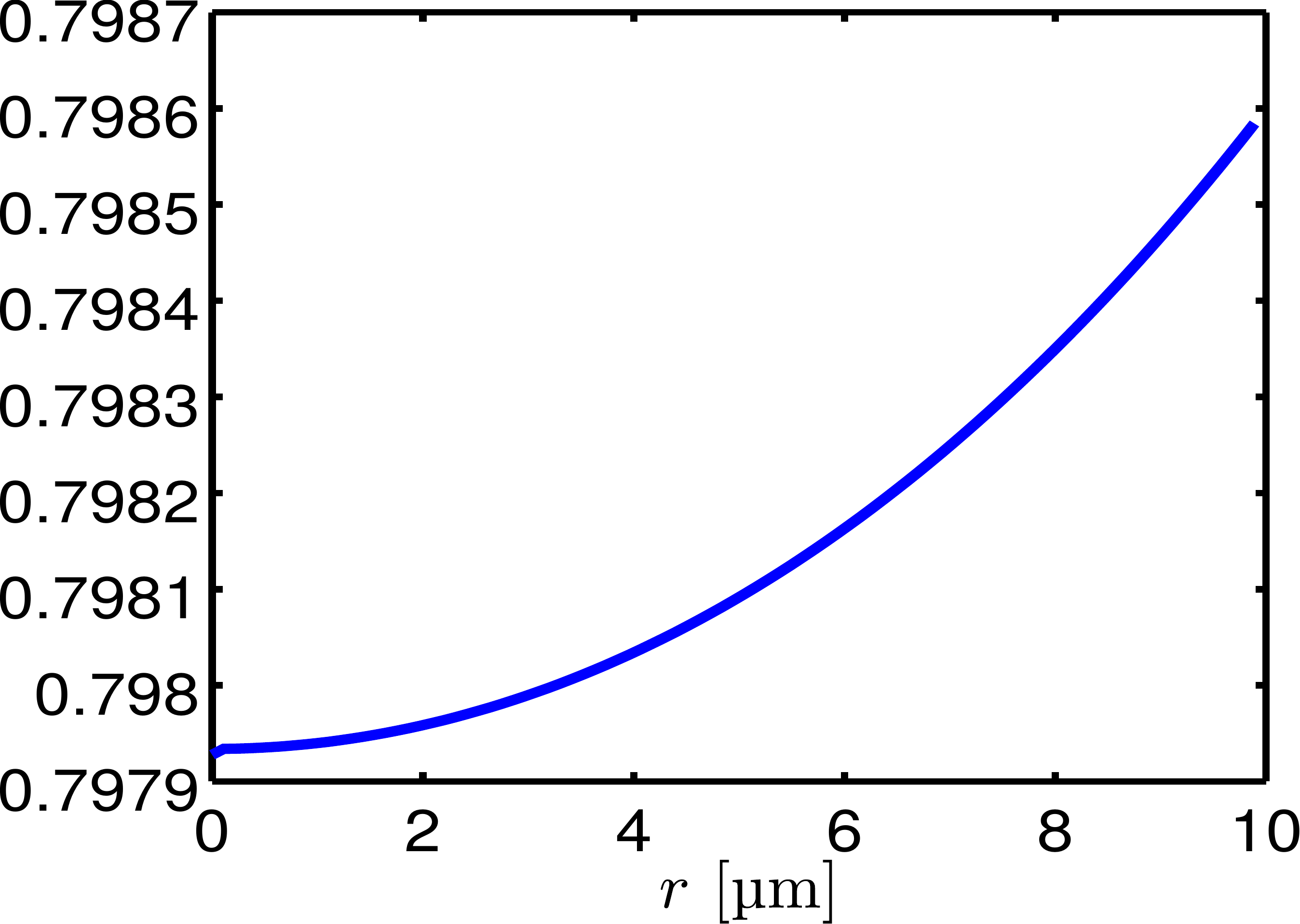}
\label{fig:nspace}} 
\hfill
\subfigure[$hm^3(r)$]
{\includegraphics[width=0.462\textwidth]
{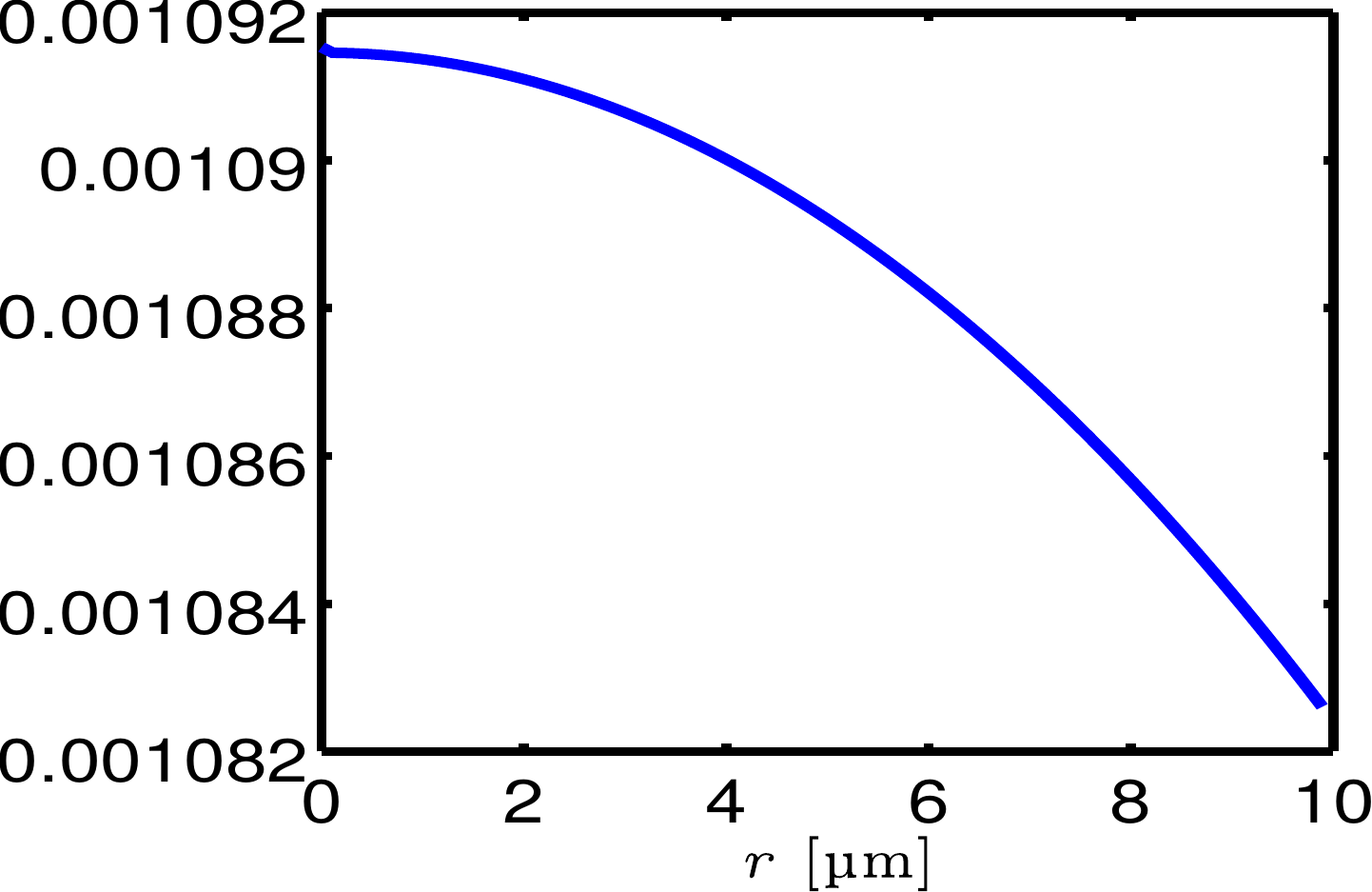}
\label{fig:mspace}}
\caption{Computed spatial distributions of the \ce{K+} and \ce{Na+} channel
opening probabilities for every point of the boundary $\Gamma_{cell}$
representing the cell membrane, \SI{10}{\milli\second} after the onset of a 
depolarizing voltage-clamp pulse at $V_{cell} = \SI{15}{\milli\volt}$.}
\label{fig:gating_space} 
\end{figure}

\begin{figure}[h!]
\centering 
\subfigure[$\varphi$]
{\includegraphics[width=0.45\textwidth]{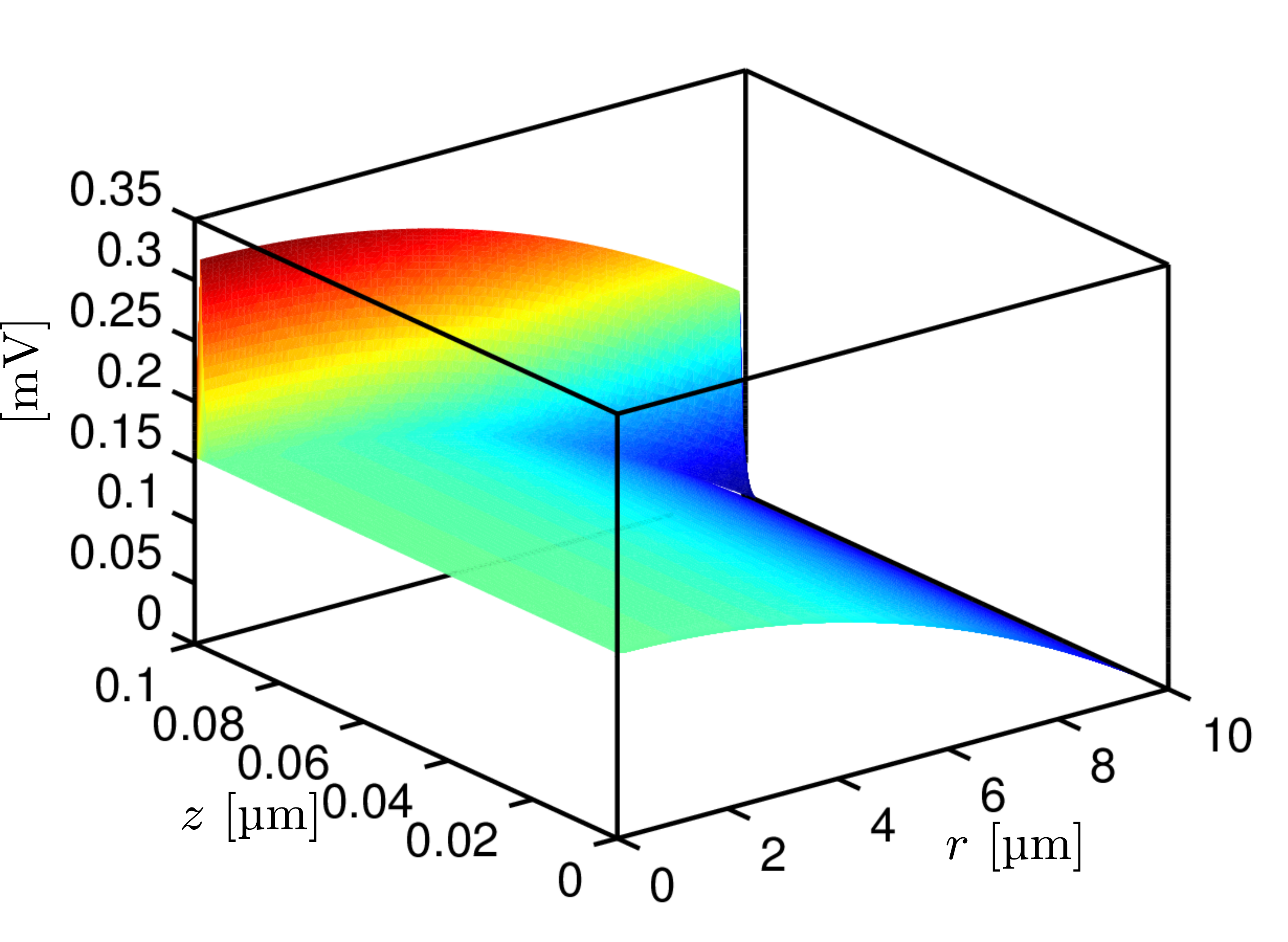}
\label{fig:phiact}}
\hfill
\subfigure[$c_{K}$]
{\includegraphics[width=0.45\textwidth]{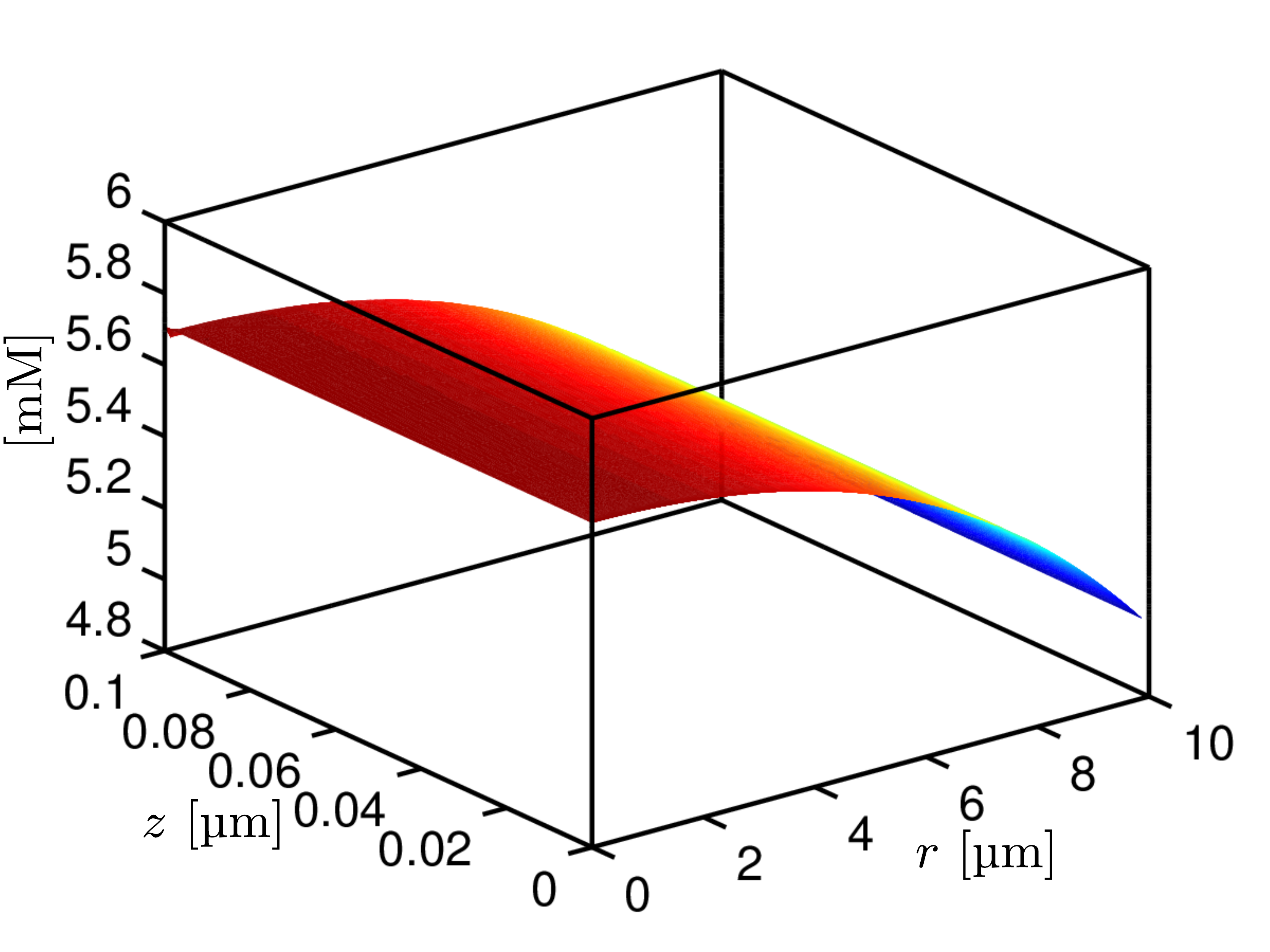}
\label{fig:kact}}
\subfigure[$c_{Na}$]
{\includegraphics[width=0.45\textwidth]{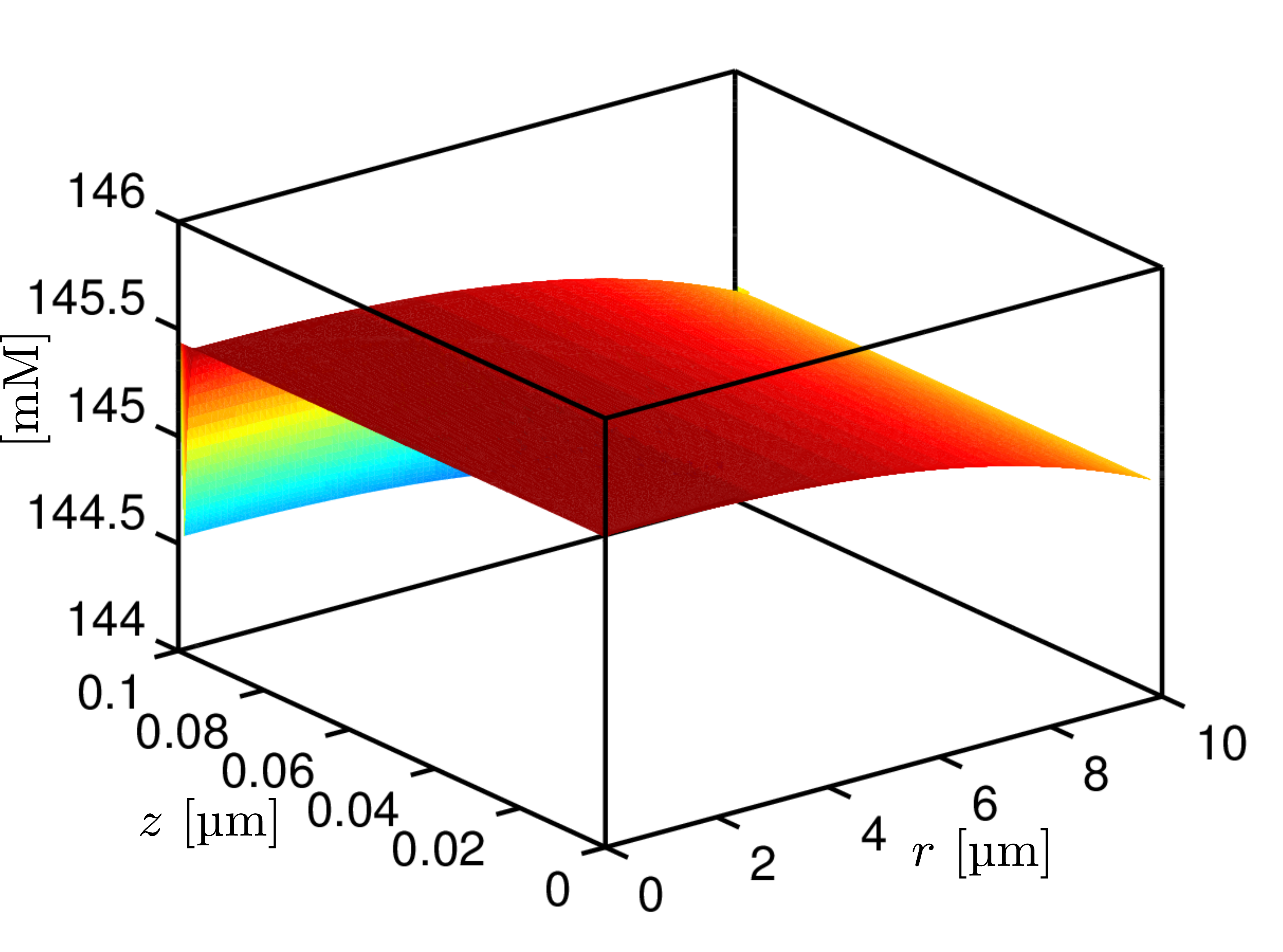}
\label{fig:naact}} 
\caption{Spatial distributions of the potential and of the potassium
and sodium concentrations in the electrolyte cleft under the cell.
 Results obtained after 10\,ms of a depolarizing pulse
keeping $V_{cell}=$15\,mV.
}
\label{fig:distrib_active} 
\end{figure}

Fig.~\ref{fig:gating_space} shows the computed opening probability profiles of 
\ce{K+} and \ce{Na+} channels, respectively, \SI{10}{\milli\second} after the onset of 
the stimulation, and as expected we can observe that the potassium channels are open while 
the sodium channels are already inactivated. 
The limited spatial variability of the obtained profiles 
agrees with the fact that
in the considered patch clamp configuration, the potential is modified 
in the whole cell domain, so the channels on the membrane respond 
almost uniformly.
Fig.~\ref{fig:distrib_active} shows the spatial distributions 
of potential and ion concentrations of the species
with active channels (\ce{K^+} and \ce{Na^+}). Notice the occurrence
of sharp boundary layers at the cell-electrolyte interface due
to sodium channel inactivation.

\begin{figure}[h!]
\centering 
\subfigure[Gating variables]
{\includegraphics[width=0.45\textwidth]
{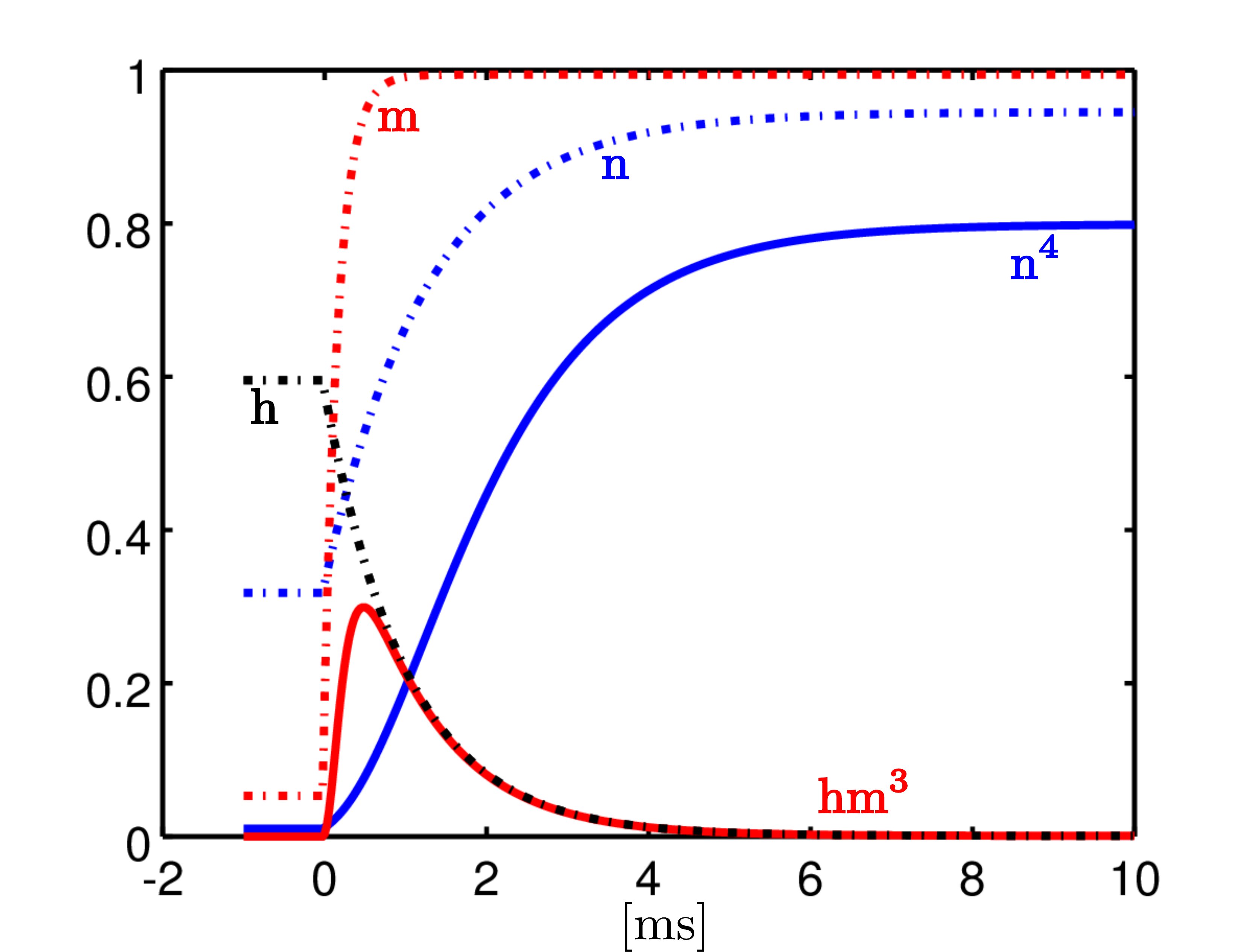}
\label{fig:gating_var}}
\hfill
\subfigure[Transmembrane currents]
{\includegraphics[width=0.45\textwidth]
{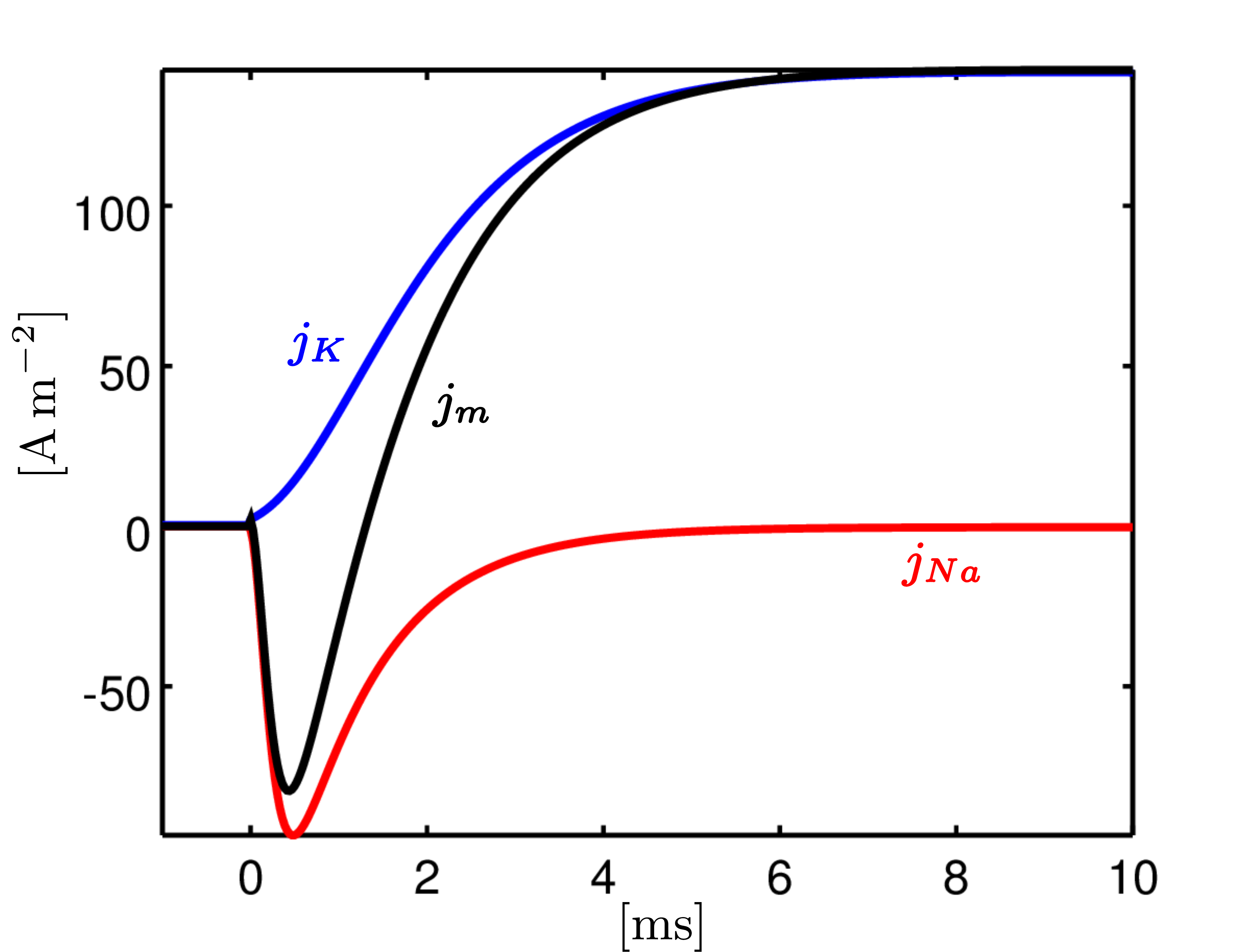}
\label{fig:trans_currents}} 
\caption{Time evolution of the gating variables and ion channel opening probabilities (left), and of the 
ion current densities (right).
Results averaged over $\Gamma_{cell}$ and obtained 
with a depolarizing voltage pulse at $V_{cell}=15\,mV$.}
\label{fig:active_15mv}
\end{figure}

In order to compare our results with the voltage clamp experimental measurements
on cells with active channels, we compute the integral average 
of the channel opening probabilities $n^4$ and $hm^3$
over $\Gamma_{cell}$ and we show their evolution in Fig.~\ref{fig:gating_var}. 
These results are in excellent agreement with the behavior 
characteristic of the classic HH model
(see. e.g.,~\cite{ermentrout2010,hille2001ion,keener2009mathematical}), 
but, while this latter is usually applied in a lumped manner,
our model has the feature of describing possible spatial inhomogeneities
in the activation (cf. Fig.~\ref{fig:distrib_active}). 
In Fig.~\ref{fig:trans_currents} we report the 
average ion channel transmembrane current densities 
computed according to~\eqref{eq:tmc_active}, and, as expected, 
the typical current profile after a voltage-clamp depolarization
is recovered. After the onset of the voltage step, an inward
current is induced by the opening of \ce{Na+} channels, which are 
eventually inactivated, and for longer times the \ce{K+} channels 
open and determine a stationary outward flux.
\begin{figure}[h!]
\centering 
\subfigure[$g_{K}=n^4\overline{g}_{K}$]
{\includegraphics[width=0.45\textwidth]
{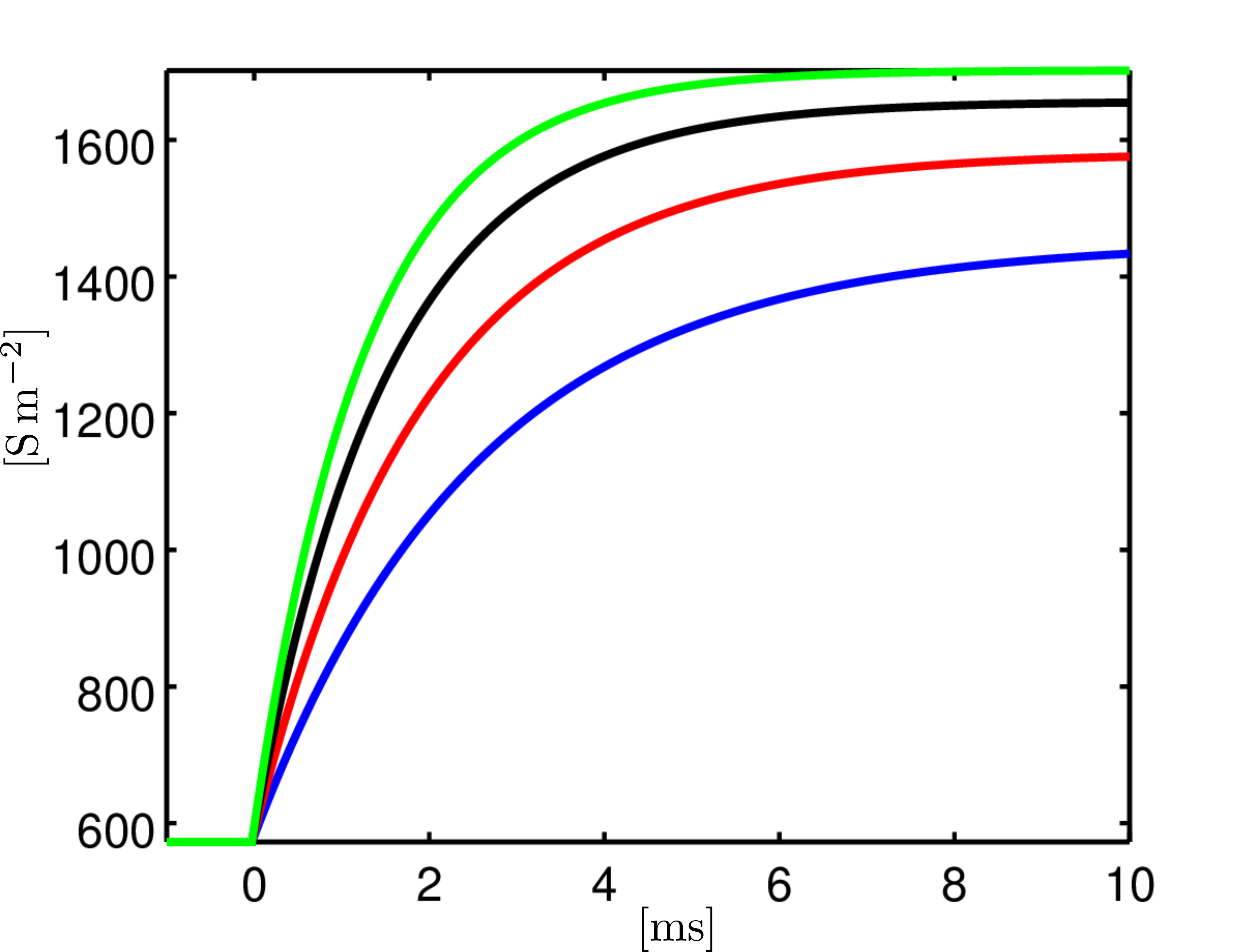}
\label{fig:gk4}} 
\hfill
\subfigure[$g_{Na}=hm^3\overline{g}_{Na}$]
{\includegraphics[width=0.45\textwidth]
{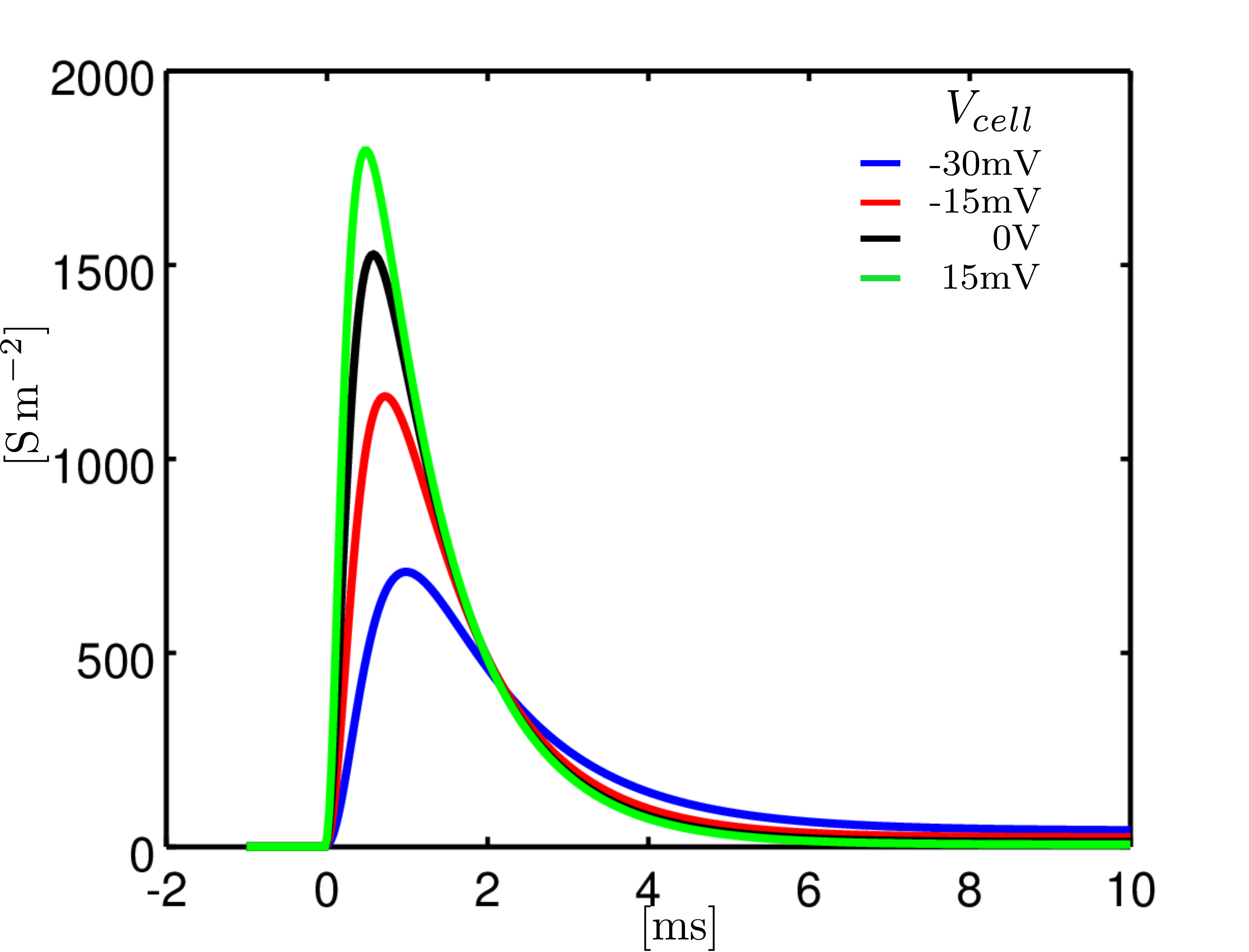}
\label{fig:gna4}} 
\caption
{\label{fig:different_conduct} Temporal variation of the integral mean over 
the boundary $\Gamma_{cell}$ of the conductances of the potassium and of the sodium 
channels, computed as in~\eqref{eq:tmc_active}. Results obtained with four different depolarizing pulses
keeping $V_{cell}=-30,-15, 0, +15$\,mV.}
\end{figure}

As a final result,  we consider the effects of different depolarizations
of the cell over the ion transmembrane currents. In Fig.~\ref{fig:different_conduct}
we show the variation in time of potassium and sodium conductances
in correspondence of four different values of $V_{cell}$. We observe that
$g_{Na}$ turns on more rapidly than $g_{K}$. Moreover, the Na$^+$ channels
begin to close before depolarization is turned off, whereas the 
K$^+$ channels
remain open as long as the membrane is depolarized. As in the classic HH theory,
when the cell is depolarized the Na$^+$ channels
 switch from the resting (closed) to the activated (open)
state and then, if depolarization is maintained, the channel switches to the inactivated state again.

\subsection{Reduced order models}\label{sec:reduced_order_validation}
In this section we illustrate the results obtained solving the reduced
models introduced in Sect.~\ref{sec:hierarchy}.
\subsubsection{Validation of the 2.5D model against the 3D PNP model}
In this section we verify the accuracy of the reduced model
of Sect.~\ref{sub:Model-reduction} 
in the study of the same biophysical configuration analyzed 
in Sect.~\ref{sec:voltage_clamp_validation} with the 
3D axisymmetric PNP model.
Since the geometrical setting is axisymmetric,
we can reduce ourselves to a one dimensional
manifold $\Omega_{rad}$, describing the variation 
of the quantities of interest along the sole radial direction. 
The considered domain $\Omega_{rad}$ is the union of two different
parts $\Omega_{r}^{cell}\cup\Omega_{r}^{ef}$ (the 
part where the cell is attached and the free part of extracellular
fluid). For the approximation of the boundary layers
introduced in Sect.~\ref{sub:Model-reduction}, 
based on the simulations of the previous sections and on 
asymptotic analysis (see~\cite{mori2006three}), 
we set $H\simeq2\lambda_{Debye} \simeq$ \SI{1.6}{\nano\meter}.

\begin{figure}[h!]
\centering 
\subfigure[$\varphi$]
{\includegraphics[width=0.45\textwidth]{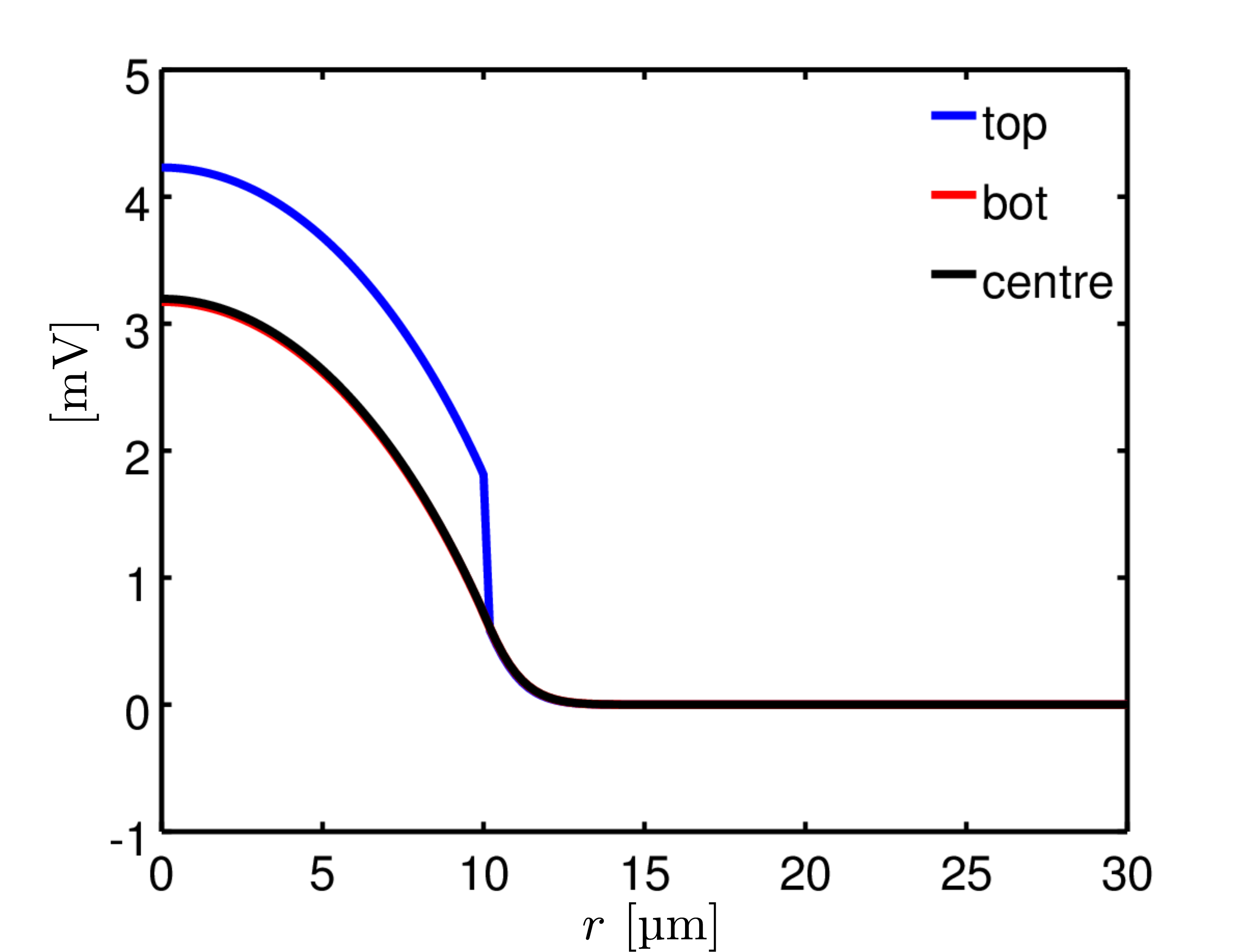}
\label{fig:pot_lambda}} 
\hfill
\subfigure[Cl]
{\includegraphics[width=0.45\textwidth]{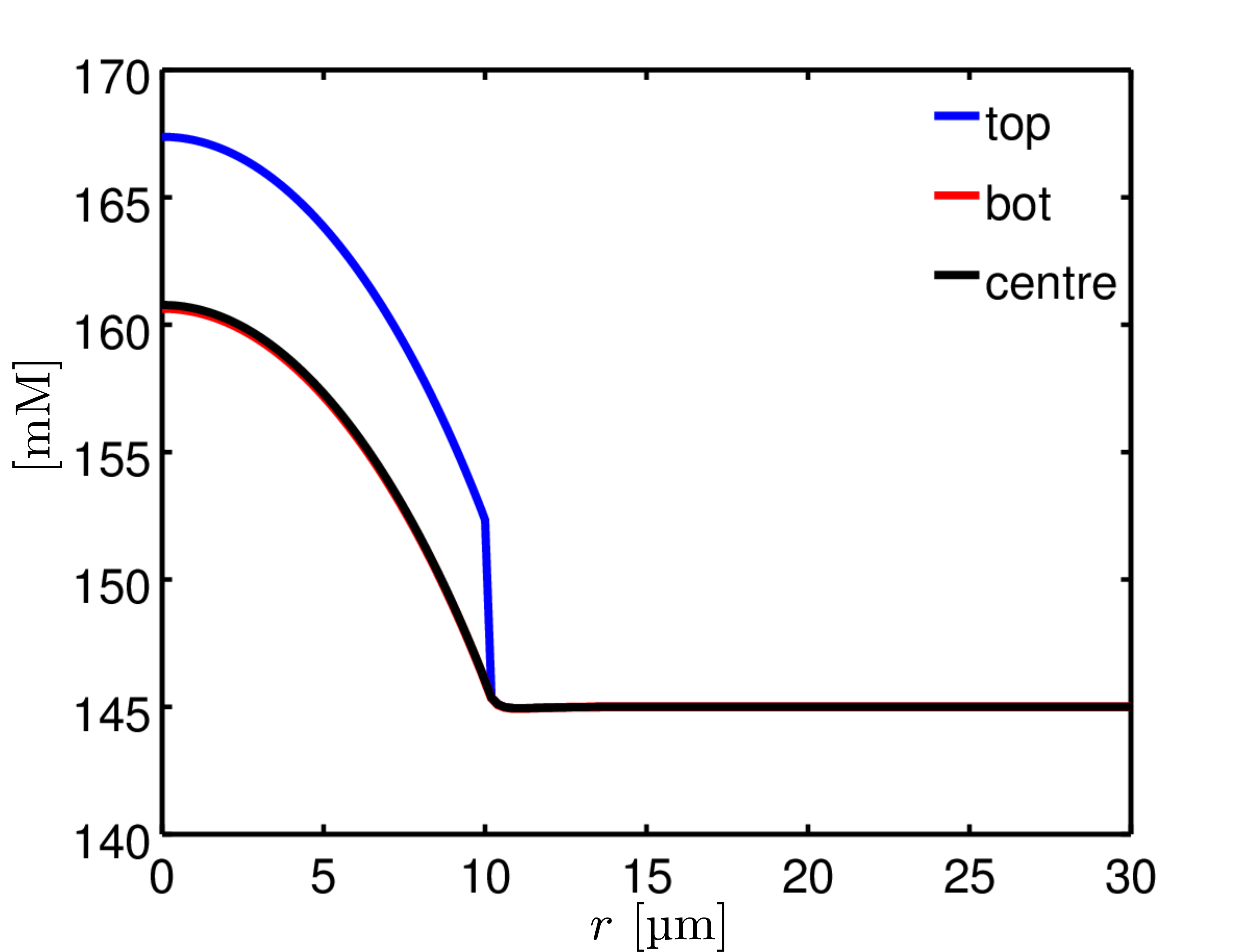}
\label{fig:cCl_lambda}} 
\hfill
\subfigure[K]
{\includegraphics[width=0.45\textwidth]{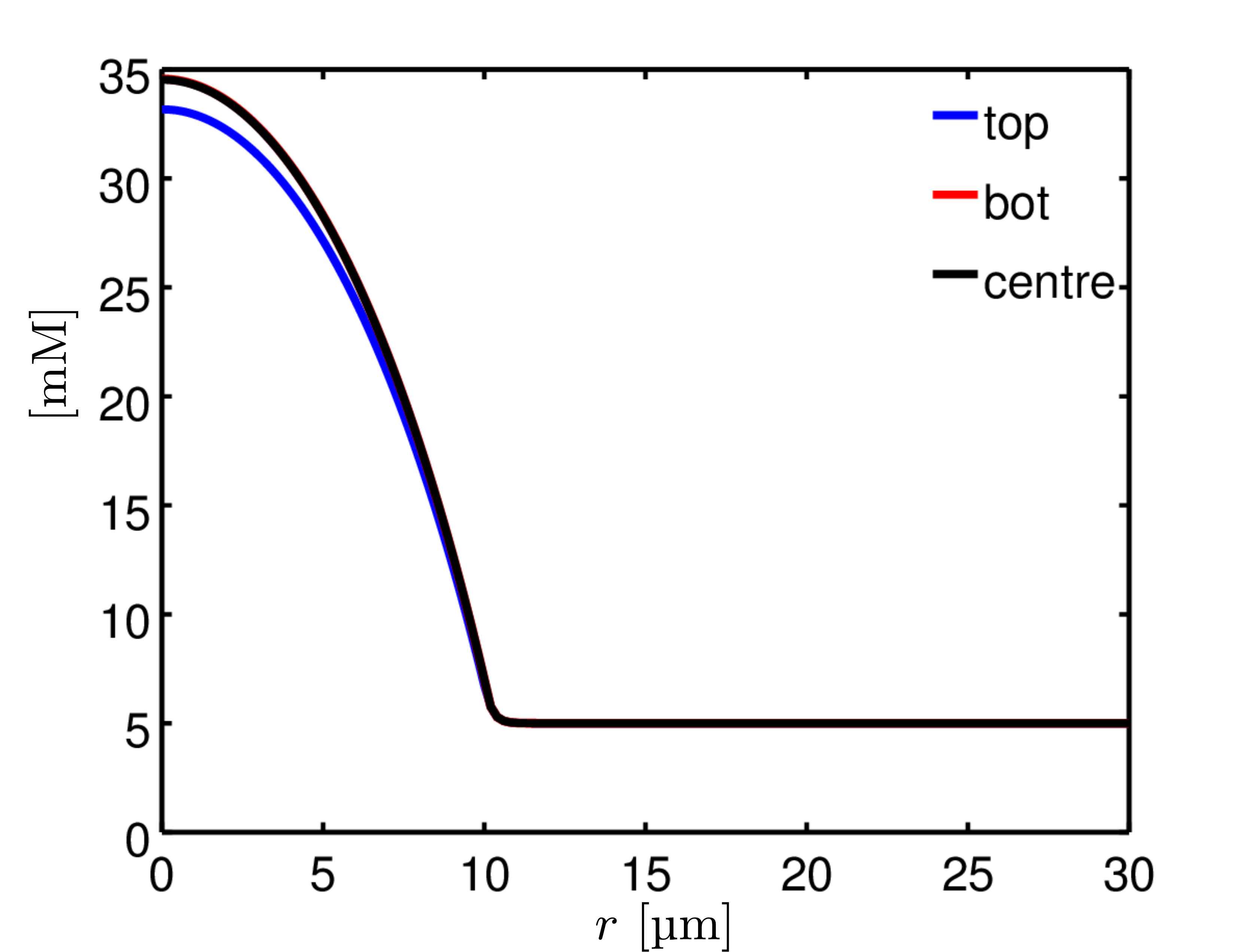}
\label{fig:cK_lambda}} 
\hfill
\subfigure[Na]
{\includegraphics[width=0.45\textwidth]{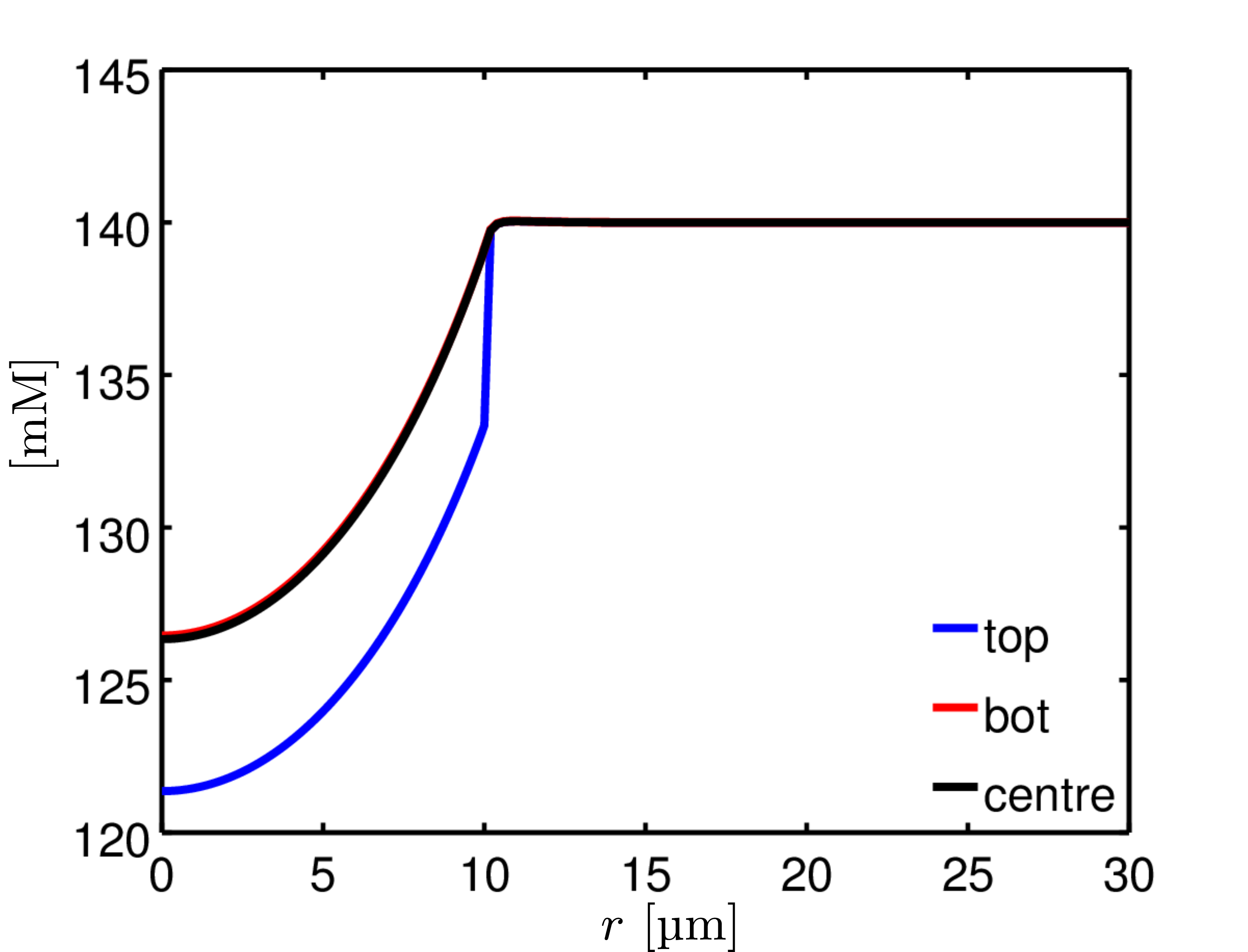}
\label{fig:cNa_lambda}} 
\caption{Spatial distribution
of $\,\overline{\varphi}$ and of $\,\overline{c}_{i}$ in the domain
$\Omega_{rad}=\Omega_{r}^{cell}\cup\Omega_{r}^{ef}$ (results obtained
with $|\Omega_{r}^{ef}|=2\cdot
|\Omega_{r}^{cell}|=2\cdot R_{cell}=$
\SI{20}{\micro\meter}). To account for the boundary layers: distributions
of the top and the bottom values $\varphi_{top}$, $\varphi_{bot}$,
$c_{i}^{top}$ and $c_{i}^{bot}$.}
\label{fig:reduction_not_justc}
\end{figure}
Simulation results obtained with the reduced model are shown 
in Fig.~\ref{fig:reduction_not_justc}, where we immediately notice that the capacitive couplings obtained with 3D simulations of Sect.~\ref{sec:voltage_clamp_validation} are well reproduced.
As observed in Sect.~\ref{sub:Model-reduction} 
and in the spatial distributions of Sect.~\ref{sub:Cell-to-chip}, 
the quantities of interest in $\Omega_{r}^{ef}$ 
are not expected to sensibly vary along $z$. 
This solution behavior is confirmed by the radial distributions
of Fig.~\ref{fig:reduction_not_justc}: potential $\overline{\varphi}$
and concentrations $\overline{c}_{i}$ are perfectly superimposed
on the distributions of the top and bottom quantities in the free part. We observe
a very fast decay in $\Omega_{r}^{ef}$, as expected and as the one obtained in
the simulation results shown in Fig.~\ref{fig:pabst_withcouplings}.

\begin{figure}[h!]
\centering
\subfigure[$\varphi(r=0,z)$]
{\includegraphics[width=0.45\textwidth]{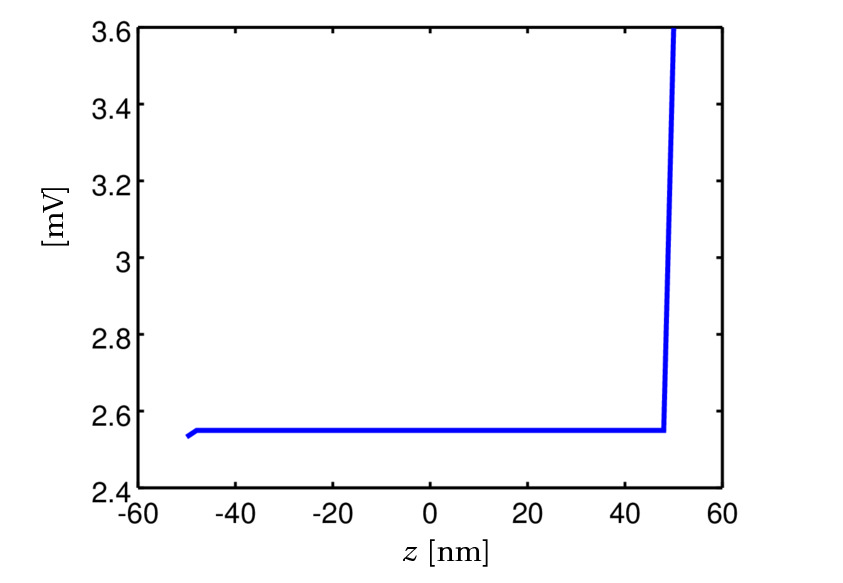}
\label{fig:pot_reduceddomain_centre}} 
\hfill
\subfigure[$c_K(r=0,z)$]
{\includegraphics[width=0.45\textwidth]{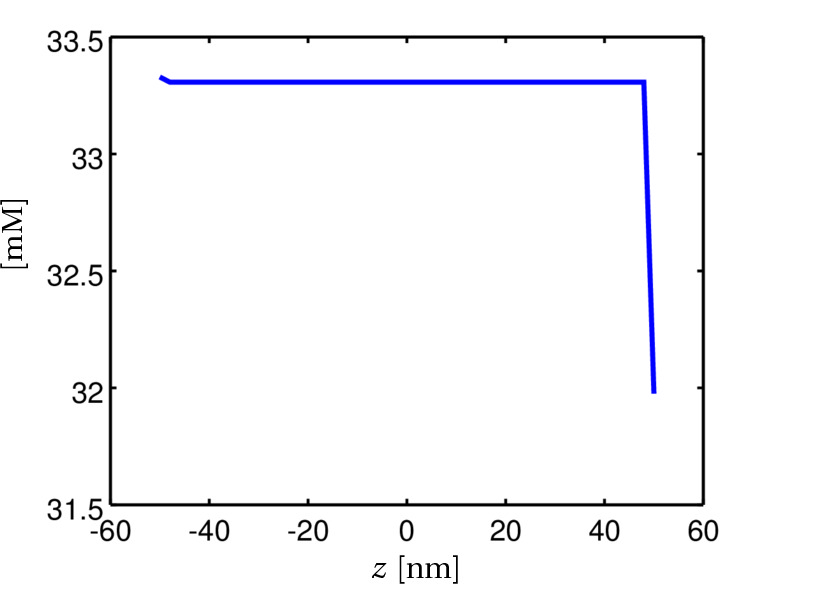}
\label{fig:concK_reduceddomain_centre}} 
\caption{Distributions along the $z$-direction
of potential and potassium concentration at $r=0$ using
\eqref{eq:post_zrec}.}
\label{fig:z_dir_recons}
\end{figure}

Having computed the averaged and the top and bottom values 
of the dependent variables, we can also reconstruct their 
$z$-dependence, at a fixed point $\bar{r}$, using the following
post-processing formulas:
\begin{subequations}\label{eq:post_zrec}
\begin{align}
& \varphi(\bar{r},z) =  \overline{\varphi}(z,\bar{r})+
\left.\left(\dfrac{\varphi_{top}(\bar{r},z)-
\overline{\varphi}(\bar{r},z)}{H}\right)\right|_{z\in\Omega_{1}}+
\left.\left(\dfrac{\overline{\varphi}(\bar{r},z)-\varphi_{bot}(\bar{r},z)}{H}
\right)\right|_{z\in\Omega_{3}}\label{eq:phi_z_post}\\
& c_{i}(\bar{r},z) =  \overline{c}_{i}
\exp\left(-z_{i}\dfrac{\varphi(\bar{r},z)-
\overline{\varphi}(\bar{r},z)}{V_{th}}\right),\label{eq:c_z_post}
\end{align}
\end{subequations}
where $\Omega_{1}$ and $\Omega_{3}$ are the boundary
layer subdomains (see Fig.~\ref{fig:boundary_layer}). The reconstruction
of $\varphi(r=0,z)$ and of $c_{K}(r=0,z)$ is illustrated in Fig.
\ref{fig:z_dir_recons} and we see that the major term in the electrolyte
system behavior is the capacitive coupling with the cell membrane.

The above obtained results allow us to conclude that the reduced 
mathematical model of Sect.~\ref{sub:Model-reduction} is valid, because 
its predictions favorably agree with those of the 3D model, but
with a much smaller amount of degrees of freedom.

\subsubsection{Validation of the 2.5D model against exact solutions}
In this section we analyze a biophysical setting similar to that
presented in~\cite{Pabst2007}, where the authors derive the
exact profiles of potential and ion concentrations
in radial coordinates under
the assumption of axial symmetry. Their model refers to the middle
plane of a cleft between a cell and a substrate as in the general
setup described in Sect.~\ref{sec:hierarchy}, but under 
the following simplifications:
\begin{itemize}
\item ions are assumed to flow only in the radial direction, so that $j_{i,\phi}=j_{i,z}=0$, and no spatial dependence on $\phi$ and
$z$ is assumed. The radial coordinate $r$ is taken 
in the interval $[0,R]$ ($R$ being the cell radius);
\item the influx of $\text{K}^{+}$ ion charge
per volume and per time is given by $\lambda_{K}=j_{K}^{tm}/\delta_{j}$,
where $\delta_{j}$ is the cleft width and $j_{K}^{tm}$ is the potassium
current density through the membrane, here assumed to be constant. 
Therefore in Eq.~\eqref{eq:continuity_2d} we only have $f^{top}_{K}=\lambda_{K}$;
\item there is no influence on the flux of ions inside the cleft from the
two interfaces (the cell-cleft and the chip-cleft interfaces), neglecting
the capacitive couplings described in Sect.~\ref{sub:Boundary-and-initial}. 
This gives $g^{top}=g^{bot}=0$ in Eq.~\eqref{eq:poiss_2d};
\item only the stationary case is considered, meaning that all quantities
are independent of time.
\end{itemize}
\begin{figure}[h!]
\centering
\subfigure[Potential]
{\includegraphics[width=0.42\textwidth]{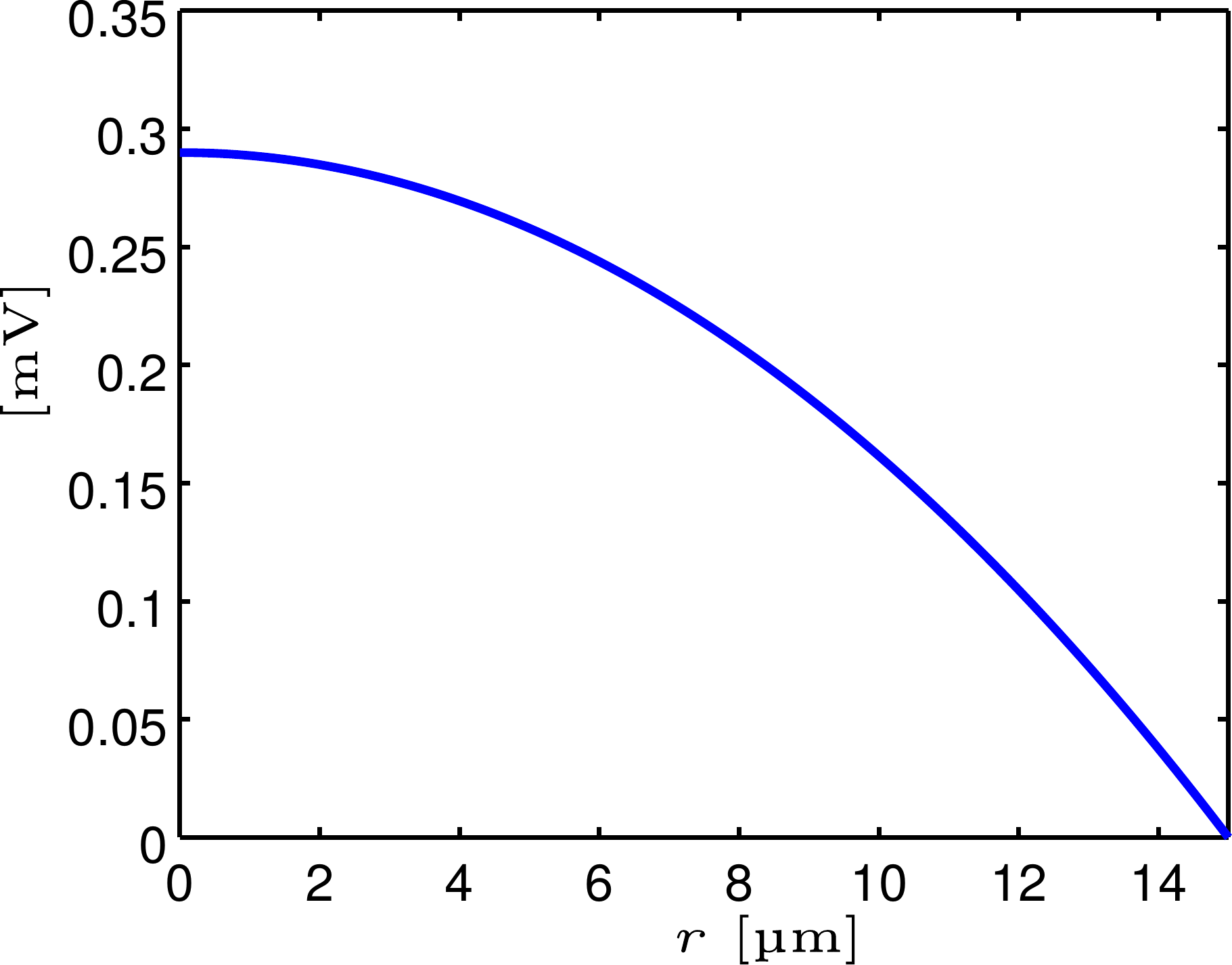}
\label{fig:pot1Dpabst}}
\hfill
\subfigure[Concentrations]
{\includegraphics[width=0.4\textwidth]{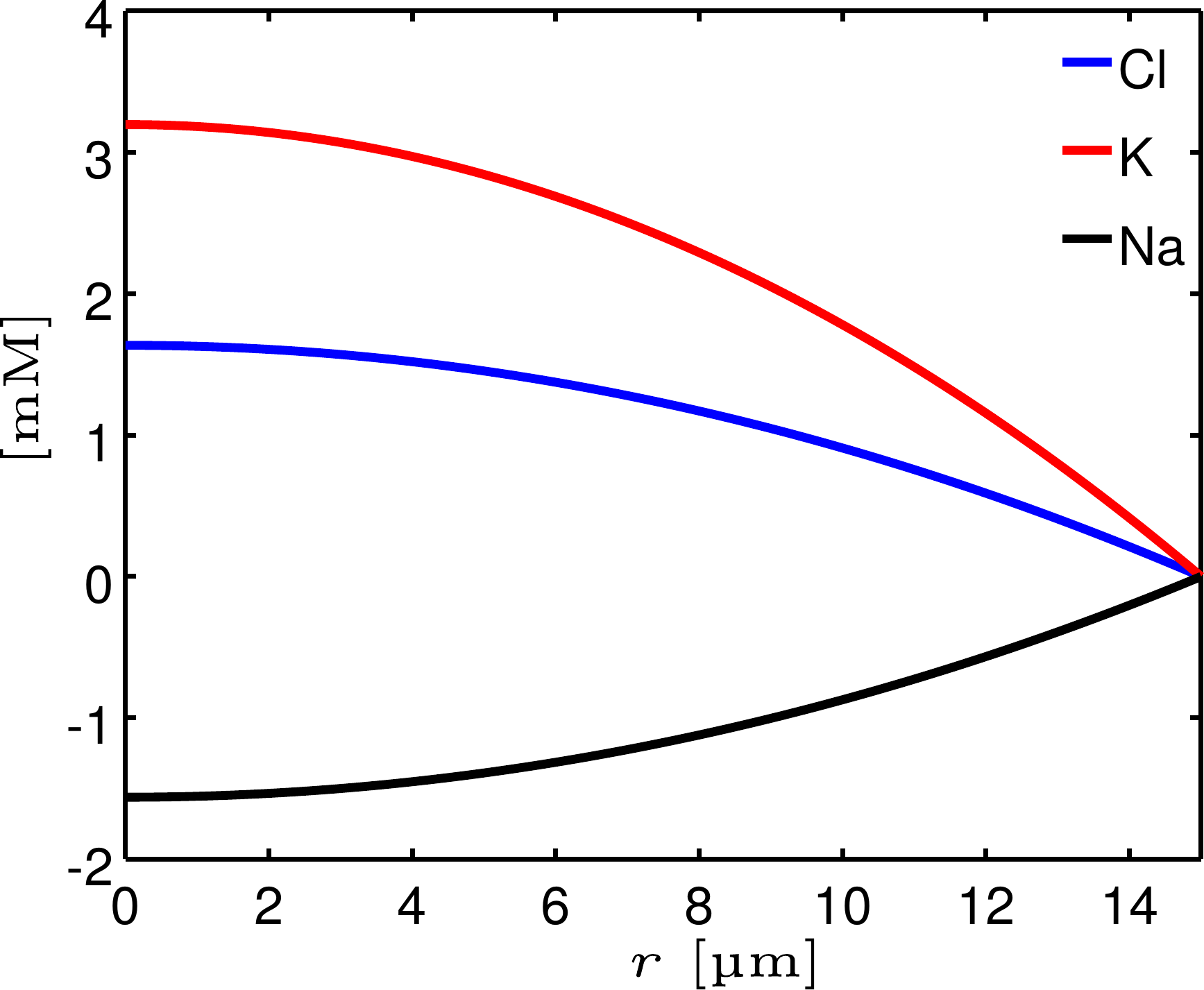}
\label{fig:concion1Dpabst}}
\caption[Radial distribution of the potential and of the concentrations]
{\label{fig:pabst}On the left: radial profile of the
potential $\varphi\left(r\right)$. On the right: radial profile of
the changes of ion concentrations with respect to their bath values
$c_{i}\left(r\right)-c_{i}^{bath}$. Results obtained with a source
term $\lambda_{K}=$ \SI{11}{\pico\ampere\per\square\micro\meter}
and a cell radius $R=$ \SI{15}{\micro\metre} as in~\cite{Pabst2007}.}
\end{figure}

\begin{figure}[h!]
\centering 
\includegraphics[width=0.45\textwidth]{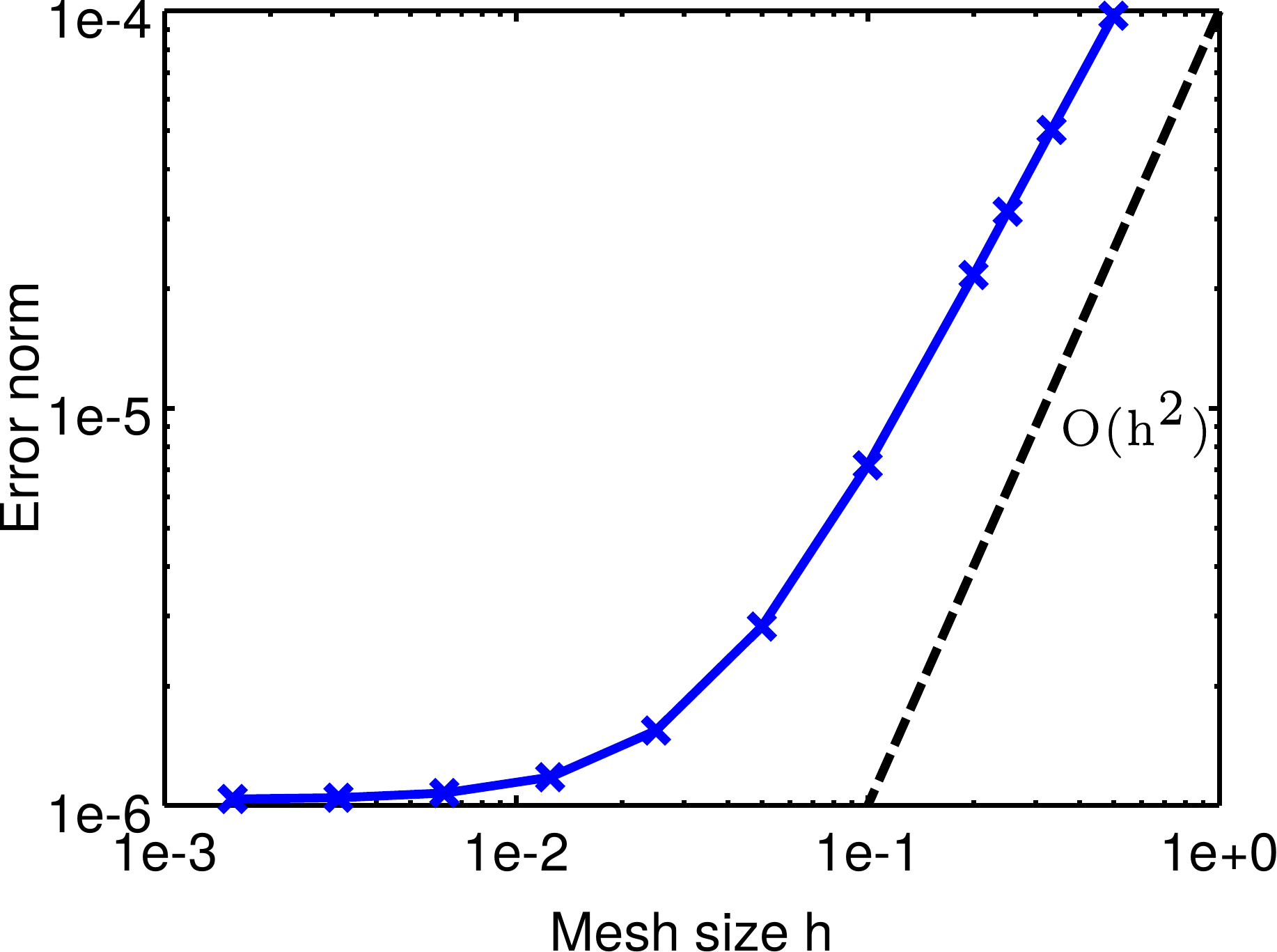}
\caption{Log-log plot of the maximum norm of the discretization error
$\varphi - \varphi_h$.}
\label{fig:convergence_pabst}
\end{figure}

Fig.~\ref{fig:pabst} shows the computed distributions of potential
and ion concentrations under the cell. By inspection on the 
analytical solutions reported in~\cite{Pabst2007} it can be seen that
our results are in excellent agreement with the latter solutions.
Convergence of the finite element solution as a function of the 
radial mesh size $h$ is reported in Fig.~\ref{fig:convergence_pabst}
where a quadratic rate can be observed before the occurrence of 
error saturation due to the nonlinear solver tolerance. This result,
obtained in the solution of a nonlinear problem, 
confirms the convergence analysis of 
Sect.~\ref{sec:convergence_analysis}, carried out on a linear
model problem.

Under
the hypotheses illustrated above, the potential variation at $r=0$
is quite small (less than \SI{1}{mV}) and the absolute changes
of ion concentrations for $\text{Cl}^{-}$ and $\text{Na}^{+}$ are
quite small too, except for $\text{K}^{+}$ ion concentration, 
from 5 mM to 8 mM.
\begin{figure}[h!]
\centering
\subfigure[Potential $\varphi(t,r=0)$]
{\includegraphics[width=0.42\textwidth]{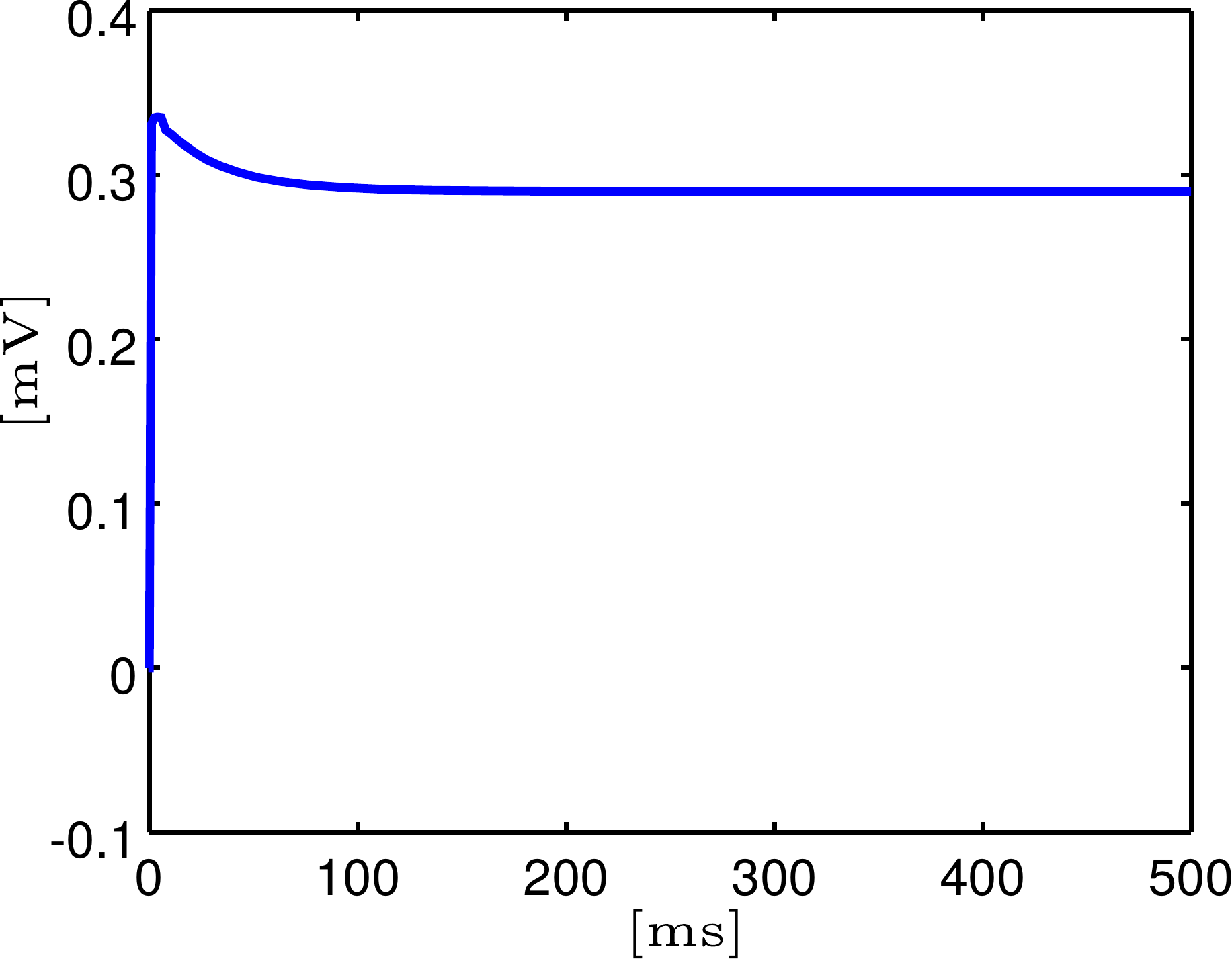}
\label{fig:pot1Dpabst_time}}
\hfill
\subfigure[Concentrations $c_i(t,r=0)-c_i^{bath}$]
{\includegraphics[width=0.4\textwidth]{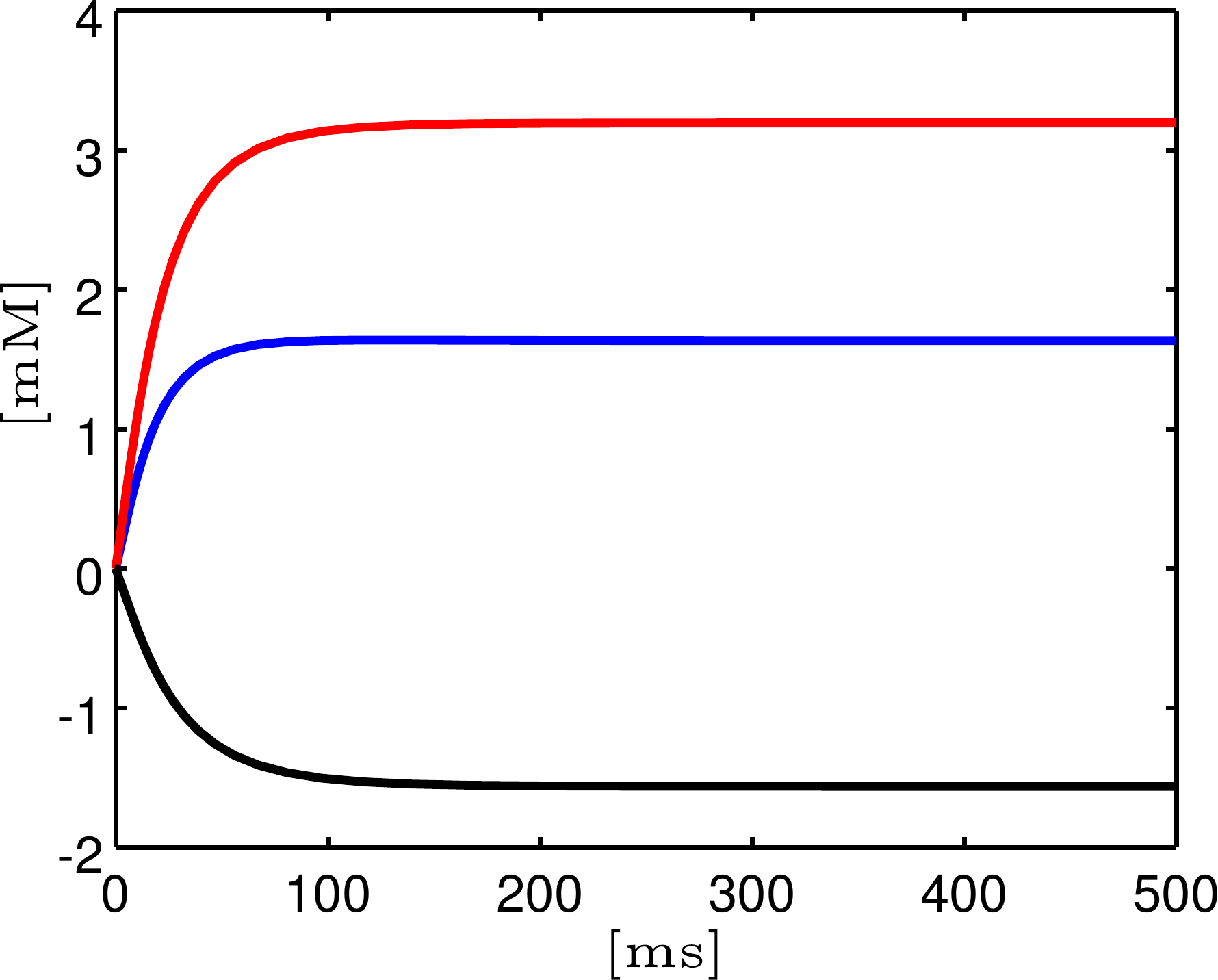}
\label{fig:concion1Dpabst_time}}
\caption[Time variation of the potential and of the concentrations]
{\label{fig:pabst_time}Time variation of the potential
and of the concentrations with respect to their bath values, at the
center of the junction.}
\end{figure}

We also conduct a time dependent simulation of this experimental setup,
redefining the transmembrane current as $\lambda_{K}\left(t\right)=\lambda_{K}H\left(t\right)$,
where $\lambda_{K}$ is the constant current used in the static 
simulation and $H\left(t\right)$ is the Heaviside function.
With this mathematical definition of potassium injection, we consider 
an instantaneous opening of the $\text{K}^{+}$ channels
at $t=0$, which leads to a time variation of the quantities $\varphi$
and $c_{i}$ (see Fig.~\ref{fig:pabst_time}, where we show the variation
of these functions evaluated at $r=0$). Transients are exhausted
in about \SI{150}{ms}, in agreement with the results of~\cite{wrobel2005cell}, 
while the steady-state values of $\varphi$ and $c_{i}$ agree well with those computed in the static case shown in Fig.~\ref{fig:pabst}.
\subsubsection{Validation of the Area-Contact model}
\label{sec:area_contact_simul}
\begin{figure}[h!]
\centering 
\subfigure[Intracellular potential $V_{cell}$]
{\includegraphics[width=0.4\textwidth]{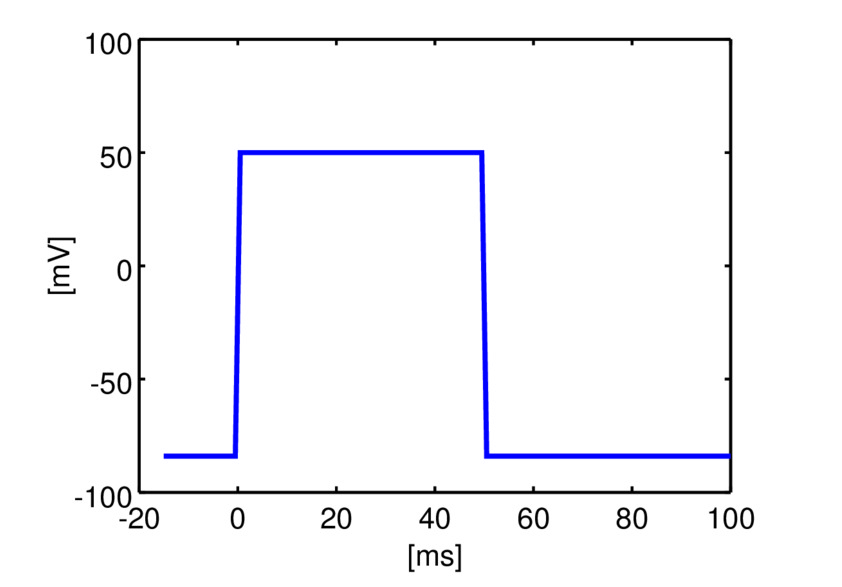}
\label{fig:pot_intracell}} 
\quad
\subfigure[$V_J$ electrical model]
{\includegraphics[width=0.38\textwidth]{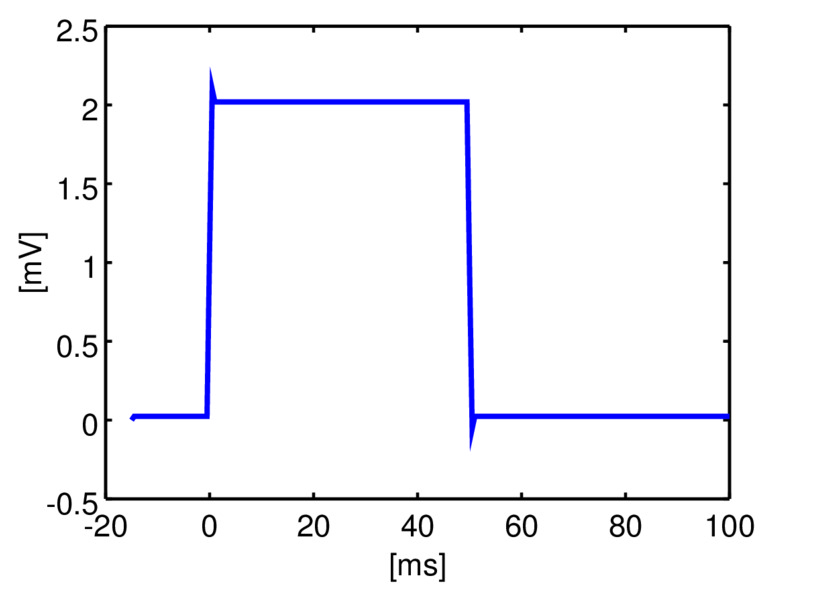}
\label{fig:potcleft_electricalmodel}}
\subfigure[$V_J$ electrodiffusion model]
{\includegraphics[width=0.4\textwidth]{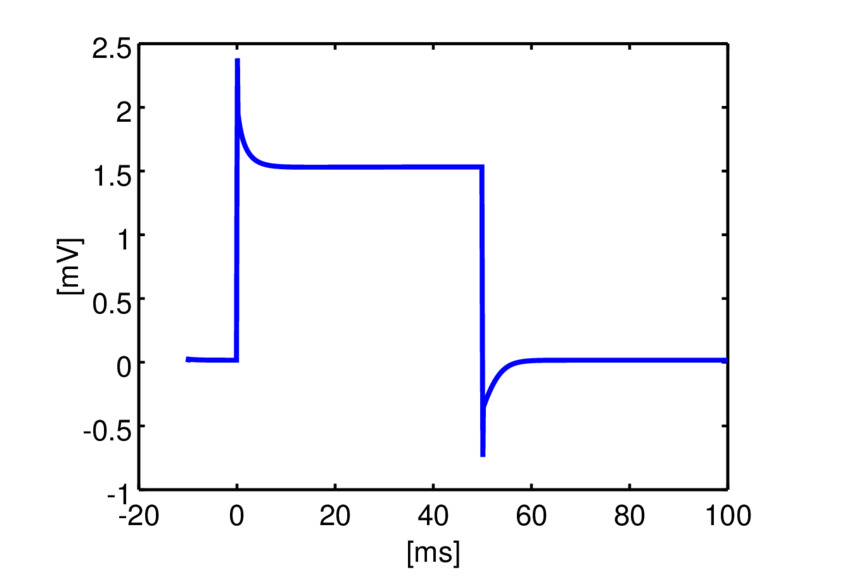}
\label{fig:electrodiff_mean}}
\quad
\subfigure[Nernst potentials $V_{J0}^i$]
{\includegraphics[width=0.38\textwidth]
{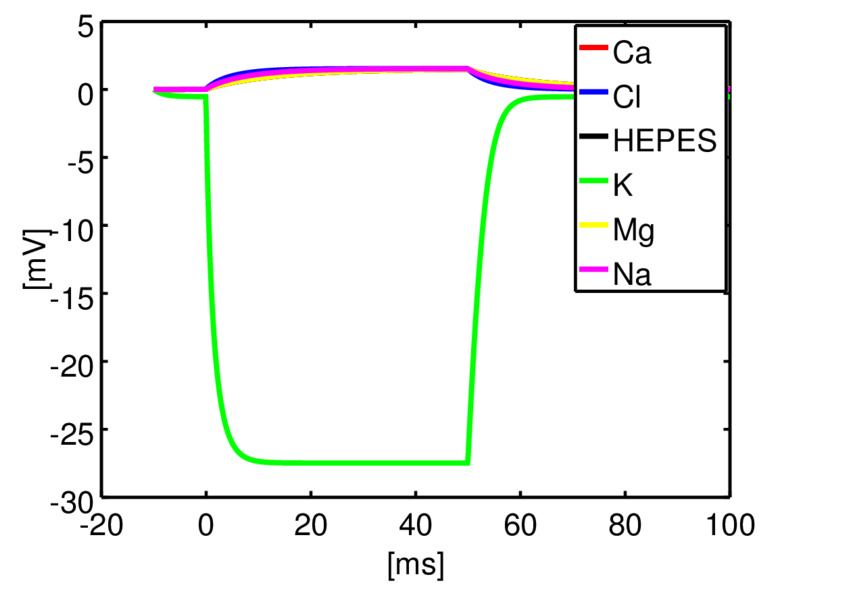}
\label{fig:Nernst_pots}} 
\caption
{(a): depolarizing
pulse of intracellular potential $V_{cell}$. (b): integral mean
of the cleft potential $\varphi(x,y)$ obtained with the electrical 
model. (c): integral mean
of the cleft potential $\varphi(x,y)$ obtained with the 
electrodiffusion model. (d): changes of the Nernst potentials
$V_{J0}$ between junction and bath.}
\label{fig:BF_models} 
\end{figure}

In this concluding section we compare the results obtained 
with the Area-Contact (A-C) model proposed in 
Sect.~\ref{sub:A-simplified-version} with the results 
of~\cite{Brittinger2005}, which we refer to for all physical data 
and details of the electrical equivalent circuits used to
determine the time evolution of ion concentrations in the cleft.
To this purpose, we conduct a first simulation considering 
given concentrations $c_i$ constant in time and only 
solving Eq.~\eqref{eq:BF_model_phi} (2D electrical model). 
Then, we conduct a second simulation 
accounting for ion dynamics, by adding 
an ODE system for the ionic concentrations (2D electrodiffusion model).
In both cases the integral mean $V_J$ is computed as in~\eqref{eq:VJ_def}.
The considered electrophysiological 
experiment is a voltage clamp stimulation with
the depolarizing pulse shown by Fig.~\ref{fig:pot_intracell} 
and the values of model parameters are the same as in~\cite{Brittinger2005}.
\begin{figure}[h!]
\centering 
\subfigure[$\varphi$ electrical model]
{\includegraphics[width=0.4\textwidth]{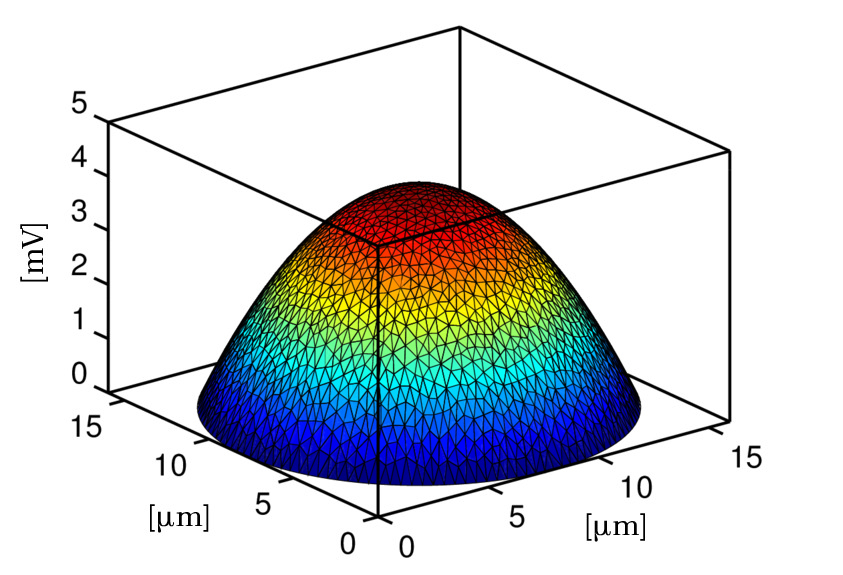}
\label{fig:pot_AC_10ms}} 
\quad
\subfigure[$\varphi$ electrodiffusion model]
{\includegraphics[width=0.4\textwidth]{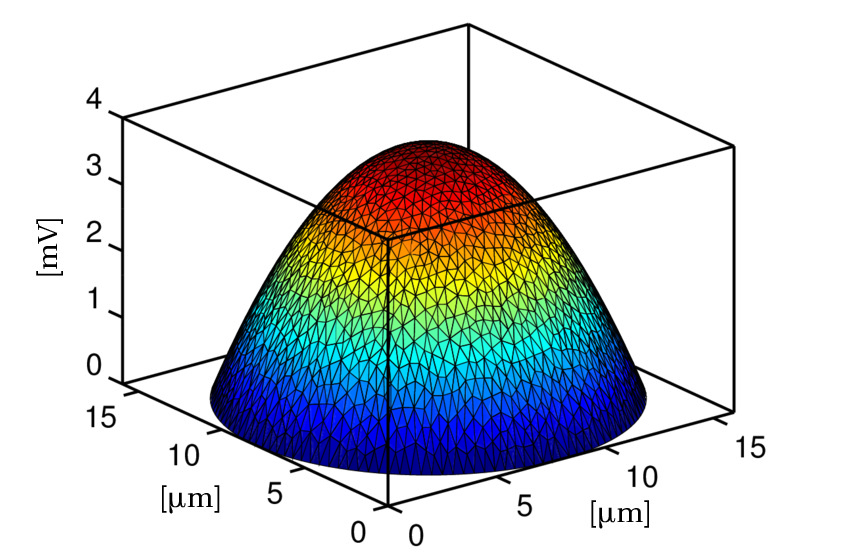}
\label{fig:pot_AC_035ms}} 
\caption
{\label{fig:BF_spatial_dis}Spatial distribution
of the cleft potential $\varphi$ in the circular domain at 
$t=0.35$ ms (cell just depolarized: on the left, electrical model, 
on the right, electrodiffusion model}.
\end{figure}

As demonstrated by Fig.~\ref{fig:potcleft_electricalmodel}, 
the 2D electrical model
accounts only for the fast response of the system,
with a dynamics determined by the electrical time constant 
$\tau=\left(C_{M}+C_{S}\right)/\sigma\simeq$ \SI{1.0944}{\micro\siemens}
($\sigma=q\sum_{i=1}^{M}|z_{i}|\mu_{i}c_{i}^{bath}\delta_{J}$ being 
the global cleft conductance): when the cell is depolarized, almost
instantaneously the potential goes to a value around \SI{2}{mV}.
The 2D electrodiffusion model, instead, describing
the time variation of $c_{i}$, accounts also for the slow component,
as shown by Figs.~\ref{fig:electrodiff_mean} and~\ref{fig:Nernst_pots}:
both potential
and concentrations have transients with a time constant in the
order of milliseconds, as expected (we have expressed the changes
of extracellular ion concentrations in the junction as Nernst
potentials between junction and bath). The integral mean of the electrical
potential $\varphi$ increases fast to a value around \SI{2.5}{mV}
and subsequently decays to a stationary level around \SI{1.5}{\milli\volt}
 and the potassium concentration increases from $5$ mM to $17$ mM, giving
a Nernst potential $V_{J}^{K}\simeq$ \SI{-27}{mV} in the junction.
Notably, all results are in excellent agreement with those reported by Brittinger and Fromherz in~\cite{Brittinger2005}.
Finally, the spatial distributions of the potential $\varphi$ computed
by the 2D electrical and electrodiffusion models, and 
used to determine $V_J$ in~\eqref{eq:VJ_def}, are reported  
in Fig.~\ref{fig:BF_spatial_dis}. The resulting parabolic shape
is in very good agreement with the behavior shown by the 
3D results of Sect.~\ref{sub:A-simplified-version}.

\section{Conclusions and Future Perspectives}\label{sec:conclusions}

In this article we have addressed the mathematical modeling
and numerical simulation of ion electrodiffusion 
in bio-hybrid devices. This subject is
of paramount importance in the wider scientific context of neuroelectronics,
where the main aim is to actually realize devices consisting of the
integration of biological tissues with solid-state integrated electronic
circuits.

In this treatise we have illustrated a suitable mathematical characterization
of bio-electronic interfaces, investigating different possible modeling
hypotheses on the coupling between the two different environments (cell
and electronic device) and on the derivation of model dimensional
reductions, performed to decrease the computational simulation effort.
A hierarchy of multiscale models has been therefore presented and
extensively validated with a broad range of numerical computations,
obtaining sensible results and comparing them with literature and
experiments. 

This mathematical description has also been applied to
complex configurations and has proved to be able to simulate the interactions
between multiple cells and multiple devices. Even if the present work
is not a faithful copy of a real-world biophysical setting, 
it can be considered a first step for the construction of 
mathematical models to be used in the design of actual devices.

Clearly, future research is needed to provide a better description
of the complex multiscale/multiphysics problem object of our
investigation. Among possible developments, we mention: 
\begin{itemize}
\item a more accurate modeling of the electronic
substrate, which can be useful in studying different types of stimulation,
for example a different polarization of the chip influencing the cell;
\item a model for the chemical binding mechanism of the ions to the electronic
substrate is also required, in order to fully describe the EOSFET
device; 
\item 
a coupling between electro-chemical and fluid-mechanical systems,
in order to account for the forces due to pressure differences and
flow in the aqueous medium;
\item a more realistic description of the problem geometry, with
full three-dimensional computations including the intracellular fluid
can be useful to faithfully reproduce the entire phenomena;
\item an application of the computational model to the 
simulation of electrophysiological experiments in 
current-clamp conditions.
\end{itemize}
The above mentioned improvements, particularly, the study of the
current-clamp protocol, should give the realistic chance to go further
in the study of the interactions between multiple cells, maybe introducing
a neural network and simulating a whole brain slice, as in the experimental
results of~\cite{hutzler2006high,zeck2001noninvasive}. 

\section*{Acknowledgments}
Matteo Porro and Riccardo Sacco were 
supported by Gruppo Nazionale per il Calcolo Scientifico
of Istituto Nazionale di Alta Matematica ``F.\ Severi''.
Thierry Nieus was partially supported by the SI-CODE project 
of the Future and Emerging Technologies (FET) programme within 
the Seventh Framework Programme for Research of The European 
Commission, under FET-Open grant number: FP7-284553.
\bibliographystyle{plain}
\bibliography{PaperBioHybrid.bib}

\end{document}